\documentclass[a4paper, 11pt]{article}

\usepackage{amsmath}
\numberwithin{equation}{section}
\usepackage{amsthm}
\usepackage{amssymb}
\usepackage{empheq}
\usepackage{mathtools}
\usepackage{physics}
\usepackage[dvipsnames]{xcolor}
\usepackage{tikz}
\usepackage{dsfont}
\usepackage[linkcolor=black, urlcolor=black, colorlinks=true, citecolor=blue]{hyperref}
\usepackage[title]{appendix}

\usepackage[
  a4paper,
  left=2.8cm,
  right=2.8cm,
  top=3cm,
  bottom=3cm
]{geometry}

\newcommand{\C}{\mathbb{C}}
\newcommand{\N}{\mathbb{N}}

\newcommand{\R}{\mathbb{R}}

\newcommand{\1}{\mathds{1}}

\newcommand{\cB}{\mathcal{B}}
\newcommand{\cD}{\mathcal{D}}
\newcommand{\cE}{\mathcal{E}}
\newcommand{\cF}{\mathcal{F}}

\newcommand{\cH}{\mathcal{H}}

\newcommand{\cU}{\mathcal{U}}
\newcommand{\cL}{\mathcal{L}}

\newcommand{\cP}{\mathcal{P}}

\newcommand{\cS}{\mathcal{S}}
\newcommand{\cT}{\mathcal{T}}
\newcommand{\cW}{\mathcal{W}}
\newcommand{\cX}{\mathcal{X}}

\newcommand{\mi}{\mathrm{i}}
\newcommand{\me}{\mathrm{e}}
\newcommand{\md}{\mathrm{d}}

\providecommand{\keywords}[1]{{\small \textbf{Keywords:} #1}}

\newtheorem{theorem}{Theorem}[section]
\newtheorem{corollary}[theorem]{Corollary}
\newtheorem{proposition}[theorem]{Proposition}
\newtheorem{lemma}[theorem]{Lemma}

\theoremstyle{definition}

\newtheorem{example}[theorem]{Example}
\newtheorem{remark}[theorem]{Remark}

\newtheorem{assumption}{Assumption}

\DeclareMathOperator{\dom}{Dom}
\DeclareMathOperator{\ran}{Ran}
\DeclareMathOperator{\gap}{gap}
\DeclareMathOperator{\myspan}{span}

\usepackage[backend=bibtex, style=alphabetic]{biblatex}
\usepackage{authblk}

\title{Spectral Analysis for Gaussian Quantum Markov Semigroups}
\author[1]{Franco Fagnola}
\author[2]{Zheng Li}

\affil[1]{Dipartimento di Matematica, Politecnico di Milano, Piazza Leonardo da Vinci 32, Milano, 20133, Italy.}
\affil[2]{School of Mathematics and Statistics, Central South University,\protect\linebreak
Changsha, 410083, China.}

\affil[ ]{\href{mailto:franco.fagnola@polimi.it}{franco.fagnola@polimi.it}, \href{mailto:zheng.li@csu.edu.cn}{zheng.li@csu.edu.cn}}

\date{}

\usepackage{changes}
\definechangesauthor[name={Zheng}, color=red]{ZL}

\allowdisplaybreaks

\begin{document}

    \maketitle

    \begin{abstract}

        Let $(\cT_t)_{t\geq0}$ be a Gaussian quantum Markov semigroup with a faithful normal invariant state $\rho$. For every $s\in[0,1]$, the $s$-embedding associated with $\rho$ induces a contraction semigroup $(T_t^{(s)})_{t\geq0}$ on the Hilbert--Schmidt space $\cB_2(\mathsf h)$; let $L^{(s)}$ denote its generator. Without assuming symmetry or quantum detailed balance, we determine the full spectra of $L^{(s)}$ and $L^{(s)*}$: they consist, respectively, of the non-negative integer combinations of the eigenvalues of the phase-space drift matrix and of their complex conjugates. We also diagonalize the self-adjoint closure of $L^{(s)*}+L^{(s)}$ and obtain an explicit formula for the spectral gap. Using quantum characteristic functions, we represent $T_t^{(s)}$, up to unitary equivalence, as a complex Gaussian integral operator and prove that its kernel is square-integrable for every $t>0$. Hence $T_t^{(s)}$ is a Hilbert--Schmidt operator on $\cB_2(\mathsf h)$ for every $t>0$, and every induced semigroup is immediately compact. Consequently, $L^{(s)}$ and $L^{(s)*}$ have compact resolvent for all $s\in[0,1]$, and their point spectra exhaust their full spectra. These results provide a full quantum counterpart of classical Ornstein--Uhlenbeck spectral theory.

    \end{abstract}


    \vskip1cm

    {\small
        \tableofcontents
    }

    \section{Introduction} \label{section-introduction}

        Let $\mathsf h$ be a bosonic Fock space, let $\cB(\mathsf h)$ be the von Neumann algebra of bounded operators on $\mathsf h$, and let $\cB_1(\mathsf h)$ denote the trace-class operators. Throughout, $\N:=\{0,1,2,\ldots\}$. A quantum Markov semigroup (QMS) $(\cT_t)_{t\geq0}$ on $\cB(\mathsf h)$ is a weak$^*$-continuous semigroup of normal, completely positive, identity-preserving maps. It acts on observables in the Heisenberg picture. The predual semigroup $(\cT_{*t})_{t\geq0}$ acts on trace-class operators and preserves density matrices in the Schr{\"o}dinger picture. It is characterized by
    \begin{equation} \label{eq-introduction-predual-semigroup}
        \Tr\!\left(\cT_{*t}(\psi)X\right)
        =\Tr\!\left(\psi\,\cT_t(X)\right),
        \quad \psi\in\cB_1(\mathsf h),\quad X\in\cB(\mathsf h).
    \end{equation}
    QMSs provide an operator-algebraic framework for Markovian open-system dynamics. Their generators describe dissipation and relaxation; in infinite dimension, however, these generators are generally unbounded, so their domains are an essential part of the analysis \cite{Da77,ChFa98,SiHoWe17}.

    Gaussian QMSs---also called quasi-free QMSs---form a distinguished family of bosonic QMSs \cite{DEMOEN197927,Po22,BaWe24}. In the one-mode setting considered here, $\mathsf h=\Gamma(\C)$, and the annihilation and creation operators $a$ and $a^\dagger$ satisfy $[a,a^\dagger]=\1$. For $z\in\C$, the associated Weyl operator is
        \[
            W(z):=\exp\!\left\{z a^\dagger-\overline z a\right\}.
        \]
    A normal state $\rho$ is called Gaussian if its quantum characteristic function $\chi_\rho $, defined in \eqref{eq-introduction-characteristic-function}, is a Gaussian function, and a QMS is called Gaussian if its predual semigroup preserves Gaussian states. Note that every normal state on $ \cB(\mathsf{h}) $ is uniquely represented by a density matrix. Throughout the paper, we identify a normal state with its density matrix and use the same symbol for both. By \cite{agredo2021gaussian}, such a Gaussian QMS $ (\cT_t)_{t \ge 0} $ acts on Weyl operators according to
    \begin{equation} \label{eq-introduction-gaussian-weyl-action}
        \cT_t(W(z)) = \exp\!\left\{-\frac{1}{2}\int_0^t
        \Re\!\left(\overline{\me^{rZ}z}\,C(\me^{rZ}z)\right)\md r + \mi\int_0^t
        \Re\!\left(\overline\zeta\,\me^{rZ}z\right)\md r\right\}
        W(\me^{tZ}z).
    \end{equation}
    Here $ Z $ and $ C $ are real-linear operators on $ \C $ and $\zeta\in\C$. We call $Z$ the \emph{drift operator}, $C$ the \emph{diffusion operator}, and $\zeta$ the displacement parameter. Their expressions in terms of the coefficients of the Gorini--Kossakowski--Sudarshan--Lindblad (GKLS) representation are given in \eqref{eq-introduction-drift-operator} and \eqref{eq-introduction-diffusion-operator}. Appendix~\ref{section-appendix-identification} fixes the real phase-space identification and its matrix convention in \eqref{eq-identification-real-linear-operators}. The real $2\times2$ matrices representing $Z$ and $C$, as complex linear operators on $\mathbb{C}^2$, are denoted by $\mathbf Z$ and $\mathbf C$ and are called, respectively, the \emph{drift matrix} and the \emph{diffusion matrix}. Thus the infinite-dimensional operator evolution is controlled by finite-dimensional drift, diffusion, and displacement data. Gaussian QMSs are therefore the open-system counterparts of the quadratic Hamiltonians that generate Bogoliubov dynamics in closed bosonic systems \cite{bruneau2007bogoliubov,girotti2024gaussian}. Their analysis is more involved because, in the generalized GKSL representation \eqref{eq-gaussian-gksl-generator}, the Hamiltonian $ H $ in \eqref{eq-Hamiltonian-with-parameters} and the noise operator $ L_\ell $ in \eqref{eq-kraus-operator-a-adagger} give rise to the two generally noncommuting operators $ - \mi H $ and $ - \frac{1}{2} \sum_\ell L^*_\ell L_\ell $, which form the anti-self-adjoint and self-adjoint contributions to the operator underlying the generalized GKSL construction \cite[Section 2]{agredo2021gaussian}. The noise operators $ L_\ell $ also enter through the completely positive part of the generator.

    This finite-dimensional description is important both mathematically and physically. Gaussian states are determined by their first and second moments and include coherent, squeezed, and thermal states. They provide tractable models for quantum optics and quantum systems with continuous variables; see \cite{parthasarathy2010gaussian,holevo2011probabilistic}. At the same time, the Gaussian class retains genuinely quantum features: the canonical variables do not commute, phase space carries a symplectic structure, and the natural invariant state is usually nontracial. These features make the choice of a Hilbert-space realization nontrivial.

    Spectral information describes several dynamical scales at once. The real parts of generator eigenvalues give decay rates, whereas their imaginary parts give oscillation frequencies. A strictly positive spectral gap, defined for the induced generators below in \eqref{eq-definition-spectral-gap}, separates the stationary vector from the decaying sector and yields an exponential rate of convergence to equilibrium \cite{frigerio1977quantum,BERTINI2024110340,ChKiKo21}. Eigenvalues clustered near the origin identify metastable time scales and the corresponding long-lived manifold \cite{macieszczak2016metastability}. For quadratic bosonic Lindbladians, spectral data also underlie normal-master-mode descriptions of the evolution \cite{prosen2010quantization,Te19}.

    The spectrum has further operator-theoretic uses. Once the semigroup is immediately compact, the spectral mapping theorem recovers the nonzero spectrum of each semigroup operator from the generator spectrum and yields a discrete modal expansion \cite[Corollary IV.3.12, Corollary V.3.2]{engel2000one}. In a non-self-adjoint problem, left and right eigenvectors provide the natural coordinates for perturbation theory \cite{li2014perturbative}. These applications require more than a formal list of polynomial eigenvalues: one must specify a closed realization of the generator, control its domain, and determine whether any continuous or residual spectrum remains.

    The classical benchmark is the Ornstein--Uhlenbeck semigroup on $L^2$ of an invariant Gaussian measure. Under stability and nondegeneracy assumptions, the spectrum of its generator consists of the non-negative integer combinations of the eigenvalues of the drift matrix \cite{metafune2002spectrum}; the corresponding Fokker--Planck generator has an analogous spectral structure \cite{arnold2018rates}. Complete $L^2$ spectra are also known for several symmetric quantum Ornstein--Uhlenbeck models. Cipriani, Fagnola, and Lindsay computed the spectrum and eigenspaces of a symmetric semigroup on the CCR algebra \cite{cipriani2000spectral}. Ko and Park obtained discrete spectra for a family of symmetric generators by unitary equivalence with number operators \cite{ko2004construction}. Bertini, De Sole, and Posta used the spectral resolution of a symmetric Bose Ornstein--Uhlenbeck generator to study trace-distance ergodicity \cite{BERTINI2024110340}. More recently, symmetric Gaussian QMSs have been characterized and related to direct sums of Ornstein--Uhlenbeck-type components \cite{poletti2026symmetric}.

    \enlargethispage{\baselineskip}

    For Gaussian QMSs without a symmetry or detailed-balance assumption, the optimal $L^2$ convergence rate can be computed explicitly from the drift and diffusion matrices \cite{fagnola2025spectral}. That result also shows how the rate depends on the chosen noncommutative $L^2$ embedding, in particular on whether one uses the GNS or KMS embedding. The distinction between these realizations is also important from a structural viewpoint: under the corresponding symmetry assumptions, different representation and derivation techniques become available \cite{VW23,Wir24}.

    Let $ (\cT_t)_{t \ge 0} $ be a Gaussian QMS admitting a faithful normal invariant state $ \rho $. Thus $\cT_{*t}(\rho)=\rho$ for every $t\geq0$ and $\ker\rho=\{0\}$. Let $\cB_2(\mathsf h)$ be the Hilbert space of Hilbert--Schmidt operators, with inner product $\langle X,Y\rangle_2:=\Tr(X^*Y)$. For $s\in[0,1]$, the map
    \[
        \cB(\mathsf h)\ni X\longmapsto
        \rho^{s/2}X\rho^{(1-s)/2}\in\cB_2(\mathsf h)
    \]
    is called the \emph{$s$-embedding}; its range is dense in $\cB_2(\mathsf h)$. On this range, define
    \begin{equation} \label{eq-induced-semigroup-and-qms}
        T_t^{(s)}\!\left(\rho^{s/2}X\rho^{(1-s)/2}\right)
        :=\rho^{s/2}\cT_t(X)\rho^{(1-s)/2}.
    \end{equation}
    By \cite[Theorem 2.3]{carbone2000exponential}, this prescription extends uniquely to a strongly continuous contraction semigroup $(T_t^{(s)})_{t\geq0}$ on $\cB_2(\mathsf h)$. We call it the \emph{induced semigroup} and denote its generator by $L^{(s)}$. The choices $s=0$, $s=1/2$, and $s=1$ are the Gelfand--Naimark--Segal (GNS), Kubo--Martin--Schwinger (KMS), and anti-GNS embeddings, respectively; the corresponding Hilbert-space constructions originate in the GNS approach of Frigerio \cite{frigerio1977quantum} and the standard-form approach of Cipriani \cite{cipriani1997dirichlet}. The other values of $s$ interpolate between these distinguished realizations. 

    The unitary standardization described in Appendix \ref{section-non-diagonal-gaussian-states} allows us to replace the original Gaussian QMS by a unitarily equivalent one whose invariant state is the thermal density operator
    \begin{equation} \label{eq-introduction-thermal-state}
        \rho=(1-\me^{-\beta})\me^{-\beta N},
        \quad N:=a^\dagger a,\quad \beta>0.
    \end{equation}
    Here $N$ is the number operator and $\beta$ is the inverse temperature. Section~\ref{section-preliminaries} gives the precise GKSL generator and states the assumptions under which this reduction is used.

    \paragraph{Main results and proof strategy.}
    Let $\lambda$ and $\mu$ be the two eigenvalues of the drift matrix $\mathbf Z$. They are eigenvalues of $L^{(s)}$ for every $s\in[0,1]$. We call them the \emph{fundamental eigenvalues} because they occur on embedded polynomials of degree one in $a$ and $a^\dagger$ and generate all other spectral values through non-negative integer combinations. Writing $\sigma(A)$ for the spectrum of a closed operator $A$, we prove
    \begin{equation} \label{eq-set-of-eigenvalues}
        \sigma\!\left(L^{(s)}\right)
        =\left\{n\lambda+m\mu:n,m\in\N\right\},
        \quad
        \sigma\!\left(L^{(s)*}\right)
        =\left\{n\overline\lambda+m\overline\mu:n,m\in\N\right\}.
    \end{equation}
    A matrix is called \emph{defective} if it is not diagonalizable. If $\mathbf Z$ is defective, then $\lambda=\mu$, and the set of distinct spectral values in the first formula reduces to $\{n\lambda:n\in\N\}$; the generalized eigenvectors retain the Jordan information. Formula \eqref{eq-set-of-eigenvalues} shows that the spectra are independent of $s$, even though the eigenvectors and Hilbert-space geometry depend on the embedding.

    For each $r\in\N$, the algebraic part of the proof is based on the finite-dimensional space
    \begin{equation} \label{eq-introduction-polynomial-spaces}
        \begin{aligned}
        \cX_r:=\bigl\{\rho^{s/2}P(a,a^\dagger)\rho^{(1-s)/2}:{}&
        P\text{ is a noncommutative polynomial in }a\text{ and }a^\dagger\\[-2pt]
        &\text{of total degree at most }r\bigr\}.
        \end{aligned}
    \end{equation}

    We set $\cX:=\bigcup_{r\in\N}\cX_r$. The spaces $\cX_r$ form a \textit{filtration} $\cX_0\subset\cX_1\subset\cdots$, and Lemma~\ref{lemma-density-of-domain-L-star-plus-L} shows that the embedded-polynomial space $\cX$ is dense in $\cB_2(\mathsf h)$. The degree estimate of Lemma~\ref{lemma-cL-X-still-polynomial}, transferred to the induced generator by Theorem~\ref{theorem-L-and-calL-relation}, gives the invariance
    \[
        L^{(s)}(\cX_r)\subset\cX_r
    \]
    recorded in \eqref{eq-invariance-X-r-under-induced-generator}. The quasi-derivation identity in Theorem~\ref{theorem-derivation-property-qms} plays a different role: it is used in Theorem~\ref{theorem-other-eigenvalues-L} to construct higher-degree eigenvalues from the fundamental ones, and its correction term lowers total degree by two. Consequently, relative to a basis ordered by nondecreasing total degree, $L^{(s)}\restriction\cX_r$ is upper triangular, with diagonal entries given by the non-negative integer combinations in \eqref{eq-set-of-eigenvalues}; see Theorem~\ref{theorem-other-eigenvalues-L}. The orthogonal projections $ \cP_r $ onto the filtration spaces $ \cX_r $ commute with both $ L^{(s)} $ and $ L^{(s)*} $ on their respective domains by Lemma \ref{lemma-exchangeability-induced-generator-P}. Together with the duality result in Corollary \ref{theorem-adjoint-generator-eigenvalues}, this yields the exhaustion of the point spectra in Proposition \ref{proposition-exhaustion-point-spectrum}.

    We also analyze the symmetric operator $L^{(s)*}+L^{(s)}$ on $\cX$, which is a common core by Lemma~\ref{lemma-cX-common-core}. Lemma~\ref{lemma-L-star-plus-L-essential-sa} proves that this operator is essentially self-adjoint on $\cX$, while Lemma~\ref{lemma-fundamental-eigenvalues-of-Ls*-plus-Ls} identifies its two fundamental eigenvalues, expressed in terms of the GKLS parameters $\gamma$ and $\kappa$ from \eqref{eq-introduction-drift-operator} and \eqref{eq-Hamiltonian-with-parameters}:
    \begin{equation} \label{eq-introduction-symmetric-fundamental-eigenvalues}
        \nu_s^{\pm}:=-2\gamma\pm2\abs{\kappa}
        \cosh\!\left(\frac{1-2s}{2}\,\beta\right)
    \end{equation}
    The spectrum of the self-adjoint closure consists of the values $n\nu_s^++m\nu_s^-$ with $n,m\in\N$. An operator $L$ is called \emph{maximally dissipative} if $\Re\langle x,Lx\rangle\leq0$ for every $x\in\dom L$ and $L$ has no proper dissipative extension. For any maximally dissipative operator with nontrivial kernel, we define its spectral gap by
    \begin{equation} \label{eq-definition-spectral-gap}
        \gap L:=\inf\left\{-\Re\langle x,Lx\rangle:
        x\in\dom L,\ \norm{x}=1,\ x\perp\ker L\right\}.
    \end{equation}
    Theorem~\ref{theorem-spectral-gap-compact-resolvent} and Remark~\ref{remark-exact-spectral-gap} give the exact formula
    \begin{equation} \label{eq-introduction-gap-formula}
        \gap L^{(s)}
        =\gamma-\abs{\kappa}\cosh\!\left(\frac{1-2s}{2}\,\beta\right).
    \end{equation}
    By the inequality $\gamma\geq\abs{\kappa}\cosh(\beta/2)$ proved in Lemma~\ref{lemma-elementary-inequalities-spectral-gap}, the spectral gap is maximal at the KMS embedding $s=1/2$, is strictly positive for $0<s<1$, and may vanish at the GNS and anti-GNS embeddings $s=0, 1$. Whenever the spectral gap is strictly positive, Theorem~\ref{theorem-spectral-gap-compact-resolvent} shows that the self-adjoint closure of $L^{(s)*}+L^{(s)}$ has compact resolvent. The quadratic-form comparison in Proposition~\ref{theorem-compact-resolvent-spectral-gap} then gives the same property for $L^{(s)}$ and $L^{(s)*}$.

    The preceding argument does not cover the endpoints of $[0,1]$ at which the spectral gap vanishes. To determine the full spectra for every $s\in[0,1]$, we use quantum characteristic functions. For $\psi\in\cB_1(\mathsf h)$, define
    \begin{equation} \label{eq-introduction-characteristic-function}
        \chi_\psi(z):=\Tr(\psi W(z)),
        \quad z\in\C.
    \end{equation}
    The quantum Plancherel theorem \cite[Section~5.3]{holevo2011probabilistic} extends the normalized map $\psi\mapsto\pi^{-1/2}\chi_\psi$ to a unitary operator $\cF_W:\cB_2(\mathsf h)\to L^2(\C)$, defined in \eqref{eq-unitary-characteristic-function-transform}. Corollary~\ref{corollary-characteristic-function-induced-semigroup} proves that $\cF_WT_t^{(s)}\cF_W^{-1}$ is the complex Gaussian integral operator displayed in \eqref{eq-characteristic-function-evolution-of-semigroup-Y}. The complex symmetric matrix $\mathbf G_t$ in its kernel is defined in \eqref{eq-definition-G-t}, and Lemma~\ref{lemma-strict-positivity-real-parts-G-t} gives
    \[
        \Re\mathbf G_t>0, \quad t>0, \quad s\in[0,1].
    \]
    Consequently, the Gaussian kernel is square-integrable, and $T_t^{(s)}$ is a Hilbert--Schmidt operator on $\cB_2(\mathsf h)$ for every $t>0$. Independently, Theorem~\ref{theorem-sobolev-regularization} shows that $T_t^{(s)}$ maps $\cB_2(\mathsf h)$ continuously into a first-order bosonic Sobolev space that is compactly embedded in $\cB_2(\mathsf h)$.

    It follows from Theorem~\ref{theorem-sobolev-regularization} that every induced semigroup is \emph{immediately compact}, meaning that $T_t^{(s)}$ is compact for each $t>0$. Corollary~\ref{corollary-compact-resolvent-full-spectrum} then proves that $L^{(s)}$ and $L^{(s)*}$ have compact resolvent for every $s\in[0,1]$ and that the eigenvalue sets in \eqref{eq-set-of-eigenvalues} are their full spectra, including at the endpoint embeddings where the spectral gap may vanish. The polynomial filtration identifies the spectral values, while the characteristic-function representation proves that no continuous or residual spectrum remains.

    We choose the one-mode setting to keep the notation transparent and the main spectral and compactness mechanisms clearly visible. Several components of the argument---including the polynomial filtration, the unitary standardization of the invariant state, and the characteristic-function calculations---extend naturally to finitely many modes and suggest a route to a finite-mode spectral theorem. Carrying out that program, however, requires more delicate computations and more advanced matrix-analysis methods, particularly in the positivity analysis of the complex Gaussian kernel.

    The paper is organized as follows. Section~\ref{section-preliminaries} presents the one-mode Gaussian framework. Section~\ref{section-eigenvalues-of-L} states the standing assumptions and uses the spaces in \eqref{eq-introduction-polynomial-spaces} to determine the point spectra of $L^{(s)}$ and $L^{(s)*}$. Section~\ref{section-spectrum-L-star-plus-L} studies $L^{(s)*}+L^{(s)}$, computes the spectral gap, and proves the compact-resolvent implication when the gap is positive. Section~\ref{section-immediate-compactness} derives the characteristic-function representation, proves immediate compactness, and completes the spectral determination. Appendix~A contains the domain arguments for embedded polynomials, Appendix~B collects phase-space and polynomial identities, Appendix~C performs the diagonal standardization, and Appendix~D contains the Gaussian integrations and matrix-positivity proofs.

    \section{The One-Mode Gaussian Framework} \label{section-preliminaries}

    We work on the one-mode bosonic Fock space $\mathsf h:=\Gamma(\C)$, with canonical orthonormal basis $(e_n)_{n\geq0}$ and inner product $\langle\cdot,\cdot\rangle$. Two dense subspaces of $\mathsf h$ will be used repeatedly. The finite-particle space is $\cD:=\myspan\{e_n:n\in\N\}$, and $\cE$ denotes the linear span of the exponential vectors. For $f\in\C$, the exponential vector $e(f)\in\mathsf h$ is
    \begin{equation*}
        e(f) := \sum_{k \geq 0} \frac{f^k}{\sqrt{k!}} e_k.
    \end{equation*}
    We write $\cB(\mathsf h)$, $\cB_1(\mathsf h)$, and $\cB_2(\mathsf h)$ for the bounded, trace-class, and Hilbert--Schmidt operators on $\mathsf h$, respectively, and use $\langle\cdot,\cdot\rangle_2$ for the Hilbert--Schmidt inner product. The annihilation operator $a$ and creation operator $a^\dagger$ act on the canonical basis as follows:
    \begin{equation*}
        a e_n = \sqrt{n} \, e_{n-1}, \quad n \geq 1, \quad a e_0 = 0; \quad a^\dagger e_n = \sqrt{n+1} \, e_{n+1}, \quad n \geq 0.
    \end{equation*}
    We retain the notation $N=a^\dagger a$ for the number operator from \eqref{eq-introduction-thermal-state}. We also use the momentum and position operators $p:=\mi(a^\dagger-a)/\sqrt2$ and $q:=(a^\dagger+a)/\sqrt2$, respectively. For $z\in\C$, $W(z)$ denotes the unitary Weyl operator introduced in Section~\ref{section-introduction}.

    We use the notions of Gaussian states and Gaussian QMS introduced in Section~\ref{section-introduction}. More explicitly, in one mode, a normal state $\rho$ is Gaussian if there exist $\omega\in\C$ and a real-linear, symmetric, invertible operator $S$ such that
    \begin{equation*}
        \Tr(\rho W(z))
        =\exp\!\left\{2\mi\Im\langle\omega,z\rangle
        -\frac12\Re\langle z,Sz\rangle\right\},
        \qquad z\in\C.
    \end{equation*}
    This formula displays the Gaussian form of the quantum characteristic function mentioned before. Theorems~\ref{theo-existence-inv-state-gqms} and~\ref{theo-gaussian-qms-explicit-weyl-operators} are taken from \cite{agredo2021gaussian}.

    Let $(\cT_t)_{t\geq0}$ be a Gaussian QMS on $\cB(\mathsf h)$. Its generator has the Gaussian GKSL form
    \begin{equation} \label{eq-gaussian-gksl-generator}
        \cL(X)=\mi[H,X]-\frac12\sum_{\ell=1}^2
        \left(L_\ell^*L_\ell X-2L_\ell^*XL_\ell+XL_\ell^*L_\ell\right),
    \end{equation}
    where
        \begin{align}
            H&=\Omega a^\dagger a+\frac{\kappa}{2}a^{\dagger2}
            +\frac{\overline\kappa}{2}a^2+\frac{\zeta}{2}a^\dagger
            +\frac{\overline\zeta}{2}a, \label{eq-Hamiltonian-with-parameters}\\
            L_\ell&=\overline v_\ell a+u_\ell a^\dagger,
            \qquad u_\ell,v_\ell\in\C,\quad \ell=1,2. \label{eq-kraus-operator-a-adagger}
        \end{align}
        The drift and diffusion operators in \eqref{eq-introduction-gaussian-weyl-action} are the real-linear operators
        \begin{align}
            Zz&=(\mi\Omega-\gamma)z+\mi\kappa\overline z,
            &&\text{where }\quad
            \gamma:=\frac12\sum_{\ell=1}^2
            \left(\abs{v_\ell}^2-\abs{u_\ell}^2\right),
            \label{eq-introduction-drift-operator}\\
            Cz&=\sum_{\ell=1}^2\left(
            (\abs{u_\ell}^2+\abs{v_\ell}^2)z
            +2v_\ell u_\ell\overline z\right).
            \label{eq-introduction-diffusion-operator}
        \end{align}
        The generator $\cL_*$ of the predual semigroup defined in \eqref{eq-introduction-predual-semigroup} is given by
    \begin{equation} \label{eq-gksl-gaussian-predual-generator}
        \cL_* (\psi) = - \mi [H, \psi] - \frac{1}{2} \sum_{\ell=1}^2 \left( L_\ell^* L_\ell \psi - 2 L_\ell \psi L_\ell^* + \psi L_\ell^* L_\ell \right), \quad \psi \in \dom \cL_*,
    \end{equation}
    where $H$ and $L_\ell$ are given by \eqref{eq-Hamiltonian-with-parameters} and \eqref{eq-kraus-operator-a-adagger}.

    The next theorem gives the precise parameter conditions under which the Gaussian QMS has a normal invariant state.

    \begin{theorem} \label{theo-existence-inv-state-gqms}
        A normal invariant state $ \rho $ for $ (\cT_t)_{t \ge 0} $ exists if and only if
        \begin{equation} \label{eq-cond-existence-inv-state}
            \gamma > 0, \quad \gamma^2 + \Omega^2 - \abs{\kappa}^2 > 0,
        \end{equation}
        where $\Omega$ and $\kappa$ are as in \eqref{eq-Hamiltonian-with-parameters} and $\gamma$ is defined in \eqref{eq-introduction-drift-operator}. Moreover, $\rho$ is Gaussian and is the unique normal invariant state of $(\cT_t)_{t\geq0}$.
    \end{theorem}

    The next theorem recalls the Weyl action already displayed in the introduction.

    \begin{theorem} \label{theo-gaussian-qms-explicit-weyl-operators}
        The action of $(\cT_t)_{t\geq0}$ on Weyl operators is given by \eqref{eq-introduction-gaussian-weyl-action}, where the drift operator $Z$ and the diffusion operator $C$ are defined in \eqref{eq-introduction-drift-operator} and \eqref{eq-introduction-diffusion-operator}, respectively.
    \end{theorem}

    We next record two algebraic properties of the Gaussian generator. Here $\deg$ denotes total polynomial degree in $a$ and $a^\dagger$. The first property follows from the Wick-ordering result in Theorem~\ref{theo-Wicks-theorem} and the commutation relations; the second is a quasi-derivation identity. We call \eqref{eq-iteration-formula} a quasi-derivation identity because it differs from the usual Leibniz rule by a correction term of degree at most $\deg(X)+\deg(Y)-2$.

    \begin{lemma} \label{lemma-cL-X-still-polynomial}
        Let $X$ be a polynomial of degree $n$ in $a^\dagger$ and degree $m$ in $a$. Then $\cL(X)$, interpreted algebraically through \eqref{eq-gaussian-gksl-generator}, either vanishes or is a polynomial in $a$ and $a^\dagger$ satisfying
        \begin{equation*}
            \deg (\cL (X)) \le m + n.
        \end{equation*}
    \end{lemma}
    \begin{proof}
        By Theorem \ref{theo-Wicks-theorem}, it suffices to prove the statement for the Wick normal-ordered product $ a^{\dagger n} a^m $, which has degree $ n $ in $ a^\dagger $ and degree $ m $ in $ a $.

        We first note that the generator (\ref{eq-gaussian-gksl-generator}) can be equivalently expressed as
        \begin{equation} \label{eq-gksl-equivalent-form}
            \cL (X) = \mi [H, X] - \frac{1}{2} \sum_\ell \left( [X, L_\ell^*] L_\ell - L_\ell^* [X, L_\ell] \right).
        \end{equation}
        It remains to compute $\cL(a^{\dagger n}a^m)$. For example, if $X=a^{\dagger n}a^m$, Lemma~\ref{lemma-commutation-relations-a-a-dagger-powers} gives
        \begin{align}
            [a^{\dagger n} a^m, L_\ell^*] &= [a^{\dagger n} a^m, v_\ell a^\dagger + \overline{u_\ell} a] \nonumber \\
            &= a^{\dagger n} [a^m, v_\ell a^\dagger + \overline{u_\ell} a] + [a^{\dagger n}, v_\ell a^\dagger + \overline{u_\ell} a] a^m \nonumber \\
            &= m v_\ell a^{\dagger n} a^{m-1} - n \overline{u_\ell} (a^\dagger)^{n-1} a^m, \label{eq-computations-CR-X-L-l}
        \end{align}
        so $ [X, L_\ell^*] L_\ell $ is a polynomial in $ a $ and $ a^\dagger $ and
        \begin{equation*}
            \deg \left( [X, L_\ell^*] L_\ell \right) \le m + n.
        \end{equation*}

        Similar computations apply to $L_\ell^*[X,L_\ell]$. For the Hamiltonian part, in addition to the quadratic terms, we use
        \begin{equation*}
            [a^\dagger a,X]=a^\dagger[a,X]+[a^\dagger,X]a,
        \end{equation*}
        which shows explicitly that the number-Hamiltonian term $\Omega a^\dagger a$ does not increase total degree. The same commutator expansion applies to $a^{\dagger2}$ and $a^2$, completing the proof.

    \end{proof}

    \begin{theorem} \label{theorem-derivation-property-qms}
        Let $ X, Y $ be polynomials in $ a $ and $ a^\dagger $. Then,
        \begin{equation} \label{eq-iteration-formula}
            \cL (XY) = X \cL (Y) + \cL (X) Y + \sum_\ell [X, L_\ell^*] [L_\ell, Y],
        \end{equation}
        where $ \sum_\ell [X, L_\ell^*] [L_\ell, Y] $ is a polynomial in $ a $ and $ a^\dagger $. If $ X $ and $ Y $ are non-zero, this polynomial either vanishes or satisfies
        \begin{equation} \label{eq-degrees-comparision}
            \deg \left( \sum_\ell [X, L_\ell^*] [L_\ell, Y] \right) \le \deg (X) + \deg (Y) - 2.
        \end{equation}
    \end{theorem}
    \begin{proof}
        Using \eqref{eq-gksl-equivalent-form} we compute:
        \begin{align*}
            \cL (X Y) &= \mi H X Y - \mi X Y H \\
            &\quad - \frac{1}{2} \sum_\ell \left( X[Y, L_\ell^*] L_\ell + [X, L_\ell^*] Y L_\ell - L_\ell^* X [Y, L_\ell] - L_\ell^* [X, L_\ell] Y \right) \\
            &= \mi X [H, Y] + \mi [H, X] Y \\
            &\quad - \frac{1}{2} \sum_\ell \left( X[Y, L_\ell^*] L_\ell - X L^*_\ell [Y, L_\ell] + [X, L_\ell^*] L_\ell Y - L^*_\ell [X, L_\ell] Y \right) \\
            &\quad + \sum_\ell \left( X L_\ell^* L_\ell Y - X L_\ell^* Y L_\ell - L_\ell^* X L_\ell Y + L_\ell^* X Y L_\ell \right) \\
            &= X \cL (Y) + \cL (X) Y + \sum_\ell [X, L_\ell^*] [L_\ell, Y],
        \end{align*}
        which completes the proof of \eqref{eq-iteration-formula}.

        By Theorem~\ref{theo-Wicks-theorem} and \eqref{eq-computations-CR-X-L-l}, if $X$ is not scalar, every commutator $[X,L_\ell^*]$ either vanishes or has degree at most $\deg(X)-1$; if $X$ is scalar, all these commutators vanish. The analogous statement holds for $[L_\ell,Y]$. Consequently, $\sum_\ell[X,L_\ell^*][L_\ell,Y]$ vanishes whenever $X$ or $Y$ is scalar. Otherwise, every nonzero summand satisfies
        \begin{equation*}
            \deg \left( \sum_\ell [X, L_\ell^*] [L_\ell, Y] \right) \le \deg (X) + \deg (Y) - 2.
        \end{equation*}
        The sum therefore either vanishes or has degree at most $\deg(X)+\deg(Y)-2$, which proves \eqref{eq-degrees-comparision}.
    \end{proof}

    \section{Polynomial Spectra of the Induced Generators} \label{section-eigenvalues-of-L}

    We impose the following assumptions throughout the remainder of the paper.

    \begin{assumption}\label{hyp-stability-Z}
        The drift matrix $\mathbf Z$ representing the operator $Z$ in \eqref{eq-introduction-drift-operator} is stable.
    \end{assumption}

    A matrix is called \textit{stable} if all its eigenvalues have strictly negative real parts. The matrix $\mathbf Z$ is stable if and only if \eqref{eq-cond-existence-inv-state} holds. Theorem~\ref{theo-existence-inv-state-gqms} then yields a unique normal invariant state $\rho$, which is Gaussian.

    Under the stability assumption, the unique normal invariant Gaussian state is either faithful or pure. The exceptional parameter regime in which it is pure is characterized in \cite[Proposition~5]{agredo2021gaussian}; Assumption~\ref{hyp-faithfulness} excludes precisely this regime. After the unitary standardization below, faithfulness is equivalent to a finite inverse temperature $\beta > 0 $.

    \begin{assumption}\label{hyp-faithfulness}
        The normal invariant state $ \rho $ is also faithful.
    \end{assumption}

    Appendix~\ref{section-non-diagonal-gaussian-states} shows that a unitary standardization reduces every QMS satisfying Assumptions~\ref{hyp-stability-Z}--\ref{hyp-faithfulness} to one whose invariant density operator is diagonal in the number basis. This unitary equivalence preserves all spectra, multiplicities, compactness properties, and numerical spectral gaps considered here. We therefore assume from now on that
    \begin{equation} \label{eq-diagonal-rho-parameters}
        \zeta = 0, \quad \kappa = \mi \tanh \left( \frac{\beta}{2} \right) \left( \sum_\ell u_\ell v_\ell \right),
    \end{equation}
    where
    \begin{equation} \label{eq-diagonal-rho-parameters-beta}
        \beta := - \log \left( \frac{ \sum_\ell\abs{u_\ell}^2}{\sum_\ell \abs{v_\ell}^2} \right)>0
    \end{equation}
    is the inverse temperature. The invariant density operator is then the thermal state
    \begin{equation*}
        \rho = (1 - \me^{-\beta}) \me^{-\beta N},
    \end{equation*}
    as shown in \eqref{eq-introduction-thermal-state}.

    \subsection{Eigenvalues of the Induced Generator \texorpdfstring{$ L^{(s)} $}{L(s)} }

    For a QMS satisfying Assumptions~\ref{hyp-stability-Z}--\ref{hyp-faithfulness}, let $(T_t^{(s)})_{t\geq0}$ be the induced semigroup defined by \eqref{eq-induced-semigroup-and-qms}, with $s\in[0,1]$. We begin by identifying the eigenvalues of its generator $L^{(s)}$. The following result, whose proof is given in Appendix~\ref{section-proof-theorem-L-calL}, provides the fundamental eigenvalues.

    \begin{theorem} \label{theorem-base-eigenvalue-L}
        The two eigenvalues of the drift matrix $\mathbf Z$ are eigenvalues of $L^{(s)}$. If $\mathbf Z$ is diagonalizable, the corresponding eigenvectors have the form $\rho^{s/2}X\rho^{(1-s)/2}$, where $X$ is a polynomial of degree one in $a$ and $a^\dagger$. If $\mathbf Z$ is defective, the restriction of $L^{(s)}$ to the embedded polynomials of degree at most one has an eigenvector and a generalized eigenvector associated with the repeated fundamental eigenvalue.
    \end{theorem}

    The remaining polynomial eigenvalues follow from the quasi-derivation identity in Theorem~\ref{theorem-derivation-property-qms}. To transfer that identity from the algebraic Gaussian generator $\cL$ to the induced generator $L^{(s)}$ on the Hilbert--Schmidt space $\cB_2(\mathsf h)$, we use the following domain result, proved in Appendix~\ref{section-proof-theorem-L-calL}.

    \begin{theorem} \label{theorem-L-and-calL-relation}
        Let $X$ be a polynomial in $a$ and $a^\dagger$. Then $\rho^{s/2}X\rho^{(1-s)/2}\in\dom L^{(s)}$ and
        \begin{equation} \label{eq-L-and-calL-relation-coro}
            L^{(s)} ( \rho^{s/2} X \rho^{(1-s)/2} ) = \rho^{s/2} \cL (X) \rho^{(1-s)/2}.
        \end{equation}
        Here $\cL(X)$ is interpreted algebraically through \eqref{eq-gaussian-gksl-generator}.
    \end{theorem}

    We can now generate all eigenvalues arising from embedded polynomials. Recall from \eqref{eq-introduction-polynomial-spaces} that, for $r\in\N$, $\cX_r$ is the finite-dimensional space of embedded polynomials of total degree at most $r$; let $\cP_r$ denote the orthogonal projection of $\cB_2(\mathsf h)$ onto $\cX_r$.

    \begin{theorem} \label{theorem-other-eigenvalues-L}
        Suppose first that $\mathbf Z$ is diagonalizable, and let $X$ and $Y$ be first-order polynomials in $a$ and $a^\dagger$ such that
        \begin{align}
            L^{(s)} ( \rho^{s/2} X \rho^{(1-s)/2} ) &= \lambda \rho^{s/2} X \rho^{(1-s)/2}, \label{eq-L-base-eigenvectors-1} \\
            \quad L^{(s)} (\rho^{s/2} Y \rho^{(1-s)/2}) &= \mu \rho^{s/2} Y \rho^{(1-s)/2}, \label{eq-L-base-eigenvectors-2}
        \end{align}
        where $ \lambda, \mu \in \C $. For all $ n, m \in \N $, we have
        \begin{align}
            L^{(s)} ( \rho^{s/2} X^n Y^m \rho^{(1-s)/2} ) &= (n \lambda + m \mu) \rho^{s/2} X^n Y^m \rho^{(1-s)/2} \nonumber \\
            &\quad + \rho^{s/2} R_{n, m} \rho^{(1-s)/2}, \label{eq-L-recursive-Xn-Ym}
        \end{align}
        where $R_{n,m}\in\myspan\{X^iY^j:i+j\leq n+m-2\}$ when $n+m\geq2$, and $R_{n,m}=0$ otherwise. Moreover, every number $n\lambda+m\mu$, with $n,m\in\N$, is an eigenvalue of $L^{(s)}$.

        Suppose now that $\mathbf Z$ is defective, and let $X$ and $Y$ be first-order polynomials in $a$ and $a^\dagger$ such that
        \begin{align}
            L^{(s)} ( \rho^{s/2} X \rho^{(1-s)/2} ) &= \lambda \rho^{s/2} X \rho^{(1-s)/2}, \label{eq-L-eigenvector-and-generalized-eigenvector-1} \\
            L^{(s)} ( \rho^{s/2} Y \rho^{(1-s)/2} ) &= \lambda \rho^{s/2} Y \rho^{(1-s)/2} - \rho^{s/2} X \rho^{(1-s)/2}, \label{eq-L-eigenvector-and-generalized-eigenvector-2}
        \end{align}
        where $ \lambda \in \C $. For all $ n, m \in \N $, we have
        \begin{align}
            L^{(s)} ( \rho^{s/2} X^n Y^m \rho^{(1-s)/2} ) &= (n + m ) \lambda \rho^{s/2} X^n Y^m \rho^{(1-s)/2} \nonumber \\
            &\quad - m \rho^{s/2} X^{n+1} Y^{m-1} \rho^{(1-s)/2} + \rho^{s/2} R_{n, m} \rho^{(1-s)/2},  \label{eq-L-X-n-Y-m-defective-case}
        \end{align}
        where $R_{n,m}\in\myspan\{X^iY^j:i+j\leq n+m-2\}$ when $n+m\geq2$, and $R_{n,m}=0$ otherwise. Furthermore, every $n\lambda$, with $n\in\N$, is an eigenvalue of $L^{(s)}$.
    \end{theorem}
    \begin{proof}
        We argue by induction. First suppose that $\mathbf Z$ is diagonalizable. Theorem~\ref{theorem-base-eigenvalue-L} guarantees the existence of embedded first-order polynomials satisfying \eqref{eq-L-base-eigenvectors-1} and \eqref{eq-L-base-eigenvectors-2}; these identities establish \eqref{eq-L-recursive-Xn-Ym} for $(n,m)=(1,0)$ and $(0,1)$.

        Suppose that \eqref{eq-L-recursive-Xn-Ym} holds for a given pair $(n,m)$. The quasi-derivation identity, Theorem~\ref{theorem-L-and-calL-relation}, and Lemma~\ref{lemma-X-rho-rho-hat-X} give
        \begin{align*}
            &\quad L^{(s)} ( \rho^{s/2} X^{n+1} Y^m \rho^{(1-s)/2} ) = \rho^{s/2} \cL ( X X^n Y^m ) \rho^{(1-s)/2} \\
            &= \rho^{s/2} X \cL (X^n Y^m) \rho^{(1-s)/2} + \rho^{s/2}  \cL (X) X^n Y^m \rho^{(1-s)/2} \\
            &\quad + \rho^{s/2} \sum_\ell  [X, L_\ell^*] [L_\ell, X^n Y^m] \rho^{(1-s)/2} \\
            &= \rho^{s/2} X \rho^{-s/2} \cdot L^{(s)} ( \rho^{s/2} X^n Y^m \rho^{(1-s)/2}) \\
            &\quad + L^{(s)} (\rho^{s/2} X \rho^{(1-s)/2} ) \cdot \rho^{-(1-s)/2} X^n Y^m \rho^{(1-s)/2} \\
            &\quad + \rho^{s/2} \sum_\ell [X, L_\ell^*] \left( X^n [L_\ell, Y^m] + [L_\ell, X^n] Y^m \right) \rho^{(1-s)/2} \\
            &= ((n+1) \lambda + m \mu) \rho^{s/2} X^{n+1} Y^{m} \rho^{(1-s)/2} + \rho^{s/2} X R_{n, m} \rho^{(1-s)/2} \\
            &\quad + \rho^{s/2} \sum_\ell \left( m[X, L_\ell^*] [L_\ell, Y] X^n Y^{m-1} + n[L_\ell, X][X, L_\ell^*] X^{n-1} Y^m \right) \rho^{(1-s)/2}.
        \end{align*}
        Since $ R_{n, m} \in \myspan \{ X^i Y^j : i + j \le n + m - 2 \} $, we have
        \begin{equation} \label{eq-L-eigenvectors-tail-part-1}
            X R_{n, m} \in \myspan \{ X^i Y^j : i + j \le (n+1) + m - 2 \}.
        \end{equation}
        Each of $[X,L_\ell^*]$, $[L_\ell,X]$, and $[L_\ell,Y]$ is a scalar multiple of $\1$. Together with
        \begin{equation*}
            X^n Y^{m-1}, X^{n-1} Y^m \in \myspan \{ X^i Y^j : i + j \le (n+1) + m - 2 \}
        \end{equation*}
        and \eqref{eq-L-eigenvectors-tail-part-1}, this proves \eqref{eq-L-recursive-Xn-Ym} for $(n+1,m)$. The same argument proves it for $(n,m+1)$, and induction completes the proof for all $n,m\in\N$.

        Because $X$ and $Y$ are linearly independent first-order polynomials,
        \begin{equation*}
            \cX_r = \myspan \left\{ \rho^{s/2} X^n Y^m \rho^{(1-s)/2} : n + m \le r \right\}.
        \end{equation*}
        Lemma~\ref{lemma-cL-X-still-polynomial} shows that $\cX_r$ is invariant under $L^{(s)}$. Order the following basis first by increasing total degree and, within each degree, by decreasing powers of $X$:
        \begin{equation} \label{eq-vectors-X-Y-basis-truncated}
        \rho^{s/2} \left( \1, X, Y, X^2, XY, Y^2, \cdots, X^r, X^{r-1} Y, \cdots, Y^r \right) \rho^{(1-s)/2},
        \end{equation}
        These vectors are linearly independent, and every remainder term has strictly smaller degree. Hence the matrix of $L^{(s)}\restriction\cX_r$ is upper triangular, with diagonal entries $n\lambda+m\mu$ for $n+m\leq r$. Letting $r$ vary proves that every $n\lambda+m\mu$ is an eigenvalue of $L^{(s)}$.

        When the drift matrix $\mathbf Z$ is defective, the induction leading to \eqref{eq-L-X-n-Y-m-defective-case} begins with a Jordan chain. Indeed, $\mathbf Z$ is defective if and only if $L^{(s)}$ is defective on the first-order subspace
        \[
            \myspan\bigl\{\rho^{s/2}a\rho^{(1-s)/2},
            \rho^{s/2}a^\dagger\rho^{(1-s)/2}\bigr\}.
        \]
        Hence there exist first-order polynomials $X$ and $Y$ satisfying \eqref{eq-L-eigenvector-and-generalized-eigenvector-1} and \eqref{eq-L-eigenvector-and-generalized-eigenvector-2}; in particular, $\rho^{s/2}Y\rho^{(1-s)/2}$ is a generalized eigenvector of $L^{(s)}$. These two identities establish \eqref{eq-L-X-n-Y-m-defective-case} for $(n,m)=(1,0)$ and $(0,1)$.

        Suppose that \eqref{eq-L-X-n-Y-m-defective-case} holds for a given pair $(n,m)$. The quasi-derivation identity, Theorem~\ref{theorem-L-and-calL-relation}, and Lemma~\ref{lemma-X-rho-rho-hat-X} give
        \begin{align*}
            &\quad L^{(s)} ( \rho^{s/2} X^{n} Y^{m+1} \rho^{(1-s)/2}) = \rho^{s/2} \cL (X^n Y^m Y) \rho^{(1-s)/2} \\
            &= \rho^{s/2} X^n Y^m \cL (Y) \rho^{(1-s)/2} + \rho^{s/2} \cL ( X^{n} Y^m ) Y \rho^{(1-s)/2} \\
            &\quad + \rho^{s/2} \sum_\ell [X^n Y^m, L_\ell^*] [L_\ell, Y] \rho^{(1-s)/2} \\
            &= \rho^{s/2} X^n Y^m \rho^{-s/2} \cdot L^{(s)} (\rho^{s/2} Y \rho^{(1-s)/2}) \\
            &\quad + L^{(s)} ( \rho^{s/2} X^n Y^m \rho^{(1-s)/2}) \cdot \rho^{-(1-s)/2} Y \rho^{(1-s)/2} \\
            &\quad + \rho^{s/2} \sum_\ell \left( X^n [Y^m, L_\ell^*] + [X^n, L_\ell^*] Y^m \right) [L_\ell, Y] \rho^{(1-s)/2} \\
            &= \lambda \rho^{s/2} X^n Y^{m+1} \rho^{(1-s)/2} - \rho^{s/2} X^{n} Y^m X \rho^{(1-s)/2} \\
            &\quad + (n + m) \lambda \rho^{s/2} X^n Y^{m+1} \rho^{(1-s)/2} \\
            &\quad - m \rho^{s/2} X^{n+1} Y^m \rho^{(1-s)/2} + \rho^{s/2} R_{n, m} Y \rho^{(1-s)/2} \\
            &\quad + \rho^{s/2} \sum_\ell ( m [Y, L_\ell^*] [L_\ell, Y] X^n Y^{m-1} + n [X, L_\ell^*] [L_\ell, Y] X^{n-1} Y^m) \rho^{(1-s)/2} \\
            &= (n+m+1) \lambda \rho^{s/2} X^n Y^{m+1} \rho^{(1-s)/2} - (m + 1) \rho^{s/2} X^{n+1} Y^m \rho^{(1-s)/2} \\
            &- m \rho^{s/2} [Y, X] X^n Y^{m-1} \rho^{(1-s)/2} + \rho^{s/2} R_{n, m} Y \rho^{(1-s)/2} \\
            &\quad + \rho^{s/2} \sum_\ell ( m [Y, L_\ell^*] [L_\ell, Y] X^n Y^{m-1} + n [X, L_\ell^*] [L_\ell, Y] X^{n-1} Y^m) \rho^{(1-s)/2}
        \end{align*}
        The condition on $R_{n,m}$ implies
        \begin{equation} \label{eq-a-tail-defective}
            R_{n, m} Y \in \myspan \{ X^i Y^j : i + j \le n + (m+1) - 2 \}.
        \end{equation}
        Moreover, $[Y,X]$, $[L_\ell,Y]$, $[Y,L_\ell^*]$, and $[X,L_\ell^*]$ are scalar multiples of $\1$. Together with
        \begin{equation*}
            X^n Y^{m-1}, X^{n-1} Y^m \in \myspan \{ X^i Y^j : i + j \le n + (m+1) -2 \}
        \end{equation*}
        and \eqref{eq-a-tail-defective}, this proves \eqref{eq-L-X-n-Y-m-defective-case} for $(n,m+1)$. The case $(n+1,m)$ follows in the same way, so induction proves \eqref{eq-L-X-n-Y-m-defective-case} for all $n,m\in\N$.

        With respect to the ordered basis in \eqref{eq-vectors-X-Y-basis-truncated}, the matrix of $L^{(s)}\restriction\cX_r$ is therefore upper triangular, with diagonal value $(n+m)\lambda$ on the homogeneous component of degree $n+m$. Letting $r$ vary proves that every $n\lambda$, with $n\in\N$, is an eigenvalue of $L^{(s)}$.
    \end{proof}

    \subsection{Eigenvalues of the Adjoint Generator \texorpdfstring{$ L^{(s)*} $}{L(s)*}}

    We now turn to the adjoint generator $L^{(s)*}$. Its analysis is based on the dual QMS $(\cT_t^\prime)_{t\geq0}$, characterized by
    \begin{equation} \label{eq-definition-dual-QMS}
        \Tr\!\left(\rho^{1/2}\cT_t^\prime(X)\rho^{1/2}Y\right)
        =\Tr\!\left(\rho^{1/2}X\rho^{1/2}\cT_t(Y)\right),
        \qquad X,Y\in\cB(\mathsf h).
    \end{equation}
    The existence of $(\cT_t^\prime)_{t\geq0}$ as a QMS follows from \cite[Theorem~2.1]{albeverio1978frobenius}; see also \cite[Theorem~1]{fagnola2010generators}. The adjoint semigroup $(T_t^{(s)*})_{t\geq0}$ is the induced semigroup associated with $(\cT_t^\prime)_{t\geq0}$. Since the dual QMS is again Gaussian, its quasi-derivation identity yields the fundamental eigenvalues and their non-negative integer combinations exactly as for $L^{(s)}$.

    \begin{lemma} \label{lemma-generator-of-KMS-dual-QMS}
        The generator of $(\cT^\prime_t)_{t\geq0}$ is
        \begin{equation} \label{eq-gksl-KMS-dual-generator}
            \cL^\prime (X) = \mi [H^\prime, X] - \frac{1}{2} \sum_\ell \left( L_\ell^{\prime *}  L_\ell^\prime X - 2 L_\ell^{\prime *} X L_\ell^\prime + X L_\ell^{\prime *} L_\ell^\prime \right), \quad X \in \dom \cL^\prime,
        \end{equation}
        where
        \begin{equation} \label{eq-Hamiltonian-Kraus-KMS-dual-QMS}
            H^\prime = - \Omega a^\dagger a + \frac{\kappa}{2} a^{\dagger 2} + \frac{\overline{\kappa}}{2} a^2, \quad
            L_\ell^\prime = \me^{\beta/2} \overline{u_\ell} a + \me^{-\beta/2} v_\ell a^\dagger.
        \end{equation}
    \end{lemma}
    \begin{proof}
        Equation~\eqref{eq-definition-dual-QMS} is equivalent to
        \begin{equation} \label{eq-dual-qms-and-predual-semegroup}
            \Tr \left( \rho^{1/2} \cT_t^\prime (X) \rho^{1/2} Y \right) = \Tr \left( \cT_{*t} ( \rho^{1/2} X \rho^{1/2} ) Y \right), \quad \forall X, Y \in \cB (\mathsf{h}),
        \end{equation}
        and therefore the dual QMS and the predual semigroup satisfy
        \begin{equation*}
            \rho^{1/2} \cT_t^\prime (X) \rho^{1/2} = \cT_{*t} (\rho^{1/2} X \rho^{1/2}), \quad \forall X \in \cB (\mathsf{h}).
        \end{equation*}
        If $X\in\cB(\mathsf h)$ and $\rho^{1/2}X\rho^{1/2}\in\dom\cL_*$, differentiating \eqref{eq-dual-qms-and-predual-semegroup} at $t=0$ yields
        \begin{equation} \label{eq-dual-generator-and-predual-generator}
            \rho^{1/2} \cL^\prime (X) \rho^{1/2} = \cL_* (\rho^{1/2} X \rho^{1/2}).
        \end{equation}
        Set $G:=-\mi H-\frac12\sum_\ell L_\ell^*L_\ell$. Equation~\eqref{eq-gksl-gaussian-predual-generator} gives
        \begin{equation} \label{eq-rho-1/2-dual-L-rho-1/2}
            \rho^{1/2} \cL^\prime (X) \rho^{1/2} = G \rho^{1/2} X \rho^{1/2} + \sum_\ell L_\ell \rho^{1/2} X \rho^{1/2} L_\ell^* + \rho^{1/2} X \rho^{1/2} G^*,
        \end{equation}
        Since $H$ and $L_\ell$ are polynomials in $a$ and $a^\dagger$, Lemma~\ref{lemma-a-a-dagger-exponential-of-N} gives polynomials $G^\prime$ and $L_\ell^\prime$ satisfying
        \begin{equation} \label{eq-G-L-dual}
            G^{\prime *} = \rho^{-1/2} G \rho^{1/2} , \quad L^{\prime *}_\ell = \rho^{-1/2} L_\ell \rho^{1/2}.
        \end{equation}
        Consequently, \eqref{eq-rho-1/2-dual-L-rho-1/2} simplifies to
        \begin{equation*}
           \cL^\prime (X) = G^{\prime *} X + \sum_\ell L_\ell^{\prime *} X L_\ell^\prime + X G^\prime.
        \end{equation*}
        Defining $H^\prime:=\mi G^\prime+\frac{\mi}{2}\sum_\ell L_\ell^{\prime *}L_\ell^\prime$ gives \eqref{eq-gksl-KMS-dual-generator}. Direct computation from \eqref{eq-G-L-dual} yields
        \begin{align*}
            G^\prime &= ( \mi \Omega - \frac{1}{2} \sum_\ell ( \abs{u_\ell}^2 + \abs{v_\ell}^2 )) a^\dagger a - \sum_\ell \frac{\abs{u_\ell}^2}{2} \1 \\
            &\quad + \frac{\me^{-\beta}}{2} ( \mi \kappa - \sum_\ell u_\ell v_\ell ) a^{\dagger 2} + \frac{\me^{\beta}}{2} (\mi \overline{\kappa} - \sum_\ell \overline{u_\ell} \overline{v_\ell}) a^2, \\
            L_\ell^\prime &= \me^{\beta/2} \overline{u_\ell} a + \me^{-\beta/2} v_\ell a^\dagger, \\
            L^{\prime *}_\ell L^{\prime}_\ell &= ( \me^\beta \abs{u_\ell}^2 + \me^{-\beta} \abs{v_\ell}^2 ) a^\dagger a + \me^{-\beta} \abs{v_\ell}^2 \1 + u_\ell v_\ell a^{\dagger 2} + \overline{u_\ell} \overline{v_\ell} a^2.
        \end{align*}
        Combining these identities with \eqref{eq-diagonal-rho-parameters} and \eqref{eq-diagonal-rho-parameters-beta}, we obtain
        \begin{align*}
            H^\prime &= \left(-\Omega + \frac{\mi}{2} \sum_\ell ( (\me^{-\beta} - 1 ) \abs{v_\ell}^2 + (\me^{\beta} - 1) \abs{u_\ell}^2 ) \right) a^\dagger a \\
            &\quad + \frac{1}{2} \left( - \kappa \me^{-\beta} + \mi (1 - \me^{-\beta} ) \sum_\ell u_\ell v_\ell \right) a^{\dagger 2} \\
            &\quad + \frac{1}{2} \left( - \overline{\kappa} \me^{\beta} + \mi (1 - \me^\beta) \sum_\ell \overline{u_\ell} \overline{v_\ell} \right) a^2 \\
            &\quad + \frac{\mi}{2} \left( \sum_\ell (\me^{-\beta} \abs{v_\ell}^2 - \abs{u_\ell}^2 ) \right) \1 \\
            &= - \Omega a^\dagger a + \frac{\kappa}{2} a^{\dagger 2} + \frac{\overline{\kappa}}{2} a^2,
        \end{align*}
        which proves \eqref{eq-Hamiltonian-Kraus-KMS-dual-QMS}.
    \end{proof}

    Equation~\eqref{eq-Hamiltonian-Kraus-KMS-dual-QMS} shows that $\cL^\prime$ has Gaussian GKSL form, with $\deg(H^\prime)\leq2$ and $\deg(L_\ell^\prime)\leq1$. Thus the dual QMS is Gaussian.

    \begin{corollary} \label{corollary-drift-matrix-of-dual-QMS-diagonal-case}
        The drift matrix of $ (\cT_t^\prime)_{t \ge 0} $, denoted by $ \mathbf{Z}^\prime $, is given by
        \begin{equation} \label{eq-drift-matrix-the-dual-qms}
            \mathbf{Z}^\prime = \begin{pmatrix}
                - \gamma - \Im \kappa & \Re \kappa + \Omega \\
                \Re \kappa - \Omega & - \gamma + \Im \kappa
            \end{pmatrix}.
        \end{equation}
        Consequently,
        \begin{equation*}
            \mathbf{Z}^\prime = \mathbf{Z}^T,
        \end{equation*}
        where $\mathbf Z^T$ is the transpose of the real drift matrix $\mathbf Z$.
    \end{corollary}
    \begin{proof}
        Compare \eqref{eq-Hamiltonian-Kraus-KMS-dual-QMS} with \eqref{eq-Hamiltonian-with-parameters} and \eqref{eq-kraus-operator-a-adagger}. Passing from $\cL$ to $\cL^\prime$ replaces $\Omega$ by $-\Omega$ and leaves $\kappa$ unchanged; moreover, $\zeta=0$ because the invariant state is diagonal. The parameters $v_\ell$ and $u_\ell$ become $\me^{\beta/2}u_\ell$ and $\me^{-\beta/2}v_\ell$, respectively. Equations~\eqref{eq-introduction-drift-operator} and \eqref{eq-diagonal-rho-parameters-beta} give
        \begin{equation*}
            \frac{1}{2} \sum_\ell \left( \abs{ \me^{\beta/2} u_\ell }^2 - \abs{ \me^{-\beta/2} v_\ell }^2 \right) = \gamma = \frac{1}{2} \sum_\ell \left( \abs{v_\ell}^2 - \abs{u_\ell}^2 \right),
        \end{equation*}
        which proves \eqref{eq-drift-matrix-the-dual-qms}.
    \end{proof}

    The adjoint semigroup $(T_t^{(s)*})_{t\geq0}$ is the semigroup induced by the dual QMS with embedding parameter $1-s$.

    \begin{lemma} \label{lemma-adjoint-semigroup-and-the-dual-qms}
        We have
        \begin{equation*}
            T_t^{(s)*} \left( \rho^{(1-s)/2} X \rho^{s/2} \right) = \rho^{(1-s)/2} \cT^\prime_t (X) \rho^{s/2}, \quad \forall X \in \cB (\mathsf{h}).
        \end{equation*}
    \end{lemma}
    \begin{proof}
        For $X,Y\in\cB(\mathsf h)$,
        \begin{align*}
            &\quad \left\langle \rho^{(1-s)/2} \cT_t^\prime (X^*) \rho^{s/2}, \rho^{s/2} Y \rho^{(1-s)/2} \right\rangle_2 \\
            &= \Tr \left( \rho^{s/2} \cT_t^\prime (X) \rho^{(1-s)/2} \rho^{s/2} Y \rho^{(1-s)/2} \right) \\
            &= \Tr \left( \rho^{1/2} \cT_t^\prime (X) \rho^{1/2} Y \right) = \Tr \left( \rho^{1/2} X \rho^{1/2} \cT_t (Y) \right) \\
            &= \Tr \left( \rho^{s/2} X \rho^{(1-s)/2} T_t^{(s)} (\rho^{s/2} Y \rho^{(1-s)/2}) \right) \\
            &= \left\langle T_t^{(s) *} (\rho^{(1-s)/2} X^* \rho^{s/2}), \rho^{s/2} Y \rho^{(1-s)/2} \right\rangle_2,
        \end{align*}
        The density of $\rho^{s/2}\cB(\mathsf h)\rho^{(1-s)/2}$ in $\cB_2(\mathsf h)$ completes the proof.
    \end{proof}

    Lemma~\ref{lemma-adjoint-semigroup-and-the-dual-qms} allows us to apply Theorems~\ref{theorem-base-eigenvalue-L} and~\ref{theorem-other-eigenvalues-L} to $L^{(s)*}$, using the parameters of the dual QMS. Corollary~\ref{corollary-drift-matrix-of-dual-QMS-diagonal-case} therefore gives the following result.
    \begin{corollary} \label{theorem-adjoint-generator-eigenvalues}
        The eigenvalues of $L^{(s)*}$ associated with embedded-polynomial eigenvectors are the complex conjugates of the corresponding eigenvalues of $L^{(s)}$.
    \end{corollary}

    \subsection{Exhaustion of Point Spectra of \texorpdfstring{$ L^{(s)} $}{L(s)} and \texorpdfstring{$ L^{(s)*} $}{ L(s)* } }

    We now show that the polynomial eigenvalues exhaust the point spectra of $L^{(s)}$ and $L^{(s)*}$ and that every eigenspace is finite-dimensional. Recall that $\cX_r\subset\cB_2(\mathsf h)$ is the space of embedded polynomials of degree at most $r$, and that $\cP_r$ is the orthogonal projection onto $\cX_r$.

    Lemma~\ref{lemma-cL-X-still-polynomial} and Theorem~\ref{theorem-L-and-calL-relation} imply
    \begin{equation} \label{eq-invariance-X-r-under-induced-generator}
        L^{(s)} (\cX_r) \subset \cX_r.
    \end{equation}
    Applying Lemma~\ref{lemma-adjoint-semigroup-and-the-dual-qms} and Theorem~\ref{theorem-L-and-calL-relation} to the dual QMS, and using the invariance of polynomial degree under modular conjugation, also gives $L^{(s)*}(\cX_r)\subset\cX_r$.

    For a closed operator $A$ on a Hilbert space, we say that a closed subspace $M$ \emph{reduces} $A$ if the orthogonal projection $P_M$ maps $\dom A$ into itself and $AP_Mx=P_MAx$ for every $x\in\dom A$.

    \begin{lemma} \label{lemma-exchangeability-induced-generator-P}
        We have
        \begin{equation*}
            L^{(s)} \cP_r = \cP_r L^{(s)}, \quad L^{(s)^*} \cP_r = \cP_r L^{(s)^*}
        \end{equation*}
        on $ \dom L^{(s)} $ and $ \dom L^{(s)^*} $, respectively.
    \end{lemma}
    \begin{proof}
        Let $x\in\dom L^{(s)}$ and $y\in\cX_r$. Then
        \begin{align*}
            \langle \cP_r L^{(s)} x, y \rangle_2 &=
            \langle L^{(s)} x, y \rangle_2 =
            \langle  x, L^{(s)*} y \rangle_2 \\
            &= \langle \cP_r x, L^{(s)*} y \rangle_2 = \langle L^{(s)} \cP_r x,  y \rangle_2.
        \end{align*}
        By \eqref{eq-invariance-X-r-under-induced-generator}, $L^{(s)}\cP_r x$ belongs to $\cX_r$. Since $y\in\cX_r$ is arbitrary, $L^{(s)}\cP_r=\cP_rL^{(s)}$ on $\dom L^{(s)}$. The adjoint identity follows by the same argument.
    \end{proof}

    Recall from \eqref{eq-introduction-polynomial-spaces} and the definition following it that $\cX$ is the space of all embedded polynomials in $a$ and $a^\dagger$. It is dense in $\cB_2(\mathsf h)$ by Lemma~\ref{lemma-density-of-domain-L-star-plus-L}, so $\cP_r\to\1$ strongly as $r\to\infty$.

    \begin{lemma} \label{lemma-cX-common-core}
        $ \cX $ is a common core for $ L^{(s)} $ and $ L^{(s)*} $.
    \end{lemma}
    \begin{proof}
        Let $x\in\dom L^{(s)}$. Strong convergence of $\cP_r$ and Lemma~\ref{lemma-exchangeability-induced-generator-P} give
        \begin{equation*}
            \cP_r x \rightarrow x, \quad L^{(s)} \cP_r x = \cP_r L^{(s)} x \rightarrow L^{(s)} x.
        \end{equation*}
        Thus $\cX$ is a core for $L^{(s)}$. The same argument applies to $L^{(s)*}$, so $\cX$ is a common core.
    \end{proof}

    We write $\sigma_p$ for the point spectrum.

    \begin{proposition} \label{proposition-exhaustion-point-spectrum}
        If $\mathbf Z$ is diagonalizable, then
        \begin{equation*}
            \sigma_p (L^{(s)}) = \{ n \lambda + m \mu : n, m \in \N \}, \quad \sigma_p (L^{(s)*}) = \{ n \overline{\lambda} + m \overline{\mu} : n, m \in \N \}.
        \end{equation*}
        If $\mathbf Z$ is defective, then
        \begin{equation*}
            \sigma_p (L^{(s)}) = \{ n \lambda : n \in \N \}, \quad \sigma_p (L^{(s)*}) = \{ n \overline{\lambda} : n \in \N \}.
        \end{equation*}
        Every eigenvalue in these sets has finite multiplicity.
    \end{proposition}
    \begin{proof}
        Suppose that a nonzero $x\in\cB_2(\mathsf h)$ satisfies $L^{(s)}x=\alpha x$ for some $\alpha\notin\{n\lambda+m\mu:n,m\in\N\}$. Then $x\in\dom L^{(s)}$, and Lemma~\ref{lemma-exchangeability-induced-generator-P} gives
        \begin{equation*}
            L^{(s)} \cP_r x = \cP_r L^{(s)} x = \alpha \cP_r x.
        \end{equation*}
        Because $\alpha$ is not an eigenvalue of $L^{(s)}\restriction\cX_r$, it follows that $\cP_r x=0$ for every $r$. Strong convergence of $\cP_r$ then gives $x=0$, a contradiction. The same argument applies to $L^{(s)*}$ and, with the corresponding spectral list, to the defective case.

        For finite multiplicity, observe that the projections commute with both generators, so each $\cX_r\ominus\cX_{r-1}$ reduces $L^{(s)}$. Its spectrum on this space is given by the diagonal terms of total degree $r$. Since $\Re\lambda<0$ and $\Re\mu<0$,
        \[
            -\Re(n\lambda+m\mu)\geq(n+m)
            \min\{-\Re\lambda,-\Re\mu\},
        \]
        and a fixed complex number can occur at only finitely many total degrees. Every eigenspace is therefore finite-dimensional. The same argument applies to $L^{(s)*}$. In the defective case, the diagonal values on $\cX_r\ominus\cX_{r-1}$ are $r\lambda$ for $L^{(s)}$ and $r\overline\lambda$ for $L^{(s)*}$. Since $\Re\lambda<0$, a fixed eigenvalue occurs in at most one total degree, and its eigenspace is again finite-dimensional.
    \end{proof}

\begin{figure}[ht]
        \centering
        \begin{tikzpicture}
            \def\n{6}
            \def\lamre{-1}
            \def\lamim{0.35}
            \def\mure{-1}
            \def\muim{-0.35}

            % Draw axes (optional)
            \draw[->, darkgray] (-5,0) -- (1,0) node[right] {$ \Re $};
            \draw[->, darkgray] (0,-2) -- (0,2) node[above] {$ \Im $};

            % Scatter points
            \foreach \i in {0,...,5} {
                \foreach \j in {0,...,5} {
                    \pgfmathtruncatemacro{\cond}{\i + \j < 5 ? 1 : 0}
                    \ifnum\cond=1
                        \pgfmathsetmacro{\x}{\lamre * \i + \mure * \j}
                        \pgfmathsetmacro{\y}{\lamim * \i + \muim * \j}
                        \filldraw[black] (\x,\y) circle (1.1pt);
                    \fi
                }
            }
        \end{tikzpicture}
        \caption{Point spectrum of $L^{(s)}$ with nonreal complex-conjugate fundamental eigenvalues.}
        \label{fig:scatter}
    \end{figure}
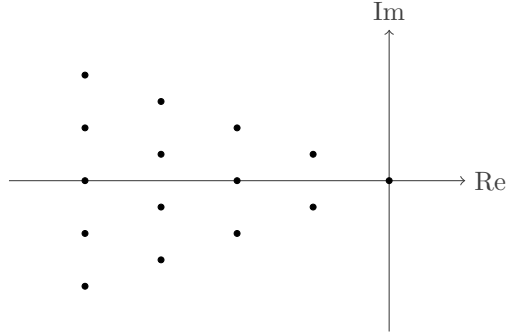

    \begin{remark}
        In the diagonalizable case, the eigenvalues of $\mathbf Z$ are either two real negative numbers or a nonreal complex-conjugate pair. In the latter case, every eigenvalue of $L^{(s)}$ is simple; Figure~\ref{fig:scatter} illustrates the resulting spectral pattern.
    \end{remark}

    % ===================================================================

    % NEW SECTION

    % ===================================================================

    \section{The Symmetric Part and the Spectral Gap} \label{section-spectrum-L-star-plus-L}

    We now study the symmetric part of the induced generator. Its eigenvalues retain the polynomial structure found in Section~\ref{section-eigenvalues-of-L}, while its largest nonzero eigenvalue determines the spectral gap. Its self-adjoint closure also provides the comparison operator used to prove compactness of the resolvent of $L^{(s)}$.

    For a general generator of a strongly continuous contraction semigroup, the domain of the generator and that of its adjoint can have trivial intersection; see \cite{arlinskiui2020everything}. Here, however, Lemma~\ref{lemma-cX-common-core} shows that the embedded-polynomial space $\cX$ is a common core for both operators. We therefore define $L^{(s)*}+L^{(s)}$ initially on $\cX$. Until the closure is introduced in \eqref{eq-definition-closure-of-L-star-plus-L}, this notation means $(L^{(s)*}+L^{(s)})\restriction\cX$.

    The following lemma computes its fundamental eigenvalues by direct algebraic calculation.

    \begin{lemma} \label{lemma-fundamental-eigenvalues-of-Ls*-plus-Ls}
        On the first-order embedded-polynomial space,

        \begin{align*}
            (L^{(s)*} + L^{(s)}) ( \rho^{s/2} a \rho^{(1-s)/2} ) &= - 2 \gamma \rho^{s/2} a \rho^{(1-s)/2} - \mi \kappa ( 1 + \me^{(2 s - 1) \beta} ) \rho^{s/2} a^\dagger \rho^{(1-s)/2}, \\
            ( L^{(s)*} + L^{(s)} ) (\rho^{s/2} a^\dagger \rho^{(1-s)/2} ) &= - 2 \gamma \rho^{s/2} a^\dagger \rho^{(1-s)/2} + \mi \overline{\kappa} ( 1 + \me^{(1 - 2s) \beta} ) \rho^{s/2} a \rho^{(1-s)/2},
        \end{align*}

        Consequently, the fundamental eigenvalues of $L^{(s)*}+L^{(s)}$ are
        \begin{equation} \label{eq-L-star-plus-L-fundamental-eigenvalues}
            - 2 \gamma \pm 2 \abs{\kappa} \cosh \left( (1-2s)\beta / 2 \right).
        \end{equation}
    \end{lemma}

    Unlike the fundamental eigenvalues of $L^{(s)}$, one of the fundamental eigenvalues in \eqref{eq-L-star-plus-L-fundamental-eigenvalues} may vanish. In that case, every spectral value of the self-adjoint closure has infinite multiplicity.

    \begin{theorem} \label{theorem-other-eigenvalues-L-star-plus-L}
        Let $X$ and $Y$ be first-order polynomials in $a$ and $a^\dagger$ such that
        \begin{align} \label{eq-L-star-plus-L-base-eigenvalues-1}
            (L^{(s)*} + L^{(s)}) (\rho^{s/2} X \rho^{(1-s)/2}) &= \lambda \rho^{s/2} X \rho^{(1-s)/2}, \\
            \quad (L^{(s)*} + L^{(s)}) (\rho^{s/2} Y \rho^{(1-s)/2}) &= \mu \rho^{s/2} Y \rho^{(1-s)/2}, \label{eq-L-star-plus-L-base-eigenvalues-2}
        \end{align}
        where $\mu\leq\lambda\leq0$ and $(\lambda,\mu)\neq(0,0)$. Then, for all $n,m\in\N$,
        \begin{equation} \label{eq-L-star-plus-L-with-residual-part}
            (L^{(s)*} + L^{(s)})( \rho^{\frac{s}{2}} X^n Y^m \rho^{\frac{1-s}{2}}) = (n \lambda + m \mu) \rho^{\frac{s}{2}} X^n Y^m \rho^{\frac{1-s}{2}} + \rho^{\frac{s}{2}} R_{n, m} \rho^{\frac{1-s}{2}},
        \end{equation}
        where $R_{n,m}\in\myspan\{X^iY^j:i+j\leq n+m-2\}$ when $n+m\geq2$, and $R_{n,m}=0$ otherwise. The eigenvalues of $(L^{(s)*}+L^{(s)})\restriction\cX$ are precisely the numbers $n\lambda+m\mu$, with $n,m\in\N$.
    \end{theorem}
    \begin{proof}
        Lemma~\ref{lemma-fundamental-eigenvalues-of-Ls*-plus-Ls} guarantees the existence of embedded polynomials satisfying \eqref{eq-L-star-plus-L-base-eigenvalues-1} and \eqref{eq-L-star-plus-L-base-eigenvalues-2}. Hence \eqref{eq-L-star-plus-L-with-residual-part} holds for $(n,m)=(1,0)$ and $(0,1)$. Assuming it for a given pair $(n,m)$, we obtain
        \begin{align}
            &\quad ( L^{(s)*} + L^{(s)} ) ( \rho^{s/2} X^{n+1} Y^m \rho^{(1-s)/2} ) \nonumber \\
            &= L^{(s)*} ( \rho^{(1-s)/2} \rho^{(2s-1)/2} X^{n+1} Y^m \rho^{(1-2s)/2} \rho^{s/2}) + L^{(s)} ( \rho^{s/2} X^{n+1} Y^m \rho^{(1-s)/2} ) \nonumber \\
            &= \rho^{(1-s)/2} \cL^\prime ( \rho^{(2s-1)/2} X X^{n} Y^m \rho^{(1-2s)/2} ) \rho^{s/2} + \rho^{s/2} \cL ( X X^{n} Y^m ) \rho^{(1-s)/2} \nonumber \\
            &= \rho^{s/2} X \rho^{-s/2} \cdot \rho^{(1-s)/2} \cL^\prime ( \rho^{(2s-1)/2} X^n Y^m \rho^{(1-2s)/2} ) \rho^{s/2} \nonumber \\
            &\quad + \rho^{(1-s)/2} \cL^\prime ( \rho^{(2s-1)/2} X \rho^{(1-2s)/2} ) \rho^{s/2} \cdot \rho^{(s-1)/2} X^n Y^m \rho^{(1-s)/2} \nonumber \\
            &\quad + \sum_\ell \rho^{(1-s)/2} [ \rho^{(2s-1)/2} X \rho^{(1-2s)/2}, L_\ell^{\prime *} ] [ L_\ell^\prime, \rho^{(2s-1)/2} X^n Y^m \rho^{(1-2s)/2} ] \rho^{s/2} \nonumber \\
            &\quad + \rho^{s/2} X \rho^{-s/2} \cdot \rho^{s/2} \cL ( X^n Y^m ) \rho^{(1-s)/2} \nonumber \\
            &\quad + \rho^{s/2} \cL(X) \rho^{(1-s)/2} \cdot \rho^{(s-1)/2} X^n Y^m \rho^{(1-s)/2} \nonumber \\
            &\quad + \sum_\ell \rho^{s/2} [X, L_\ell^*] [ L_\ell, X^n Y^m] \rho^{(1-s)/2} \nonumber \\
            &= \rho^{s/2} X \rho^{-s/2} \cdot (L^{(s)*} + L^{(s)}) ( \rho^{s/2} X^n Y^m \rho^{(1-s)/2} ) \nonumber \\
            &\quad + ( L^{(s)*} + L^{(s)} ) (\rho^{s/2} X \rho^{(1-s)/2}) \cdot \rho^{(s-1)/2} X^n Y^m \rho^{(1-s)/2}  \nonumber \\
            &\quad + \sum_\ell \rho^{(1-s)/2} [ \rho^{(2s-1)/2} X \rho^{(1-2s)/2}, L_\ell^{\prime *} ] [ L_\ell^\prime, \rho^{(2s-1)/2} X^n Y^m \rho^{(1-2s)/2} ] \rho^{s/2} \nonumber \\
            &\quad + \sum_\ell \rho^{s/2} [X, L_\ell^*] [ L_\ell, X^n Y^m] \rho^{(1-s)/2} \nonumber \\
            &= ( (n+1) \lambda + m \mu ) \rho^{s/2} X^{n+1} Y^m \rho^{(1-s)/2}  \nonumber \\
            &\quad + \sum_\ell \rho^{(1-s)/2} [ \rho^{(2s-1)/2} X \rho^{(1-2s)/2}, L_\ell^{\prime *} ] [ L_\ell^\prime, \rho^{(2s-1)/2} X^n Y^m \rho^{(1-2s)/2} ] \rho^{s/2}  \nonumber\\
            &\quad + \sum_\ell \rho^{s/2} [X, L_\ell^*] [ L_\ell, X^n Y^m] \rho^{(1-s)/2} + \rho^{s/2} X R_{n, m} \rho^{(1-s)/2}. \label{eq-action-of-L*-plus-L-on-X-n-Y-m}
        \end{align}
        As in the proof of Theorem~\ref{theorem-other-eigenvalues-L}, the remainder in \eqref{eq-action-of-L*-plus-L-on-X-n-Y-m} has total degree at most $(n+1)+m-2$. The analogous computation with $Y$ completes the induction and produces every value $n\lambda+m\mu$.

        It remains to prove exhaustion. Each space $\cX_r\ominus\cX_{r-1}$ reduces $(L^{(s)*}+L^{(s)})\restriction\cX$. The restriction to this finite-dimensional Hilbert space is symmetric and therefore has an orthonormal eigenbasis; its eigenvalues are exactly the diagonal values of total degree $r$ computed above. Taking the union of these orthonormal eigenbases yields an orthonormal basis of $\cB_2(\mathsf h)$. Thus the displayed values exhaust the point spectrum on $\cX$.
    \end{proof}

    The orthonormal eigenbasis constructed above yields essential self-adjointness.

    \begin{lemma} \label{lemma-L-star-plus-L-essential-sa}
    The operator $(L^{(s)*}+L^{(s)})\restriction\cX$ is essentially self-adjoint.
    \end{lemma}
    \begin{proof}
        By Theorem~\ref{theorem-other-eigenvalues-L-star-plus-L}, the eigenvectors of $(L^{(s)*}+L^{(s)})\restriction\cX$ form an orthonormal basis of $\cB_2(\mathsf h)$. Their finite linear span is a dense set of analytic vectors. Nelson's analytic-vector theorem \cite[Theorem~X.39, p.~202]{reed1975methods} therefore implies essential self-adjointness.
    \end{proof}

    From this point onward, the notation $L^{(s)*}+L^{(s)}$ refers to the self-adjoint closure
    \begin{equation} \label{eq-definition-closure-of-L-star-plus-L}
        L^{(s)*} + L^{(s)} := \overline{ L^{(s)*} + L^{(s)} \restriction \cX }.
    \end{equation}
    The following result shows that this closure introduces no new eigenvalues and identifies its full spectrum.

    \begin{lemma}
        We have
        \begin{equation} \label{eq-no-additional-eigenvalues-after-closure}
            \sigma_p ( L^{(s)*} + L^{(s)} ) = \sigma_p ( (L^{(s)*} + L^{(s)})\restriction \cX ) = \{ n \lambda + m \mu : n, m \in \N \},
        \end{equation}
        and
        \begin{equation} \label{eq-point-spectrum-and-full-spectrum-L-star-plus-L}
            \sigma ( L^{(s)*} + L^{(s)} ) = \sigma_p ( L^{(s)*} + L^{(s)} ).
        \end{equation}
    \end{lemma}
    \begin{proof}
        By Theorem~\ref{theorem-other-eigenvalues-L-star-plus-L}, there is an orthonormal basis $(f_k)_{k\in\N}$ of $\cB_2(\mathsf h)$, contained in $\cX$, such that
        \begin{equation*}
            ( L^{(s) *} + L^{(s)} ) f_k = \alpha_k f_k,
        \end{equation*}
        where
        \begin{equation*}
            \{ \alpha_k : k \in \N \} = \{ n \lambda + m \mu : n, m \in \N \},
        \end{equation*}
        with eigenvalues repeated according to their multiplicities. Moreover, $\cX=\myspan\{f_k:k\in\N\}$. Under the unitary identification $f_k\mapsto e_k$ of $\cB_2(\mathsf h)$ with $\ell^2(\N)$, the restriction $(L^{(s)*}+L^{(s)})\restriction\cX$ becomes the diagonal operator on finitely supported sequences with diagonal entries $(\alpha_k)_{k\in\N}$. Its closure is the maximal diagonal operator, whose domain consists of the sequences for which $(\alpha_kc_k)_k$ is square-summable; see \cite[Section~1.5, Exercise~12]{schmudgen2012unbounded}. Consequently,
        \begin{equation} \label{eq-domain-L-star-plus-L-explicit}
            \dom ( L^{(s)*} + L^{(s)} ) = \left\{ \sum_{k \ge 0} c_k f_k : \sum_{k \ge 0} ( 1 +  \abs{\alpha_k}^2 ) \abs{c_k}^2 < + \infty \right\},
        \end{equation}
        and
        \begin{equation*}
            \left( L^{(s)*} + L^{(s)} \right) \left( \sum_{k \ge 0} c_k f_k \right) = \sum_{k \ge 0} \alpha_k c_k f_k.
        \end{equation*}

        Suppose that $(L^{(s)*}+L^{(s)})x=\alpha x$ for a nonzero vector $x=\sum_{k\geq0}c_kf_k$ in the domain given by \eqref{eq-domain-L-star-plus-L-explicit}. Then $(\alpha_k-\alpha)c_k=0$ for every $k$. Since some $c_k$ is nonzero, $\alpha=\alpha_k$ for at least one $k$. Conversely, every $f_k$ remains an eigenvector of the closure. This proves \eqref{eq-no-additional-eigenvalues-after-closure}.

        Moreover, \cite[Example~2.2]{schmudgen2012unbounded} gives
        \begin{equation*}
            \sigma ( L^{(s)*} + L^{(s)} ) = \overline{ \{ \alpha_k : k \in \N \} } = \overline{ \{ n \lambda + m \mu : n, m \in \N \} }.
        \end{equation*}
        This set is already closed. If $\lambda<0$, then $\mu\leq\lambda<0$, and every bounded interval contains only finitely many values $n\lambda+m\mu$. If $\lambda=0$, then $\mu<0$, and the set of distinct values is $\{m\mu:m\in\N\}$. Consequently,
        \begin{equation*}
            \sigma ( L^{(s)*} + L^{(s)} ) = \sigma_p ( L^{(s)*} + L^{(s)} ) = \{ n \lambda + m \mu : n, m \in \N \},
        \end{equation*}
        which proves \eqref{eq-point-spectrum-and-full-spectrum-L-star-plus-L}. When $\lambda=0$, closure may enlarge an eigenspace by admitting infinite square-summable combinations of eigenvectors with the same eigenvalue, but it introduces no new eigenvalue. When $\lambda,\mu<0$, every eigenspace is finite-dimensional and is therefore unchanged by closure.
    \end{proof}

    The next lemma identifies the invariant vectors of every induced semigroup. Frigerio proved the case $s=0$ in \cite[proof of Lemma~2]{frigerio1977quantum}; the argument below works for all $s\in[0,1]$ and does not use the Gaussian structure of $\rho$.

    \begin{lemma}
        Let $\rho$ be the unique faithful normal invariant state of a QMS $(\cT_t)_{t\geq0}$ on $\cB(\mathsf h)$. Then the invariant-vector space of $(T_t^{(s)})_{t\geq0}$ is $\C\rho^{1/2}$.
    \end{lemma}
    \begin{proof}
        Under this assumption, the fixed-point space of $(\cT_t)_{t\geq0}$ is $\C\1$, and the invariant-vector space of the predual semigroup is $\C\rho$ \cite[Lemma~2]{frigerio1977quantum}.

        The definition \eqref{eq-induced-semigroup-and-qms} gives
        \begin{equation*}
            T_t^{(s)} ( \rho^{1/2} ) = \rho^{s/2} \cT_t (\1)  \rho^{(1-s)/2} = \rho^{1/2}.
        \end{equation*}

        Conversely, suppose that $Y\in\cB_2(\mathsf h)$ and $T_t^{(s)}Y=Y$. Then $T_t^{(s)*}Y=Y$, because
        \begin{align*}
            \norm{ T_t^{(s)*} (Y) - Y }_2^2 &= \norm{T_t^{(s)*} (Y)}_2^2 + \norm{Y}_2^2 - 2 \Re \langle T_t^{(s)*} Y, Y \rangle_2 \\
            &\le 2 \norm{Y}_2^2 - 2 \Re \langle Y, T_t^{(s)} Y \rangle_2 = 0.
        \end{align*}
        H\"older's inequality shows that $\rho^{(1-s)/2}Y^*\rho^{s/2}$ is trace class. For every $X\in\cB(\mathsf h)$,
        \begin{equation*}
            \Tr ( \rho^{(1-s)/2} Y^* \rho^{s/2} \cT_t (X) ) = \Tr ( Y^* T_t^{(s)} ( \rho^{s/2} X \rho^{(1-s)/2} ) ) = \Tr ( \rho^{(1-s)/2} Y^* \rho^{s/2} X ),
        \end{equation*}
        so $\rho^{(1-s)/2}Y^*\rho^{s/2}$ is invariant under the predual semigroup and hence belongs to $\C\rho$. Faithfulness of $\rho$ then implies $Y\in\C\rho^{1/2}$.
    \end{proof}

    Consequently, $\ker L^{(s)}=\ker L^{(s)*}=\C\rho^{1/2}$.

    \begin{theorem} \label{theorem-spectral-gap-compact-resolvent}
        The following statements are equivalent:
        \begin{enumerate}
            \itemsep0em
            \item \label{item-spectral-gap-theorem-1} $ L^{(s)} $ has a spectral gap.
            \item \label{item-spectral-gap-theorem-2} $ L^{(s)*} $ has a spectral gap.
            \item \label{item-spectral-gap-theorem-3} Both fundamental eigenvalues of $L^{(s)*}+L^{(s)}$ are strictly negative; equivalently, $0$ is a simple eigenvalue of $L^{(s)*}+L^{(s)}$.
            \item \label{item-spectral-gap-theorem-4} The self-adjoint operator $L^{(s)*}+L^{(s)}$ has compact resolvent.
        \end{enumerate}
    \end{theorem}
    \begin{proof}
        The equivalence of \eqref{item-spectral-gap-theorem-1} and \eqref{item-spectral-gap-theorem-2} follows from \cite[Lemma~2.6]{frigerio1993asymptotic}.

        We first prove that \eqref{item-spectral-gap-theorem-1} implies \eqref{item-spectral-gap-theorem-3}. By definition of the spectral gap \eqref{eq-definition-spectral-gap} and Lemma \ref{lemma-cX-common-core},
        \begin{align}
            \gap L^{(s)} &= \inf \left\{ - \Re \langle x, L^{(s)} x \rangle_2 : x \in \dom L^{(s)}, \norm{x}_2 = 1, x \perp \rho^{1/2} \right\} \nonumber \\
            &= \inf \left\{ - \frac{1}{2} \langle x, (L^{(s)*} + L^{(s)} ) x \rangle_2 : x \in \cX, \norm{x}_2 = 1, x \perp \rho^{1/2} \right\}. \label{eq-gap-L-s-equivalent-form}
        \end{align}
        To justify the second equality, choose graph-norm approximants in $\cX$ and subtract their $\rho^{1/2}$ components. The resulting approximants remain in $\cX$, are orthogonal to $\rho^{1/2}$, and have the same images under $L^{(s)}$. By Theorem~\ref{theorem-other-eigenvalues-L-star-plus-L}, there is an orthonormal basis $(f_k)_{k\in\N}$ of $\cB_2(\mathsf h)$ consisting of eigenvectors of $L^{(s)*}+L^{(s)}$, with corresponding eigenvalues $(\lambda_k)_{k\in\N}$. In particular, $\lambda_0=0$ and $f_0=\rho^{1/2}$. If $0$ were not simple, there would be an eigenvector $f_m\perp\rho^{1/2}$, with $m\geq1$, such that $(L^{(s)*}+L^{(s)})f_m=0$. Equation~\eqref{eq-gap-L-s-equivalent-form} would then give
        \begin{equation*}
            \gap L^{(s)} \le - \frac{1}{2} \left\langle f_m, (L^{(s)*} + L^{(s)}) f_m \right\rangle_2  = 0,
        \end{equation*}
        contradicting the assumption that $L^{(s)}$ has a spectral gap.

        Conversely, suppose that \eqref{item-spectral-gap-theorem-3} holds. Then $\lambda_k<0$ for every $k\geq1$. If $x\in\cX$ satisfies $\norm{x}_2=1$ and $x\perp\rho^{1/2}$, write $x=\sum_{k\geq0}x_kf_k$, where $\sum_{k\geq0}\abs{x_k}^2=1$ and $x_0=0$. Then
        \begin{equation} \label{eq-spectral-gap-by-eigenvalues}
            \left\langle x, ( L^{(s)*} + L^{(s)} ) x \right\rangle_2 = \left\langle \sum_{j \ge 0} x_j f_j, \left(L^{(s)*} + L^{(s)}\right) \sum_{k \ge 0} x_k f_k \right\rangle_2 = \sum_{k \ge 1} \lambda_{k} \abs{x_k}^2.
        \end{equation}
        Each $\lambda_k$, $k\geq1$, has the form $n\lambda+m\mu$ with $n+m\geq1$. Since $\mu\leq\lambda<0$,
        \[
            n\lambda+m\mu\leq(n+m)\lambda\leq\lambda<0.
        \]
        Equations~\eqref{eq-gap-L-s-equivalent-form} and \eqref{eq-spectral-gap-by-eigenvalues} therefore imply
        \begin{equation*}
            \gap L^{(s)} \ge - \lambda / 2 > 0.
        \end{equation*}
        Thus \eqref{item-spectral-gap-theorem-3} implies \eqref{item-spectral-gap-theorem-1}.

        Assume \eqref{item-spectral-gap-theorem-3}. Both fundamental eigenvalues are then strictly negative. Theorem~\ref{theorem-other-eigenvalues-L-star-plus-L} shows that every eigenspace is finite-dimensional and that the eigenvalues tend to $-\infty$. The resolvent is therefore the operator-norm limit of finite-rank spectral truncations and is compact; see \cite[Proposition~5.12]{schmudgen2012unbounded}. This proves \eqref{item-spectral-gap-theorem-4}.

        Finally, assume \eqref{item-spectral-gap-theorem-4}. If $0$ were not simple, one fundamental eigenvalue would vanish. Theorem~\ref{theorem-other-eigenvalues-L-star-plus-L} would then give an infinite-dimensional eigenspace, which is incompatible with compact resolvent. Hence \eqref{item-spectral-gap-theorem-4} implies \eqref{item-spectral-gap-theorem-3}.
    \end{proof}

    \begin{remark} \label{remark-exact-spectral-gap}
        Equations~\eqref{eq-gap-L-s-equivalent-form} and \eqref{eq-spectral-gap-by-eigenvalues} show that the spectral gap of $L^{(s)}$ is determined by the fundamental eigenvalue of $L^{(s)*}+L^{(s)}$ closest to zero. Thus \eqref{eq-L-star-plus-L-fundamental-eigenvalues} gives
        \begin{equation*}
            s \mapsto \gap L^{(s)} = \gamma - \abs{\kappa} \cosh \left( (1-2s)\beta / 2 \right).
        \end{equation*}
        This function is increasing on $[0,1/2]$ and decreasing on $[1/2,1]$. At the KMS embedding $s=1/2$,
        \begin{equation*}
            \gap L^{(1/2)} = \gamma - \abs{\kappa} > 0
        \end{equation*}
        by \eqref{eq-kms-spectral-gap-strict-positive}. Moreover, \eqref{eq-gns-spectral-gap-positive} gives $ \abs{\kappa}\cosh(\beta/2)\leq\gamma $. If $\kappa\neq0$ and $0<s<1$, then $\abs{1-2s}<1$, and therefore
        \[
            \abs{\kappa}\cosh\!\left(\frac{1-2s}{2}\,\beta\right)
            <\abs{\kappa}\cosh(\beta/2)\leq\gamma.
        \]
        If $\kappa=0$, the spectral gap equals $\gamma>0$. Hence the spectral gap is maximal at $s=1/2$, is strictly positive for every $0<s<1$, and may vanish only at the GNS and anti-GNS endpoints $s=0,1$.

    \end{remark}

    We next transfer compact resolvent from the self-adjoint comparison operator to the induced generator. The same argument applies to the adjoint generator.

    \begin{proposition} \label{theorem-compact-resolvent-spectral-gap}
        If $L^{(s)}$ has a spectral gap, then both $L^{(s)}$ and $L^{(s)*}$ have compact resolvent.
    \end{proposition}
    \begin{proof}
        Fix $\lambda>0$. Since all three operators generate contraction semigroups or are nonpositive self-adjoint, $\lambda$ belongs to the resolvent sets of $L^{(s)}$, $L^{(s)*}$, and $L^{(s)*}+L^{(s)}$.

        By Lemma~\ref{lemma-L-star-plus-L-essential-sa}, $\cX$ is an operator core, and hence a form core, for the self-adjoint operator $L^{(s)*}+L^{(s)}$; see \cite[Problem~16, p.\,314]{reed1980methods}. It is consequently a form core for $\lambda\1-(L^{(s)*}+L^{(s)})$.

        Lemma~\ref{lemma-cX-common-core} shows that $\cX$ is also a core for $\lambda\1-L^{(s)}$. By \cite[Theorem~X.25 and its corollary]{reed1975methods}, it is therefore a form core for the positive self-adjoint operator
        \begin{equation*}
            (\lambda \1 - L^{(s)*})(\lambda \1 - L^{(s)}) =: \abs{\lambda \1 - L^{(s)}}^2.
        \end{equation*}

        For every $x\in\cX$,
        \begin{align}
            \left\langle x, \abs{\lambda \1 - L^{(s)}}^2 x \right\rangle_2  &= \left\langle x, (\lambda^2 \1 - \lambda (L^{(s)*} + L^{(s)}) + L^{(s)*} L^{(s)}) x  \right\rangle_2 \nonumber \\
            &= \left\langle x, \lambda (\lambda \1 - (L^{(s)*} + L^{(s)})) x \right\rangle_2 + \left\langle L^{(s)} x, L^{(s)} x \right\rangle_2 \nonumber \\
            &\ge \left\langle x, \lambda (\lambda \1 - (L^{(s)*} + L^{(s)})) x \right\rangle_2 \ge \lambda^2 \left\langle x,  x \right\rangle_2. \label{eq-form-norms-dominance}
        \end{align}
        Inequality~\eqref{eq-form-norms-dominance} shows that the form norm of $\lambda(\lambda\1-(L^{(s)*}+L^{(s)}))$ is controlled by that of $\abs{\lambda\1-L^{(s)}}^2$. Since $\cX$ is a common form core,
        \begin{equation*}
            \dom \abs{\lambda \1 - L^{(s)}} \subset \dom \left( \lambda \1 - (L^{(s)*} + L^{(s)}) \right)^{1/2},
        \end{equation*}
        and therefore,
        \begin{equation*}
            \abs{\lambda \1 - L^{(s)}}^2 \ge \lambda (\lambda \1 - (L^{(s)*} + L^{(s)})).
        \end{equation*}
        The order relation for inverses \cite[Corollary~10.12]{schmudgen2012unbounded} now gives
        \begin{equation} \label{eq-crucial-operator-inequality-upper-bound-of-resolvent-of-L}
            \abs{\lambda \1 - L^{(s)}}^{-2} \le \lambda^{-1} \left( \lambda \1 - (L^{(s)*} + L^{(s)}) \right)^{-1}.
        \end{equation}

        Because $L^{(s)}$ has a spectral gap, Theorem~\ref{theorem-spectral-gap-compact-resolvent} makes the operator on the right-hand side compact. Inequality~\eqref{eq-crucial-operator-inequality-upper-bound-of-resolvent-of-L} then implies that $\abs{\lambda\1-L^{(s)}}^{-2}$ is compact. Since
        \[
            \abs{(\lambda\1-L^{(s)*})^{-1}}^2
            =\abs{\lambda\1-L^{(s)}}^{-2},
        \]
        the positive square root $\abs{(\lambda\1-L^{(s)*})^{-1}}$ is compact. The polar decomposition shows that $(\lambda\1-L^{(s)*})^{-1}$ is compact, and its adjoint $(\lambda\1-L^{(s)})^{-1}$ is compact as well.
    \end{proof}

    Proposition~\ref{theorem-compact-resolvent-spectral-gap} shows that a spectral gap implies compact resolvent for both $L^{(s)}$ and $L^{(s)*}$. Consequently, in that case, their spectra consist entirely of the discrete eigenvalues identified in Theorem~\ref{theorem-other-eigenvalues-L} and Corollary~\ref{theorem-adjoint-generator-eigenvalues}. At an endpoint embedding where the spectral gap vanishes, this argument gives no information about compactness. Section~\ref{section-immediate-compactness} resolves that case by a quantum characteristic-function method.

    % ===================================================================

    % NEW SECTION

    \section{Characteristic Functions, Regularization, and Compactness} \label{section-immediate-compactness}

    This section completes the spectral analysis of $L^{(s)}$. Section~\ref{section-spectrum-L-star-plus-L} proves compactness of the resolvent when the spectral gap is strictly positive. Since the spectral gap may vanish at the GNS or anti-GNS embedding, we now prove directly that $(T_t^{(s)})_{t\geq0}$ is immediately compact for every $s\in[0,1]$.

    Quantum characteristic functions provide the representation used in this section. The quantum Plancherel theorem \cite[Lemma 2.7]{keyl2016schwartz} identifies $\cB_2(\mathsf h)$ unitarily with $L^2(\C)$ through normalized quantum characteristic functions. After this identification, the induced semigroup is a complex Gaussian integral operator. We compute its kernel on a natural dense subspace, prove that the formula extends to all of $L^2(\C)$, and show that the kernel is square-integrable for every $t>0$. It follows that $T_t^{(s)}$ is a Hilbert--Schmidt operator on $\cB_2(\mathsf h)$. We also establish the stronger regularization property
    \begin{equation*}
        T_t^{(s)}\big(\cB_2(\mathsf h)\big)\subset\cH,
        \qquad t>0,
    \end{equation*}
    where the Sobolev-type space $\cH$ is compactly embedded in $\cB_2(\mathsf h)$.

    Immediate compactness then shows that the polynomial eigenvalues obtained in Section~\ref{section-eigenvalues-of-L} exhaust the spectra of $L^{(s)}$ and $L^{(s)*}$ for every $s\in[0,1]$.

    \subsection{Characteristic-Function Representation of the Induced Semigroup}

    Throughout this section, we use the phase-space convention from Appendix~\ref{section-appendix-identification} and identify $\C$ with $\R^2$ by
    \begin{equation*}
        z=x+\mi y \longleftrightarrow \mathbf{z}:=\begin{pmatrix}x\\y\end{pmatrix},
        \qquad
        \mathbf J:=\begin{pmatrix}0&1\\-1&0\end{pmatrix}.
    \end{equation*}
    Then $\mathbf z^T\mathbf J\mathbf w=\Im\langle z,w\rangle$. We write $\md\mathbf z=\md x\,\md y$ for Lebesgue measure. Accordingly, $L^2(\C)$ means $L^2(\C,\md\mathbf z)$, identified with $L^2(\R^2,\md x\,\md y)$, and $\cS(\C)$ denotes the Schwartz space. Unless limits are displayed, all phase-space integrals are taken over $\C\simeq\R^2$. As above, $\cE$ is the dense linear span of the exponential vectors in $\mathsf h$.

    We recall that the induced semigroup is explicitly given on the dense subspace
    \[
        \rho^{s/2}\cB(\mathsf h)\rho^{(1-s)/2}
        \subset\cB_2(\mathsf h)
    \]
    by
    \begin{equation*}
        T_t^{(s)} ( \rho^{s/2} X \rho^{(1-s)/2} ) = \rho^{s/2} \cT_t (X) \rho^{(1-s)/2},
    \end{equation*}
    for $X\in\cB(\mathsf h)$. Although the strongly continuous contraction semigroup $(T_t^{(s)})_{t\geq0}$ acts on all of $\cB_2(\mathsf h)$, this formula does not describe its action explicitly outside that dense subspace.

    Every operator $\rho^{s/2}X\rho^{(1-s)/2}$ in this dense subspace is trace class for $s\in[0,1]$.

    For $\psi\in\cB_1(\mathsf h)$, its quantum characteristic function $\chi_\psi$ is defined in \eqref{eq-introduction-characteristic-function}.
    The quantum Plancherel theorem \cite[Section 5.3]{holevo2011probabilistic} gives
    \begin{equation} \label{eq-quantum-plancherel}
        \langle \varphi,\psi\rangle_2
        =\frac1\pi\int\overline{\chi_\varphi(z)}\chi_\psi(z)\,\md\mathbf z
    \end{equation}
    for $\varphi,\psi\in\cB_1(\mathsf h)$. Consequently, there is a unique unitary operator
    \begin{equation} \label{eq-unitary-characteristic-function-transform}
        \cF_W:\cB_2(\mathsf h)\longrightarrow L^2(\C,\md\mathbf z),
        \qquad
        \cF_W\psi:=\pi^{-1/2}\chi_\psi
        \quad\bigl(\psi\in\cB_1(\mathsf h)\bigr).
    \end{equation}
    See also \cite[Lemma 2.7]{keyl2016schwartz}. We call $\cF_W$ the normalized quantum Fourier--Plancherel transform. For $f\in\cS(\C)$, its inverse is
    \begin{equation*}
        \cF_W^{-1}f
        =\frac{1}{\sqrt{\pi}}\int f(z)W(-z)\,\md\mathbf z.
    \end{equation*}
    Thus, whenever $\chi_Y\in\cS(\C)$, Holevo's inversion formula \cite[Corollary~5.3.5]{holevo2011probabilistic}, written for the unnormalized quantum characteristic function, becomes
    \begin{equation*}
        Y=\frac{1}{\pi}\int\chi_Y(z)W(-z)\,\md\mathbf z,
    \end{equation*}
    where the integral is understood weakly. For a general Hilbert--Schmidt operator $Y$, the notation $\chi_Y$ denotes the $L^2$-equivalence class supplied by \eqref{eq-unitary-characteristic-function-transform}; it need not have a pointwise trace representation.

    We will prove that
    \begin{equation*}
        \cF_W T_t^{(s)}\cF_W^{-1}:L^2(\C)\longrightarrow L^2(\C)
    \end{equation*}
    is a complex Gaussian integral operator. This representation describes $T_t^{(s)}$ on the entire Hilbert--Schmidt space. When $s=1/2$ or $\kappa=0$, it reduces to a real Gaussian convolution followed by multiplication by a strictly decaying real Gaussian.

    Recall that $\mathbf S=\coth(\beta/2)\1$ is the covariance matrix of the invariant density operator $\rho$, and set $\mathbf\Sigma:=\coth(\beta/4)\1$ and $\mathbf S_t:=\int_0^t\me^{r\mathbf Z^T}\mathbf C\me^{r\mathbf Z}\,\md r$. Define
    \begin{equation*}
        \mathbf{Q}_t := \mathbf{S}_t + \me^{t \mathbf{Z}^T} \mathbf{\Sigma} \me^{t \mathbf{Z}}
        = \mathbf{S} + \me^{t \mathbf{Z}^T} ( \mathbf{\Sigma} - \mathbf{S} ) \me^{t \mathbf{Z}}.
    \end{equation*}
    The second equality follows by integrating the invariant-state Lyapunov equation $\mathbf Z^T\mathbf S+\mathbf S\mathbf Z+\mathbf C=0$ from \eqref{eq-lyapunov-equations}.
    Moreover, $\sqrt{\tanh(\beta/4)}\,\rho^{1/2}$ is a normalized Gaussian state with covariance matrix $\mathbf\Sigma$. Hence $\mathbf Q_t$ is the covariance matrix of the Gaussian state $\sqrt{\tanh(\beta/4)}\,\cT_{*t}(\rho^{1/2})$.

    \begin{lemma} \label{lemma-weighted-weyl-pairing}
        For all $s\in[0,1]$ and $ z, w \in \C $, we have
        \begin{equation} \label{eq-innerproduct-two-Weyl-operators}
            \Tr ( \rho^{\frac{s}{2}} W(z) \rho^{\frac{1-s}{2}} W(w) ) = \sqrt{\coth(\beta/4)} \exp\!\left\{-\frac{1}{2} \mathbf{z}^T \mathbf{\Sigma} \mathbf{z} - \frac{1}{2} \mathbf{w}^T \mathbf{\Sigma} \mathbf{w} - \mathbf{z}^T \csch(\beta/4) \mathbf{\Theta}_s \mathbf{w}\right\},
        \end{equation}
        where
        \begin{equation} \label{definition-Theta-s}
            \mathbf{\Theta}_s := \cosh\!\left(\frac{1-2s}{4}\,\beta\right)\1
            -\sinh\!\left(\frac{1-2s}{4}\,\beta\right)\mi\mathbf J.
        \end{equation}
    \end{lemma}
    \begin{proof}
        The calculation is given in Appendix~\ref{section-characteristic-function-calculations}.
    \end{proof}
    \gdef\proofweightedweylpairs{%
    \begin{proof}[Proof of Lemma~\ref{lemma-weighted-weyl-pairing}]
        Using Parthasarathy's resolution of the identity in terms of coherent states \cite[Proposition 2.11]{parthasarathy2010gaussian}, for $ z, w \in \C $ we have
        \begin{align*}
            &\quad \Tr(\rho^{\frac{s}{2}} W(z) \rho^{\frac{1-s}{2}} W(w) ) \\
            &= \frac{1}{\pi} \int \me^{- \abs{\xi}^2} \langle e(\xi), \rho^{\frac{s}{2}} W(z) \rho^{\frac{1-s}{2}} W(w) e(\xi) \rangle \md \xi \\
            &= \frac{(1 - \me^{-\beta})^{\frac{1}{2}}}{\pi} \int \me^{- \abs{\xi}^2} \langle W(-z) \me^{-\frac{s}{2} \beta N} e(\xi), \me^{-\frac{1-s}{2} \beta N} W(w) e(\xi) \rangle \md \xi \\
            &= \frac{(1 - \me^{-\beta})^{\frac{1}{2}}}{\pi} \me^{ - \frac{1}{2} (\abs{w}^2 + \abs{z}^2)  - \me^{- \frac{1-s}{2} \beta} \langle z, w \rangle } \int \exp\!\left\{ - (1-\me^{-\frac{\beta}{2}}) \abs{\xi}^2\right\} \\
            &\quad \cdot \exp\!\left\{ - \me^{-\frac{1-s}{2} \beta} \langle z, \xi \rangle + \me^{-\frac{1}{2} \beta} \langle \xi, w \rangle + \me^{-\frac{s}{2} \beta} \langle \xi, z \rangle - \langle w, \xi \rangle \right\} \md \xi \\
            &= \frac{(1 - \me^{-\beta})^{\frac{1}{2}}}{\pi}  \me^{ - \frac{1}{2} (\abs{w}^2 + \abs{z}^2)  - \me^{- \frac{1-s}{2} \beta} \langle z, w \rangle } \int \exp\!\left\{ - (1 -\me^{-\frac{\beta}{2}}) \boldsymbol{\xi}^T \boldsymbol{\xi}\right\} \\
            &\quad \cdot \exp\!\left\{ ( - (\me^{-\frac{1-s}{2} \beta} \mathbf{z}^T + \mathbf{w}^T) (\1 + \mi \mathbf{J}) + ( \me^{-\frac{1}{2} \beta} \mathbf{w}^T + \me^{-\frac{s}{2} \beta} \mathbf{z}^T ) (\1 - \mi \mathbf{J}) ) \boldsymbol{\xi} \right\} \md \boldsymbol{\xi} \\
            &= \sqrt{\coth(\beta/4)}  \exp\!\left\{ - \frac{1}{2} (\abs{w}^2 + \abs{z}^2)  - \me^{- \frac{1-s}{2} \beta} \langle z, w \rangle \right\} \\
            &\quad \cdot \exp\!\left\{ \frac{ \langle - \me^{-\frac{1-s}{2} \beta} z - w, \me^{-\frac{1}{2} \beta} w + \me^{-\frac{s}{2} \beta} z \rangle }{ (1 - \me^{-\beta/2})} \right\} \\
            &= \sqrt{\coth(\beta/4)} \exp\!\left\{ - \left( \frac{1}{2} + \frac{\me^{-\frac{\beta}{2}}}{1 - \me^{-\frac{\beta}{2}}} \right) \abs{z}^2 - \left( \frac{1}{2} + \frac{\me^{-\frac{\beta}{2}}}{1 - \me^{-\frac{\beta}{2}}} \right) \abs{w}^2 \right. \\
            &\qquad \left. - \me^{-\frac{1-s}{2} \beta} \langle z, w \rangle - \frac{\me^{-\frac{1-s}{2} \beta} \me^{- \frac{1}{2} \beta} }{1 - \me^{-\frac{\beta}{2}}} \langle z, w \rangle - \frac{\me^{-\frac{s}{2} \beta}}{1 - \me^{-\frac{\beta}{2}}} \langle w, z \rangle \right\} \\
            &= \sqrt{\coth(\beta/4)} \exp\!\left\{ - \frac{1}{2} \mathbf{z}^T \mathbf{\Sigma} \mathbf{z} - \frac{1}{2} \mathbf{w}^T \mathbf{\Sigma} \mathbf{w} \right. \\
            &\qquad \left. - \csch(\beta/4) \mathbf{z}^T \left( \cosh(\frac{1-2s}{4} \beta) \1 - \sinh(\frac{1-2s}{4} \beta) \mi \mathbf{J} \right)  \mathbf{w} \right\}.
        \end{align*}
        By \eqref{definition-Theta-s}, the last expression is precisely \eqref{eq-innerproduct-two-Weyl-operators}.
    \end{proof}}

    Since $(\mi\mathbf J)^2=\1$, definition \eqref{definition-Theta-s} is equivalently written as
    \[
            \mathbf\Theta_s=\exp\!\left\{\frac{2s-1}{4}\,\beta\mi\mathbf J\right\}.
    \]
    Thus $\mathbf\Theta_s$ is Hermitian and positive definite, as well as complex orthogonal in the sense that
    \begin{equation} \label{eq-Theta-basic-identities}
        \mathbf\Theta_s^T=\mathbf\Theta_{1-s}=\mathbf\Theta_s^{-1},
        \qquad
        \mathbf\Theta_s^T\mathbf\Theta_s
        =\mathbf\Theta_s\mathbf\Theta_s^T=\1.
    \end{equation}
    At the three distinguished embeddings,
    \begin{equation*}
        \mathbf\Theta_0=\cosh(\beta/4)\1-\sinh(\beta/4)\mi\mathbf J,
        \quad \mathbf\Theta_{1/2}=\1,
        \quad \mathbf\Theta_1=\cosh(\beta/4)\1+\sinh(\beta/4)\mi\mathbf J.
    \end{equation*}

    \begin{corollary} \label{corollary-characteristic-function-predual-wrapped-Wz}
        The characteristic function of $ \cT_{*t} (\rho^{(1-s)/2} W(z) \rho^{s/2}) $ is
        \begin{equation*}
            \chi_{ \cT_{*t} (\rho^{(1-s)/2} W(z) \rho^{s/2}) } (y) = \sqrt{\coth(\beta/4)} \, \me^{ -\frac{1}{2} \mathbf{z}^T \mathbf{\Sigma} \mathbf{z} - \frac{1}{2} \mathbf{y}^T \mathbf{Q}_t \mathbf{y} - \mathbf{z}^T \csch(\beta/4) \mathbf{\Theta}_{1-s} \me^{t \mathbf{Z}} \mathbf{y} }.
        \end{equation*}
    \end{corollary}

    We use the term \emph{Schwartz operator} in the sense of \cite[Definition~3.1]{keyl2016schwartz}. Equivalently, an operator is Schwartz precisely when its quantum characteristic function is a Schwartz function; see \cite[Proposition~3.18]{keyl2016schwartz}. The following theorem provides the Gaussian kernel needed below.

    \begin{theorem} \label{theorem-characteristic-function}
        Let $t>0$ and $z\in\C$. There exists a unique Schwartz operator $K_{t,z}$ such that
        \begin{equation} \label{eq-sandwiched-K-defining-identity}
            \rho^{(1-s)/2}K_{t,z}\rho^{s/2}
            =\cT_{*t}\!\left(\rho^{(1-s)/2}W(z)\rho^{s/2}\right).
        \end{equation}
        In particular, $K_{t,z}$ is trace class. Its characteristic function is the Gaussian
        \begin{align}
            \chi_{K_{t,z}}(\xi)
            &= c_t \, \exp\!\left\{ - \frac{1}{2} \boldsymbol{\xi}^T \mathbf{A}_t \boldsymbol{\xi} + \mathbf{z}^T \mathbf{B}_t \boldsymbol{\xi} - \frac{1}{2} \mathbf{z}^T \mathbf{C}_t \mathbf{z} \right\} \label{eq-characteristic-function-final-result} \\
            &= c_t \, \exp\!\left\{ -\frac{1}{2} \begin{pmatrix}
            \mathbf{z}^T & \boldsymbol{\xi}^T
            \end{pmatrix} \mathbf{G}_t \begin{pmatrix}
            \mathbf{z} \\ \boldsymbol{\xi}
            \end{pmatrix}\right\} \label{eq-characteristic-function-final-result-equivalent-form}
        \end{align}
        where
        \begin{align}
            \mathbf{A}_t &= \coth(\beta/4)\mathbf{\Theta}_{1-s} (\1 + \me^{t \mathbf{Z}^T} \me^{t \mathbf{Z}} ) (\1 - \me^{t \mathbf{Z}^T} \me^{t \mathbf{Z}})^{-1} \mathbf{\Theta}_{1-s}^T, \label{eq-A-t-by-M-t} \\
            \mathbf{B}_t &= 2\coth(\beta/4)\mathbf{\Theta}_{1-s} ( \me^{t \mathbf{Z}^T} - \me^{- t \mathbf{Z}} )^{-1} \mathbf{\Theta}_{1-s}^T, \label{eq-B-t-by-M-t} \\
            \mathbf{C}_t &= \coth(\beta/4)\mathbf{\Theta}_{1-s} (\1 + \me^{t \mathbf{Z}} \me^{t \mathbf{Z}^T}) (\1 - \me^{t \mathbf{Z}} \me^{t \mathbf{Z}^T})^{-1} \mathbf{\Theta}_{1-s}^T, \label{eq-C-t-by-M-t} \\
            c_t &= \coth(\beta/4) \left( \det (\1 - \me^{t \mathbf{Z}^T} \me^{t \mathbf{Z}} ) \right)^{-1/2}, \label{eq-final-c-t}
        \end{align}
        and
        \begin{equation} \label{eq-definition-G-t}
            \mathbf{G}_t := \begin{pmatrix}
            \mathbf{C}_t & - \mathbf{B}_t \\
            - \mathbf{B}_t^T & \mathbf{A}_t
            \end{pmatrix}.
        \end{equation}
        Equivalently, $\mathbf G_t$ is the Cayley transform
        \begin{equation} \label{eq-G-t-Cayley-transform}
        \begin{split}
            \mathbf G_t=\coth(\beta/4)&\left(
            \1+\begin{pmatrix}
                0 & \mathbf\Theta_{1-s}\me^{t\mathbf Z}\mathbf\Theta_{1-s}^T\\
                \mathbf\Theta_{1-s}\me^{t\mathbf Z^T}\mathbf\Theta_{1-s}^T & 0
            \end{pmatrix}\right)\\
            &\times\left(
            \1-\begin{pmatrix}
                0 & \mathbf\Theta_{1-s}\me^{t\mathbf Z}\mathbf\Theta_{1-s}^T\\
                \mathbf\Theta_{1-s}\me^{t\mathbf Z^T}\mathbf\Theta_{1-s}^T & 0
            \end{pmatrix}\right)^{-1}.
        \end{split}
        \end{equation}
        Here the quadratic variables are real, and $\Re\mathbf G_t$ denotes the entrywise real part of the complex symmetric matrix $\mathbf G_t$. Moreover, $\Re\mathbf G_t>0$. Consequently, the formal product
        \begin{equation*}
            \rho^{-(1-s)/2}
            \cT_{*t}\!\left(\rho^{(1-s)/2}W(z)\rho^{s/2}\right)
            \rho^{-s/2}
        \end{equation*}
        is well defined on $\dom\rho^{-s/2}$ and extends uniquely to $K_{t,z}$.
    \end{theorem}

    \begin{proof}
        The proof is given in Appendix~\ref{section-characteristic-function-calculations}.
    \end{proof}
    \gdef\proofcharacteristickernel{%
    \begin{proof}[Proof of Theorem~\ref{theorem-characteristic-function}]
        Since the inverse powers of $\rho$ are unbounded, we do not assume at the outset that the formal product in the statement defines an operator. Instead, on $\cE\times\cE$ consider the sesquilinear form
        \begin{equation} \label{eq-sesquilinear-form-on-exponential-vectors}
            \mathfrak{t}_{t,z}(u,v)
            :=\left\langle \rho^{-\frac{1-s}{2}}u,
            \cT_{*t}\!\left(\rho^{\frac{1-s}{2}}W(z)\rho^{\frac{s}{2}}\right)
            \rho^{-\frac{s}{2}}v\right\rangle .
        \end{equation}
        Every exponential vector belongs to the domains of all real powers of $\rho$, so the form is well-defined. Motivated by the coherent-state resolution of the identity, we compute the following candidate characteristic function:
        \begin{align}
            &\quad \frac{1}{\pi} \int \me^{-\abs{\eta}^2}
            \mathfrak{t}_{t,z}(e(\eta),W(\xi)e(\eta))\,\md \eta \nonumber \\
            &= \frac{\me^{-\frac{1}{2} \abs{\xi}^2}}{\pi (1 - \me^{-\beta})^{\frac{1}{2}} } \int \me^{-\abs{\eta}^2 - \langle \xi, \eta \rangle } \langle e( \me^{\frac{1-s}{2} \beta} \eta ), \cT_{*t} ( \rho^{\frac{1-s}{2}} W(z) \rho^{\frac{s}{2}} ) e( \me^{\frac{s}{2} \beta} (\eta + \xi) )  \rangle \md \eta \label{eq-characteristic-function-predual-and-wrapped-Wz-first-step}
        \end{align}
        We introduce the matrix $ \mathbf{J}_s $ appearing in the first Gaussian integral:
        \begin{equation} \label{eq-definition-J-s}
            \mathbf{J}_s := \me^{\frac{s}{2} \beta} (\1 + \mi \mathbf{J}) - \me^{\frac{1-s}{2} \beta} (\1 - \mi \mathbf{J}).
        \end{equation}
        Note that $ \mathbf{J}_s = (\me^{\beta/2} - 1) \mi \mathbf{J} \csch(\beta/4) \mathbf{\Theta}_{s} $.

        Applying the Fourier inversion formula \cite[Corollary 5.3.5]{holevo2011probabilistic} together with Corollary \ref{corollary-characteristic-function-predual-wrapped-Wz}, the first Gaussian integration, with respect to $y$, gives
        \begin{align}
            &\quad \langle e( \me^{\frac{1-s}{2} \beta} \eta ), \cT_{*t} ( \rho^{\frac{1-s}{2}} W(z) \rho^{\frac{s}{2}} ) e(\me^{\frac{s}{2} \beta} (\eta + \xi) ) \rangle \nonumber \\
            &= \frac{1}{\pi} \int \chi_{\cT_{*t} (\rho^{(1-s)/2} W(z) \rho^{s/2}) } (y) \langle e(\me^{\frac{1-s}{2} \beta} \eta), W(-y) e( \me^{\frac{s}{2} \beta} (\eta + \xi) ) \rangle \md y \nonumber \\
            &= \frac{\sqrt{\coth(\beta/4)}}{\pi} \me^{-\frac{1}{2} \mathbf{z}^T \mathbf{\Sigma} \mathbf{z} + \me^{\frac{\beta}{2}} \boldsymbol{\eta}^T \boldsymbol{\eta} + \me^{\frac{\beta}{2}} \boldsymbol{\xi}^T (\1 - \mi \mathbf{J}) \boldsymbol{\eta}} \int \exp\!\left\{ -\frac{1}{2} \mathbf{y}^T ( \mathbf{Q}_t + \1 ) \mathbf{y} \right\} \nonumber \\
            &\quad \cdot \exp\!\left\{ ( - \mathbf{z}^T \csch(\beta/4) \mathbf{\Theta}_{1-s} \me^{t \mathbf{Z}} + \me^{\frac{s}{2} \beta} (\boldsymbol{\eta} + \boldsymbol{\xi})^T (\1 - \mi \mathbf{J}) - \me^{\frac{1-s}{2} \beta} \boldsymbol{\eta}^T (\1 + \mi \mathbf{J}) ) \mathbf{y} \right\}  \md \mathbf{y} \nonumber \\
            &= \frac{2 \sqrt{\coth(\beta/4)}}{\sqrt{\det (\mathbf{Q}_t + \1) }} \exp\!\left\{-\frac{1}{2} \mathbf{z}^T \mathbf{\Sigma} \mathbf{z} \right\} \nonumber \\
            &\quad \cdot \me^{\frac{1}{2} (- \mathbf{z}^T \csch(\beta/4) \mathbf{\Theta}_{1-s} \me^{t \mathbf{Z}} + \me^{\frac{s}{2} \beta} \boldsymbol{\xi}^T (\1 - \mi \mathbf{J}) ) (\mathbf{Q}_t + \1)^{-1} ( - \me^{t \mathbf{Z}^T} \csch(\beta/4) \mathbf{\Theta}_{1-s}^T \mathbf{z} + \me^{\frac{s}{2} \beta} (\1 + \mi \mathbf{J}) \boldsymbol{\xi} ) } \nonumber \\
            &\quad \cdot \exp\!\left\{ \boldsymbol{\eta}^T \left( \me^{\frac{\beta}{2}} \1 + \frac{1}{2} \mathbf{J}_s^T (\mathbf{Q}_t + \1)^{-1} \mathbf{J}_s \right) \boldsymbol{\eta} \right\} \nonumber \\
            &\quad \cdot \exp\!\left\{ \left( \me^{\frac{\beta}{2}} \boldsymbol{\xi}^T (\1 - \mi \mathbf{J})+ ( - \mathbf{z}^T \csch(\beta/4) \mathbf{\Theta}_{1-s} \me^{t \mathbf{Z}} + \me^{\frac{s}{2} \beta} \boldsymbol{\xi}^T (\1 - \mi \mathbf{J}) ) (\mathbf{Q}_t + \1)^{-1} \mathbf{J}_s \right) \boldsymbol{\eta} \right\} \label{eq-characteristic-function-predual-and-wrapped-Wz-integrand}
        \end{align}
        For the Gaussian integration with respect to $\eta$, set
        \begin{equation} \label{eq-definition-R-t}
            \mathbf{R}_t := - 2 (  ( \me^{\frac{\beta}{2}} - 1) \1 + \frac{1}{2} \mathbf{J}_s^T (\mathbf{Q}_t + \1)^{-1} \mathbf{J}_s ) = \mathbf{J}_s^T ( \frac{1 - \me^{-\beta/2}}{2} \1 - (\mathbf{Q}_t + \1)^{-1} ) \mathbf{J}_s.
        \end{equation}
        Lemma \ref{lemma-strict-positivity-real-parts-R-t} proves that $ \Re \mathbf{R}_t > 0 $. Hence, by \cite[Appendix A, Theorem 1]{Folland1989harmonic}, the Gaussian integral in $\eta$ converges for $t>0$.

        Substituting \eqref{eq-characteristic-function-predual-and-wrapped-Wz-integrand} into \eqref{eq-characteristic-function-predual-and-wrapped-Wz-first-step} and using
        \begin{equation} \label{eq-identity-exp-s-beta-one-half-Theta}
            \begin{aligned}
            \mathbf{J}_0 &= \mathbf{\Theta}_{\frac{1}{2}-s} \mathbf{J}_s, \\
            \me^{\frac{s}{2} \beta} (\1 - \mi \mathbf{J})
            &= \mathbf{\Theta}_{\frac{1}{2} - s} (\1 - \mi \mathbf{J}).
            \end{aligned}
        \end{equation}
        yields
        \begin{align}
            &\quad \frac{1}{\pi} \int \me^{-\abs{\eta}^2} \langle \rho^{-(1-s)/2} e(\eta),  \cT_{*t} ( \rho^{(1-s)/2} W(z) \rho^{s/2} ) \rho^{-s/2} W(\xi) e(\eta) \rangle \md \eta \nonumber \\
            &= \frac{\me^{-\frac{1}{2} \abs{\xi}^2}}{\pi (1 - \me^{-\beta})^{\frac{1}{2}} } \int \me^{-\abs{\eta}^2 - \langle \xi, \eta \rangle } \langle e( \me^{\frac{1-s}{2} \beta} \eta ), \cT_{*t} ( \rho^{\frac{1-s}{2}} W(z) \rho^{\frac{s}{2}} ) e( \me^{\frac{s}{2} \beta} (\eta + \xi) )  \rangle \md \eta \nonumber \\
            &= \frac{2 \me^{- \frac{1}{2} \abs{\xi}^2}}{\pi (1- \me^{-\beta})^{\frac{1}{2}}} \sqrt{ \frac{\coth(\beta/4)}{\det ( \mathbf{Q}_t + \1 )} } \exp\!\left\{-\frac{1}{2} \mathbf{z}^T \mathbf{\Sigma} \mathbf{z}\right\} \nonumber \\
            &\quad \cdot \me^{\frac{1}{2} (- \mathbf{z}^T \csch(\beta/4) \mathbf{\Theta}_{1-s} \me^{t \mathbf{Z}} + \me^{\frac{s}{2} \beta} \boldsymbol{\xi}^T (\1 - \mi \mathbf{J}) ) (\mathbf{Q}_t + \1)^{-1} ( - \me^{t \mathbf{Z}^T} \csch(\beta/4) \mathbf{\Theta}_{1-s}^T \mathbf{z} + \me^{\frac{s}{2} \beta} (\1 + \mi \mathbf{J}) \boldsymbol{\xi})} \nonumber \\
            &\quad \cdot \int \exp\!\left\{- \frac{1}{2} \boldsymbol{\eta}^T ( - 2 (  ( \me^{\frac{\beta}{2}} - 1) \1 + \frac{1}{2} \mathbf{J}_s^T (\mathbf{Q}_t + \1)^{-1} \mathbf{J}_s ) ) \boldsymbol{\eta}\right\} \nonumber \\
            &\qquad \cdot \exp\!\left\{\left( - \boldsymbol{\xi}^T \mathbf{J}_0 + (- \mathbf{z}^T \csch(\beta/4) \mathbf{\Theta}_{1-s} \me^{t \mathbf{Z}} + \me^{\frac{s}{2} \beta} \boldsymbol{\xi}^T (\1 - \mi \mathbf{J}) ) (\mathbf{Q}_t + \1)^{-1} \mathbf{J}_s \right) \boldsymbol{\eta} \right\} \md \boldsymbol{\eta} \nonumber \\
            &= \frac{2 \me^{- \frac{1}{2} \abs{\xi}^2}}{\pi (1- \me^{-\beta})^{\frac{1}{2}}} \sqrt{ \frac{\coth(\beta/4)}{\det ( \mathbf{Q}_t + \1 )} } \exp\!\left\{-\frac{1}{2} \mathbf{z}^T \mathbf{\Sigma} \mathbf{z}\right\} \nonumber \\
            &\quad \cdot \me^{\frac{1}{2} (- \mathbf{z}^T \csch(\beta/4) \mathbf{\Theta}_{1-s} \me^{t \mathbf{Z}} + \me^{\frac{s}{2} \beta} \boldsymbol{\xi}^T (\1 - \mi \mathbf{J}) ) (\mathbf{Q}_t + \1)^{-1} ( - \me^{t \mathbf{Z}^T} \csch(\beta/4) \mathbf{\Theta}_{1-s}^T \mathbf{z} + \me^{\frac{s}{2} \beta} (\1 + \mi \mathbf{J}) \boldsymbol{\xi})} \nonumber \\
            &\quad \cdot \int \me^{- \frac{1}{2} \boldsymbol{\eta}^T \mathbf{R}_t \boldsymbol{\eta} - \left( \mathbf{z}^T \csch(\beta/4) \mathbf{\Theta}_{1-s} \me^{t \mathbf{Z}} + \mathbf{\Theta}_{\frac{1}{2} - s} (\mathbf{Q}_t + \mi \mathbf{J}) \right) (\mathbf{Q}_t + \1)^{-1} \mathbf{J}_s \boldsymbol{\eta} } \md \boldsymbol{\eta} \nonumber \\
            &= \frac{4 \sqrt{\coth(\beta/4)}}{(1-\me^{-\beta})^{\frac{1}{2}} \sqrt{\det(\mathbf{Q}_t + \1) \det \mathbf{R}_t}} \exp\!\left\{-\frac{1}{2} \boldsymbol{\xi}^T \boldsymbol{\xi} \right\} \exp\!\left\{-\frac{1}{2} \mathbf{z}^T \mathbf{\Sigma} \mathbf{z}\right\} \nonumber \\
            &\quad \cdot \me^{\frac{1}{2} (- \mathbf{z}^T \csch(\beta/4) \mathbf{\Theta}_{1-s} \me^{t \mathbf{Z}} + \me^{\frac{s}{2} \beta} \boldsymbol{\xi}^T (\1 - \mi \mathbf{J}) ) (\mathbf{Q}_t + \1)^{-1} ( - \me^{t \mathbf{Z}^T} \csch(\beta/4) \mathbf{\Theta}_{1-s}^T \mathbf{z} + \me^{\frac{s}{2} \beta} (\1 + \mi \mathbf{J}) \boldsymbol{\xi})} \nonumber \\
            &\quad \cdot \exp\!\left\{ \frac{1}{2} \left( \mathbf{z}^T \csch(\beta/4) \mathbf{\Theta}_{1-s} \me^{t \mathbf{Z}} + \boldsymbol{\xi}^T \mathbf{\Theta}_{\frac{1}{2} - s} ( \mathbf{Q}_t + \mi \mathbf{J} )  \right) \right. \nonumber \\
            &\qquad \left. \cdot \mathbf{J}_s \mathbf{R}_t^{-1} \mathbf{J}_s^T (\mathbf{Q}_t + \1)^{-1} \left( \me^{t \mathbf{Z}^T} \csch(\beta/4) \mathbf{\Theta}^T_{1-s} \mathbf{z} + (\mathbf{Q}_t - \mi \mathbf{J}) \mathbf{\Theta}_{\frac{1}{2}-s}^T \boldsymbol{\xi} \right) \right\} \nonumber \\
            &= c_t \cdot \exp\!\left\{ - \frac{1}{2} \boldsymbol{\xi}^T \mathbf{A}_t \boldsymbol{\xi} + \mathbf{z}^T \mathbf{B}_t \boldsymbol{\xi} - \frac{1}{2} \mathbf{z}^T \mathbf{C}_t \mathbf{z} \right\}, \label{eq-characteristic-function-predual-and-wrapped-Wz-gaussian-form}
        \end{align}
        where
        \begin{align}
            \mathbf{A}_t &:= \1 - \me^{s \beta} (\1 - \mi \mathbf{J}) ( \mathbf{Q}_t + \1 )^{-1} (\1 + \mi \mathbf{J}) - \mathbf{\Theta}_{\frac{1}{2} - s} (\mathbf{Q}_t + \mi \mathbf{J}) \nonumber \\
            &\quad \cdot (\mathbf{Q}_t + \1)^{-1} \mathbf{J}_s \mathbf{R}_t^{-1} \mathbf{J}_s^T ( \mathbf{Q}_t + \1 )^{-1} (\mathbf{Q}_t - \mi \mathbf{J}) \mathbf{\Theta}_{\frac{1}{2} - s}^T, \label{eq-original-A-t} \\
            \mathbf{B}_t &= \csch(\beta/4) \mathbf{\Theta}_{1 - s} \me^{t \mathbf{Z}} ( \mathbf{Q}_t + \1 )^{-1} \left( - \me^{\frac{s}{2} \beta} (\1 + \mi \mathbf{J}) \right. \nonumber \\
            &\quad \left. + \mathbf{J}_s \mathbf{R}_t^{-1} \mathbf{J}_s^T (\mathbf{Q}_t + \1)^{-1} (\mathbf{Q}_t - \mi \mathbf{J}) \mathbf{\Theta}_{\frac{1}{2}-s}^T \right), \label{eq-original-B-t} \\
            \mathbf{C}_t &= \mathbf{\Sigma} - \csch(\beta/4) \mathbf{\Theta}_{1-s} \me^{t \mathbf{Z}} ( \mathbf{Q}_t + \1 )^{-1} \me^{t \mathbf{Z}^T} \csch(\beta/4) \mathbf{\Theta}_{1-s}^T \nonumber \\
            &\quad - \csch(\beta/4) \mathbf{\Theta}_{1-s} \me^{t \mathbf{Z}} (\mathbf{Q}_t + \1)^{-1} \mathbf{J}_s \mathbf{R}_t^{-1} \mathbf{J}_s^T ( \mathbf{Q}_t + \1 )^{-1} \me^{t \mathbf{Z}^T} \csch(\beta/4) \mathbf{\Theta}_{1-s}^T, \label{eq-original-C-t} \\
            & c_t = \frac{4 \sqrt{\coth(\beta/4)}}{(1-\me^{-\beta})^{\frac{1}{2}} \sqrt{\det(\mathbf{Q}_t + \1) \det \mathbf{R}_t}}. \label{eq-original-c-t}
        \end{align}

        These quantities are simplified in Lemma \ref{lemma-simplification-A-t-B-t-C-t}.

        We have therefore shown that \eqref{eq-characteristic-function-predual-and-wrapped-Wz-first-step} equals the Gaussian in \eqref{eq-characteristic-function-predual-and-wrapped-Wz-gaussian-form}, namely \eqref{eq-characteristic-function-final-result}. Its equivalent block-matrix form is \eqref{eq-characteristic-function-final-result-equivalent-form}. The strict positivity of $ \Re \mathbf{A}_t $ and $ \Re \mathbf{G}_t $ is established in Lemma \ref{lemma-strict-positivity-real-parts-G-t}.

        We next identify this Gaussian with an operator. Since $\Re\mathbf A_t>0$, the function of $\xi$ in \eqref{eq-characteristic-function-predual-and-wrapped-Wz-gaussian-form} is Schwartz. By \cite[Lemma 3.6 and Proposition 3.18]{keyl2016schwartz}, it is the characteristic function of a unique Schwartz operator, denoted by $K_{t,z}$. In particular, $K_{t,z}$ is trace class. Since $\rho^{s/2}W(y)\rho^{(1-s)/2}$ is also trace class, the product formula and \eqref{eq-quantum-plancherel} give
        \begin{align}
            \chi_{\rho^{(1-s)/2} K_{t,z} \rho^{s/2} } ( y )
            &= \Tr ( \rho^{(1-s)/2} K_{t,z} \rho^{s/2} W(y) ) \nonumber \\
            &=\Tr(K_{t,z} \rho^{s/2} W(y) \rho^{(1-s)/2}) \nonumber \\
            &= \frac{1}{\pi} \int \chi_{K_{t,z}} (- \xi) \chi_{\rho^{s/2} W(y) \rho^{(1-s)/2}} (\xi) \, \md \xi. \label{eq-characteristic-function-sandwiched-K}
        \end{align}
        Set
        \begin{equation} \label{eq-definition-D-t}
            \mathbf{D}_{t} = \mathbf{A}_t + \mathbf{\Sigma}.
        \end{equation}
        Using \eqref{eq-characteristic-function-predual-and-wrapped-Wz-gaussian-form} and \eqref{eq-innerproduct-two-Weyl-operators}, we obtain
        \begin{align*}
            \chi_{\rho^{(1-s)/2} K_{t,z} \rho^{s/2} } ( y ) &=  \frac{1}{\pi} \int c_t \exp\!\left\{ - \frac{1}{2} \boldsymbol{\xi}^T \mathbf{A}_t \boldsymbol{\xi} - \mathbf{z}^T \mathbf{B}_t \boldsymbol{\xi} - \frac{1}{2} \mathbf{z}^T \mathbf{C}_t \mathbf{z} \right\} \\
            &\quad \cdot \sqrt{\coth(\beta/4)} \exp\!\left\{ -\frac{1}{2} \mathbf{y}^T \mathbf{\Sigma} \mathbf{y} - \frac{1}{2} \boldsymbol{\xi}^T \mathbf{\Sigma} \boldsymbol{\xi} - \mathbf{y}^T \csch(\beta/4) \mathbf{\Theta}_s \boldsymbol{\xi}\right\} \md \boldsymbol{\xi} \\
            &= \frac{2 c_t \sqrt{\coth(\beta/4)}}{\sqrt{ \det ( \mathbf{A}_t + \mathbf{\Sigma} ) }} \exp\!\left\{ - \frac{1}{2}\mathbf{z}^T ( \mathbf{C}_t - \mathbf{B}_t \mathbf{D}_t^{-1} \mathbf{B}_t^T ) \mathbf{z} \right\} \\
            &\quad \cdot \exp\!\left\{-\frac12\mathbf y^T
            \left(\mathbf\Sigma-\csch^2(\beta/4)
            \mathbf\Theta_s\mathbf D_t^{-1}\mathbf\Theta_s^T\right)\mathbf y\right\}\\
            &\quad \cdot \exp\!\left\{\mathbf z^T\mathbf B_t\mathbf D_t^{-1}
            \csch(\beta/4)\mathbf\Theta_s^T\mathbf y\right\}.
        \end{align*}
        Lemma \ref{lemma-simplification-characteristic-function-sandwiched-K} and Corollary \ref{corollary-characteristic-function-predual-wrapped-Wz} give
        \begin{equation*}
            \chi_{\rho^{(1-s)/2} K_{t,z} \rho^{s/2}} (y) = \chi_{\cT_{*t} (\rho^{(1-s)/2} W(z) \rho^{s/2}) } (y), \quad \forall y \in \C.
        \end{equation*}
        The injectivity of $\cF_W$ therefore yields
        \begin{equation*}
            \rho^{(1-s)/2} K_{t,z} \rho^{s/2}
            =\cT_{*t}(\rho^{(1-s)/2}W(z)\rho^{s/2}),
        \end{equation*}
        which proves \eqref{eq-sandwiched-K-defining-identity}. Since the powers of $\rho$ have dense range, this identity also proves the uniqueness of $K_{t,z}$.

        Finally, if $v\in\dom\rho^{-s/2}$ and $u=\rho^{-s/2}v$, then
        \begin{equation*}
            \cT_{*t}(\rho^{(1-s)/2}W(z)\rho^{s/2})u
            =\rho^{(1-s)/2}K_{t,z}v.
        \end{equation*}
        Hence the formal product in the statement agrees with $K_{t,z}$ on its dense domain and therefore extends uniquely to $K_{t,z}$.
    \end{proof}}

    The preceding theorem now yields an explicit representation of the induced semigroup on the entire Hilbert--Schmidt space.

    \begin{corollary} \label{corollary-characteristic-function-induced-semigroup}
        Let $t>0$ and $Y\in\cB_2(\mathsf h)$. Then
        \begin{equation} \label{eq-characteristic-function-evolution-of-semigroup-Y}
            \chi_{T_t^{(s)} (Y)} (\mathbf{z}) = \frac{c_t}{\pi} \int \chi_{Y} (- \boldsymbol{\xi}) \exp\!\left\{ -\frac{1}{2} \begin{pmatrix}
            \mathbf{z}^T & \boldsymbol{\xi}^T
            \end{pmatrix} \mathbf{G}_t \begin{pmatrix}
            \mathbf{z} \\ \boldsymbol{\xi}
            \end{pmatrix}\right\} \, \md \boldsymbol{\xi},
        \end{equation}
        where equality is understood in $L^2(\C,\md\mathbf z)$. Moreover,
        \begin{equation*}
            T_t^{(s)}\in\cB_2\!\left(\cB_2(\mathsf h)\right)
            \qquad (t>0).
        \end{equation*}
    \end{corollary}
    \begin{proof}
        Let $X\in\cB(\mathsf h)$. Duality, Theorem~\ref{theorem-characteristic-function}, and the product formula give
        \begin{align}
            &\quad \chi_{T_t^{(s)} (\rho^{s/2} X \rho^{(1-s)/2}) } (z) \nonumber \\
            &= \Tr( \cT_t (X) \rho^{(1-s)/2} W(z) \rho^{s/2} ) \nonumber \\
            &= \Tr( X \cT_{*t} (\rho^{(1-s)/2} W(z) \rho^{s/2} ) ) \nonumber \\
            &= \Tr (X \rho^{(1-s)/2} K_{t,z} \rho^{s/2} ) \nonumber\\
            &= \Tr (\rho^{s/2} X \rho^{(1-s)/2} K_{t,z} ) \nonumber\\
            &= \frac{c_t}{\pi} \int \chi_{\rho^{s/2} X \rho^{(1-s)/2}} (- \boldsymbol{\xi}) \cdot \exp\!\left\{ -\frac{1}{2} \begin{pmatrix}
            \mathbf{z}^T & \boldsymbol{\xi}^T
            \end{pmatrix} \mathbf{G}_t \begin{pmatrix}
            \mathbf{z} \\ \boldsymbol{\xi}
            \end{pmatrix}\right\} \, \md \boldsymbol{\xi},
            \label{eq-characteristic-function-induced-semigroup-sandwiched-X}
        \end{align}
        Since $\Re\mathbf G_t>0$, the Gaussian kernel is square-integrable and defines a Hilbert--Schmidt integral operator \cite[Theorem~VI.23]{reed1980methods}. Unitarity of $\cF_W$ and density of $\rho^{s/2}\cB(\mathsf h)\rho^{(1-s)/2}$ extend \eqref{eq-characteristic-function-induced-semigroup-sandwiched-X} to every $Y\in\cB_2(\mathsf h)$. Hence $T_t^{(s)}$ is Hilbert--Schmidt on $\cB_2(\mathsf h)$.
    \end{proof}

    \begin{remark}
        When $s=1/2$ or $\kappa=0$, the matrices $\mathbf A_t$, $\mathbf B_t$, and $\mathbf C_t$ are real. Completing the square in \eqref{eq-characteristic-function-final-result} and translating the integration variable reduce \eqref{eq-characteristic-function-evolution-of-semigroup-Y} to
        \begin{equation} \label{eq-characteristic-function-induced-semigroup-real-case}
            \chi_{T_t^{(s)} (Y)} (z) = \frac{c_t}{\pi} \, \me^{ -\frac{1}{2} \mathbf{z}^T (\mathbf{C}_t - \mathbf{B}_t \mathbf{A}_t^{-1} \mathbf{B}_t^T ) \mathbf{z} } \int \chi_{Y} (- \mathbf{A}_t^{-1} \mathbf{B}_t^T \mathbf{z} -\boldsymbol{\xi}) \me^{ - \frac{1}{2} \boldsymbol{\xi}^T \mathbf{A}_t \boldsymbol{\xi} } \, \md \boldsymbol{\xi},
        \end{equation}
        where
        \begin{align*}
            \mathbf{C}_t - \mathbf{B}_t \mathbf{A}_t^{-1} \mathbf{B}_t^T  &=  \coth(\beta/4) (\1 - \me^{t \mathbf{Z}} \me^{t \mathbf{Z}^T}) (\1 +\me^{t \mathbf{Z}} \me^{t \mathbf{Z}^T})^{-1},  \\
            \mathbf{A}_t^{-1} \mathbf{B}_t^T &= - 2  (\me^{t \mathbf{Z}} + \me^{- t \mathbf{Z}^T} )^{-1}.
        \end{align*}
        Lemma~\ref{lemma-strict-positivity-real-parts-G-t} also gives $\Re\mathbf A_t>0$ and $\mathbf C_t-\mathbf B_t\mathbf A_t^{-1}\mathbf B_t^T>0$. Formula~\eqref{eq-characteristic-function-induced-semigroup-real-case} does not extend to the remaining cases: in general, $-\mathbf A_t^{-1}\mathbf B_t^T\mathbf z$ belongs to $\C^2$, whereas $\chi_Y$ is defined on $\C$, identified with $\R^2$.

        Thus, when $s=1/2$ or $\kappa=0$, $T_t^{(s)}$ acts by real Gaussian convolution followed by multiplication by a strictly decaying real Gaussian.
    \end{remark}

    The next example gives the further simplification obtained when the drift matrix is a scalar multiple of the identity. In that case $\kappa=0$, and the phase-space evolution is particularly transparent.

    \begin{example}
        Consider the Gaussian generator
        \begin{align}
            \cL (x) &= \frac{\me^\beta}{\me^\beta - 1} \left( - \frac{1}{2} a^\dagger a x + a^\dagger x a - \frac{1}{2} x a^\dagger a \right) \nonumber \\
            &\quad + \frac{1}{\me^\beta - 1} \left( -\frac{1}{2} a a^\dagger x + a x a^\dagger - \frac{1}{2} x a a^\dagger \right). \label{eq-generator-nicest-model}
        \end{align}
        The Gaussian QMS generated by \eqref{eq-generator-nicest-model} has the unique faithful normal invariant state $ \rho = (1-\me^{-\beta}) \me^{-\beta N} $, and its GNS-induced semigroup is self-adjoint. Moreover, $ \mathbf{Z} = \mathbf{Z}^T = - \frac{1}{2} \1 $, and
        \begin{align*}
            \chi_{T_t (Y) } (z) &= \frac{\coth(\beta/4)}{\pi ( 1 - \me^{-t}) } \exp\!\left\{ -\frac{1}{2} \mathbf{z}^T ( \coth (\beta/4) \tanh (t/2) \1 ) \mathbf{z}\right\} \\
            &\quad \cdot \int \chi_Y \left( \sech(t/2) \mathbf{z} - \boldsymbol{\xi} \right) \exp\!\left\{ -\frac{1}{2} \boldsymbol{\xi}^T ( \coth(\beta/4) \coth(t/2) \1) \boldsymbol{\xi}\right\} \, \md \boldsymbol{\xi}.
        \end{align*}
    \end{example}

    \subsection{Sobolev Regularization and Immediate Compactness}

    Corollary~\ref{corollary-characteristic-function-induced-semigroup} already implies immediate compactness because every Hilbert--Schmidt operator between Hilbert spaces is compact. We now prove the stronger regularization property
    \begin{equation*}
        T_t^{(s)}\big(\cB_2(\mathsf h)\big)\subset\cH,
        \qquad t>0,
    \end{equation*}
    which provides a second proof through a compact Sobolev embedding. This regularization parallels the classical Ornstein--Uhlenbeck phenomenon in which the semigroup maps $L^2$ immediately into a compactly embedded Sobolev space.

    For $Y\in\cB_2(\mathsf h)$, write $Y_{mn}:=\langle e_m,Ye_n\rangle$. Define

    \begin{equation} \label{eq-definition-first-order-derivative-sobolev-space-simplified}
        \cH
        :=\left\{Y\in\cB_2(\mathsf h):
        \sum_{m,n\geq0}(m+n+1)\abs{Y_{mn}}^2<\infty\right\}.
    \end{equation}
    Equivalently,
    \begin{equation*}
        \cH=\left\{Y\in\cB_2(\mathsf h):
        aY,\ Ya^\dagger\in\cB_2(\mathsf h)\right\}.
    \end{equation*}
    Here $Ya^\dagger$ denotes the closure of the product initially defined on $\dom a^\dagger$, while $aY\in\cB_2(\mathsf h)$ includes the condition $\ran Y\subset\dom a$. A direct calculation in the canonical basis gives
    \begin{equation*}
        \norm{aY}_2^2=\sum_{m,n\geq0}m\abs{Y_{mn}}^2,
        \qquad
        \norm{Ya^\dagger}_2^2
        =\sum_{m,n\geq0}n\abs{Y_{mn}}^2.
    \end{equation*}

    % \begin{remark} \label{remark-matrix-coefficient-Y-and-first-order-a-a-dagger}
    %     Let $ Y $ be a Hilbert-Schmidt operator. Since $ ( \ketbra{e_m}{e_n} )_{m, n \in \N} $ forms an orthonormal basis in $ \cB_2 (\mathsf{h}) $, we may write
    %     \begin{equation*}
    %         Y = \sum_{m, n \ge 0} Y_{m n} \ketbra{e_m}{e_n},
    %     \end{equation*}
    %     where $ Y_{m n} := \langle e_m, Y e_n \rangle $.
    %     A straightforward computation shows that
    %     \begin{align*}
    %         \norm{a Y}_2^2 &= \sum_{m, n \ge 0} m \abs{Y_{m n}}^2, \quad \norm{a^\dagger Y}_2^2 = \sum_{m, n \ge 0} (m + 1) \abs{Y_{m n}}^2, \\
    %         \norm{Y a}_2^2 &= \sum_{m, n \ge 0} (n+1) \abs{Y_{m, n}}^2, \quad \norm{Y a^\dagger}^2_2 = \sum_{m, n \ge 0} n \abs{Y_{m, n}}^2.
    %     \end{align*}
    % \end{remark}

    Equip $\cH$ with the inner product
    \begin{equation*}
        \langle X, Y \rangle_\cH := \langle X, Y \rangle_2 + \langle a X, a Y \rangle_2 + \langle X a^\dagger, Y a^\dagger \rangle_2.
    \end{equation*}
    The corresponding norm is
    \begin{equation} \label{eq-norm-first-order-derivative-sobolev-space}
        \norm{Y}_{\cH}^2 = \norm{Y}_2^2 + \norm{a Y}_2^2 + \norm{Y a^\dagger}_2^2 = \sum_{m, n \ge 0} (m+n+1) \abs{Y_{m n}}^2.
    \end{equation}
    Thus $\cH$ is isometrically isomorphic to the weighted sequence space $\ell^2(\N^2,m+n+1)$ and is therefore a Hilbert space.

    Let
    \begin{equation*}
        I : \cH \ni x \mapsto x \in \cB_2 (\mathsf{h}),
    \end{equation*}
    be the canonical embedding. We approximate $I$ in $\cB(\cH,\cB_2(\mathsf h))$ by finite-rank operators.

    \begin{proposition} \label{proposition-compactness-calligraphic-H}
        The embedding $ \cH \hookrightarrow \cB_2 (\mathsf{h}) $ is compact.
    \end{proposition}
    \begin{proof}
        For $Y\in\cH$, define
        \begin{equation*}
            P_M (Y) := \sum_{m + n < M} Y_{m n} \ketbra{e_m}{e_n}.
        \end{equation*}
        The operator $P_M$ has finite rank. Equation~\eqref{eq-norm-first-order-derivative-sobolev-space} gives
        \begin{align*}
            \norm{ Y - P_M Y }_2^2 &= \norm{\sum_{m+n \ge M} Y_{m n} \ketbra{e_m}{e_n} }^2_2 \\
            &= \sum_{m+n \ge M} \abs{Y_{m n}}^2 \le \sum_{m+n \ge M} \frac{m+n+1}{M} \abs{Y_{m n}}^2 \\
            &\le \frac{1}{M} \sum_{m, n \ge 0} (m+n+1) \abs{Y_{mn}}^2 = \frac{1}{M} \norm{Y}_{\cH}^2.
        \end{align*}
        Therefore,
        \begin{equation*}
            \norm{ I - P_M }_{\cB(\cH, \cB_2 (\mathsf{h}))} = \sup_{\norm{Y}_{\cH} = 1} \norm{ (I - P_M) Y }_2 \le \sup_{\norm{Y}_{\cH} = 1} \frac{1}{\sqrt{M}} \norm{Y}_\cH = \frac{1}{\sqrt{M}}.
        \end{equation*}
        Hence $P_M\to I$ in operator norm. Since every $P_M$ has finite rank, $I$ is compact.
    \end{proof}

    \begin{lemma}[Comparison with number-operator Sobolev spaces]
        \label{lemma-comparison-number-operator-sobolev-spaces}
        For $k>0$, following \cite{becker2021convergence}, define the bosonic Sobolev space of order $k$ associated with the number operator by
        \begin{equation*}
            \cW^k:=\left\{Y\in\cB_2(\mathsf h):
            (N+\1)^{k/4}Y(N+\1)^{k/4}\in\cB_2(\mathsf h)\right\}.
        \end{equation*}
        Equip it with the norm $\norm{Y}_{\cW^k}:=\norm{(N+\1)^{k/4}Y(N+\1)^{k/4}}_2$. Then every $\cW^k$ is compactly embedded in $\cB_2(\mathsf h)$, and
        \begin{equation*}
            \cW^2\hookrightarrow\cH\hookrightarrow\cW^1
        \end{equation*}
        continuously.
    \end{lemma}
    \begin{proof}
        The number-operator Sobolev norm has the coefficient representation
        \begin{equation*}
            \norm{Y}_{\cW^k}^2
            =\sum_{m,n\geq0}(m+1)^{k/2}(n+1)^{k/2}\abs{Y_{mn}}^2.
        \end{equation*}
        Since
        \begin{equation*}
            \sqrt{(m+1)(n+1)}\leq m+n+1\leq(m+1)(n+1),
        \end{equation*}
        the two continuous inclusions follow from \eqref{eq-norm-first-order-derivative-sobolev-space}. Moreover, if $m+n\geq M$, then $(m+1)(n+1)\geq M+1$. Consequently, truncation to the matrix coefficients with $m+n<M$ approximates the inclusion $\cW^k\hookrightarrow\cB_2(\mathsf h)$ in operator norm, with error at most $(M+1)^{-k/4}$. Thus this inclusion is compact.
    \end{proof}

    For the quantum characteristic function $\chi_Y$ of a Hilbert--Schmidt operator $Y$, the directional derivative $\nabla_w\chi_Y$ in a direction $w\in\C$ is understood distributionally.

    For $w\in\C$, let
    \begin{equation} \label{eq-definition-field-quadrature-p-w}
        p(w):=\frac{\mi}{\sqrt2}\left(w a^\dagger-\overline w\,a\right),
    \end{equation}
    with the self-adjoint closure of the operator initially defined on $\cD$ understood. This is the field quadrature in direction $w$; in particular, $p(1)=p$ is the momentum operator introduced in Section~\ref{section-preliminaries}.

    The relevant commutation relation is
    \begin{equation} \label{eq-commutation-relation-pw-Wz}
        [p(w), W(z)] = - \sqrt{2} \Im \langle z, w \rangle W(z).
    \end{equation}

    \begin{lemma}
        For $z,w\in\C$,
        \begin{equation} \label{eq-directional-derivative-W-z}
            \nabla_w W(z) := \left. \frac{\md}{\md t} W(z + t w) \right\vert_{t = 0} = \mi \Im \langle z, w \rangle W(z) - \mi \sqrt{2} W(z) p(w).
        \end{equation}
    \end{lemma}
    \begin{proof}
        For $u\in\dom p(w)$,
        \begin{align*}
            \frac{\md}{\md t} W(z + t w) u &= \frac{\md}{\md t} \me^{t \mi \Im \langle z, w \rangle} W(z) W(t w) u \\
            &= \mi \Im \langle z, w \rangle \me^{t \mi \Im \langle z, w \rangle} W(z) W(t w) u \\
            &\quad +  \me^{t \mi \Im \langle z, w \rangle} W(z) (-\mi \sqrt{2} p(w)) W(t w) u,
        \end{align*}
        which proves \eqref{eq-directional-derivative-W-z} on $\dom p(w)$.
    \end{proof}

    The following result is a Hilbert--Schmidt variant of \cite[Lemma~5.4.3]{holevo2011probabilistic}.

    \begin{lemma} \label{lemma-Y-pw-and-characteristic-function}
        Let $Y$ be a Hilbert--Schmidt operator with quantum characteristic function $\chi_Y$, and let $w\in\C$. If
        \begin{equation} \label{eq-characteristic-function-Y-pw}
            \frac{1}{\sqrt{2}} \left( \mi \nabla_w - \Im \langle z, w \rangle \right) \chi_Y (z)
        \end{equation}
        belongs to $L^2(\C,\md\mathbf z)$, then $Yp(w)$ has a Hilbert--Schmidt extension whose quantum characteristic function is \eqref{eq-characteristic-function-Y-pw}. If
        \begin{equation} \label{eq-characteristic-function-pw-Y}
            \frac{1}{\sqrt{2}} \left( \mi \nabla_w + \Im \langle z, w \rangle \right) \chi_Y (z)
        \end{equation}
        belongs to $L^2(\C,\md\mathbf z)$, then $\ran Y\subset\dom p(w)$ and $p(w)Y$ is Hilbert--Schmidt, with quantum characteristic function \eqref{eq-characteristic-function-pw-Y}.
    \end{lemma}
    \begin{proof}
        Suppose that \eqref{eq-characteristic-function-Y-pw} belongs to $L^2(\C,\md\mathbf z)$. By the quantum Plancherel theorem, it is the quantum characteristic function of a Hilbert--Schmidt operator $A$. For $u,v\in\cE$, consider the rank-one operator $\ketbra{u}{v}$. Then
        \begin{align}
            \chi_{\ketbra{u}{v} p(w)} (z) &= \Tr ( \ketbra{u}{v} p(w) W(z) ) = \langle v, p(w) W(z) u \rangle \nonumber \\
            &= \left\langle v, \left( -\frac{1}{\sqrt{2}} \Im \langle z, w \rangle W(z) + \frac{\mi}{\sqrt{2}} \nabla_w W(z) \right) u \right\rangle \nonumber  \\
            &= -\frac{1}{\sqrt{2}} \Im \langle z, w \rangle \chi_{ \ketbra{u}{v} } (z) + \frac{\mi}{\sqrt{2}} \nabla_w \chi_{\ketbra{u}{v} } (z). \label{eq-characteristic-function-rank-one-times-p-w}
        \end{align}
        Since $\chi_A$ is square-integrable and $\chi_{\ketbra{u}{v}}$ is a Schwartz function, integration by parts and \eqref{eq-characteristic-function-rank-one-times-p-w} give
        \begin{align*}
            \langle u, A v \rangle &= \Tr( ( \ketbra{u}{v} )^* A ) = \frac{1}{\pi} \int \overline{ \chi_{\ketbra{u}{v}} (z) } \cdot \chi_A (z) \, \md z \\
            &= \frac{1}{\sqrt{2} \pi} \int \overline{ \left( \mi \nabla_w - \Im \langle z, w \rangle \right) \chi_{\ketbra{u}{v}} (z) } \cdot \chi_Y (z) \, \md z \\
            &= \frac{1}{\pi} \int \overline{ \chi_{\ketbra{u}{v} p(w)} (z) } \cdot \chi_Y (z) \, \md z \\
            &= \Tr ( (\ketbra{u}{v} p(w) )^* Y ) = \langle u, Y p(w) v \rangle.
        \end{align*}
        Since $u,v$ are arbitrary and $\cE$ is a core for $p(w)$, $A$ is the Hilbert--Schmidt extension of $Yp(w)$. Its quantum characteristic function is \eqref{eq-characteristic-function-Y-pw}.

        The proof for $p(w)Y$ is analogous.
    \end{proof}

    \begin{remark} \label{remark-characteristic-function-Y-a-Y-a-dagger}
        We have
        \begin{equation*}
            p = p(1), \quad q = - p(\mi), \quad a = \frac{q + \mi p}{\sqrt{2}}, \quad a^\dagger = \frac{q - \mi p}{\sqrt{2}}.
        \end{equation*}
        Write $z=x+\mi y$ and define $\partial_z:=\partial_x+\mi\partial_y$ and $\partial_{\overline z}:=\partial_x-\mi\partial_y$. For a Hilbert--Schmidt operator $Y$,
        \begin{align*}
            \chi_{Y a} (z) &= \frac{1}{2} ( z - \partial_z ) \chi_Y (z), \quad \chi_{Y a^\dagger} (z) = \frac{1}{2} ( \overline{z} + \partial_{\overline{z}} ) \chi_Y (z), \\
            \chi_{a Y} (z) &= - \frac{1}{2} ( z + \partial_z ) \chi_Y (z), \quad \chi_{a^\dagger Y} (z) = \frac{1}{2} ( \partial_{\overline{z}} - \overline{z} ) \chi_Y (z),
        \end{align*}
        whenever the four products are well-defined Hilbert--Schmidt operators.
    \end{remark}

    For convenience, set
    \begin{equation}
        G_t (\mathbf{z}, \boldsymbol{\xi}) := \exp\!\left\{ -\frac{1}{2} \begin{pmatrix}
        \mathbf{z}^T & \boldsymbol{\xi}^T
        \end{pmatrix} \mathbf{G}_t \begin{pmatrix}
        \mathbf{z} \\ \boldsymbol{\xi}
        \end{pmatrix}\right\}.
    \end{equation}

    \begin{theorem} \label{theorem-sobolev-regularization}
        For every $t>0$ and $s\in[0,1]$, the map $T_t^{(s)}:\cB_2(\mathsf h)\to\cH$ is bounded. Consequently, $(T_t^{(s)})_{t\geq0}$ is immediately compact on $\cB_2(\mathsf h)$.
    \end{theorem}
    \begin{proof}
        Let $Y\in\cB_2(\mathsf h)$ and $t>0$. Corollary~\ref{corollary-characteristic-function-induced-semigroup} gives
        \begin{equation*}
            \chi_{T_t^{(s)} (Y)} (z) = \frac{c_t}{\pi} \int \chi_Y (- \boldsymbol{\xi}) G_t (\mathbf{z}, \boldsymbol{\xi}) \md \boldsymbol{\xi},
        \end{equation*}
        Because $\Re\mathbf G_t>0$, $G_t$ is a Schwartz function of the joint variable $(\mathbf z,\boldsymbol\xi)$.

        For fixed $\mathbf z$, the Cauchy--Schwarz inequality in $\boldsymbol\xi$ gives
        \begin{align*}
            \abs{ (z + \partial_z) \chi_{T_t^{(s)}(Y)} (z) }^2 &= \frac{c_t^2}{\pi^2} \abs{ \int \chi_Y (-\boldsymbol{\xi}) \cdot (z + \partial_z) G_t (\mathbf{z}, \boldsymbol{\xi}) \, \md \boldsymbol{\xi} }^2 \\
            &\le \frac{c_t^2}{\pi^2} \int \abs{ \chi_Y (-\boldsymbol{\xi}) } ^2  \, \md \boldsymbol{\xi} \cdot \int \abs{ (z + \partial_z) G_t (\mathbf{z}, \boldsymbol{\xi}) }^2 \, \md \boldsymbol{\xi} \\
            &= \frac{c_t^2}{\pi^2} \norm{ \chi_Y }^2_{L^2(\md \boldsymbol{\xi})} \cdot \int \abs{ (z + \partial_z) G_t (\mathbf{z}, \boldsymbol{\xi}) }^2 \, \md \boldsymbol{\xi},
        \end{align*}
        and hence
        \begin{align}
            \norm{ (z + \partial_z) \chi_{T_t^{(s)} (Y)} (z) }_{L^2 (\md \mathbf{z})}^2 &= \int \abs{ (z + \partial_z) \chi_{T_t^{(s)}(Y)} (z) }^2 \, \md \mathbf{z} \nonumber \\
            &\le \frac{c_t^2}{\pi^2} \norm{ \chi_Y }^2_{L^2(\md \boldsymbol{\xi})} \norm{ (z + \partial_z) G_t(\mathbf{z}, \boldsymbol{\xi}) }^2_{L^2 (\md \mathbf{z} \md \boldsymbol{\xi})}. \label{eq-estimate-partial-z-chi}
        \end{align}
        Every derivative of $G_t$ is a polynomial times $G_t$ and therefore remains a Schwartz function. The same argument gives
        \begin{equation} \label{eq-estimate-partial-overlinez-chi}
            \norm{ (\overline{z} + \partial_{\overline{z}}) \chi_{T_t^{(s)} (Y)} (z) }_{L^2 (\md \mathbf{z})}^2 \le \frac{c_t^2}{\pi^2} \norm{ \chi_Y }^2_{L^2(\md \boldsymbol{\xi})} \norm{ (\overline{z} + \partial_{\overline{z}}) G_t(\mathbf{z}, \boldsymbol{\xi}) }^2_{L^2 (\md \mathbf{z} \md \boldsymbol{\xi})}.
        \end{equation}
        Lemma~\ref{lemma-Y-pw-and-characteristic-function} and Remark~\ref{remark-characteristic-function-Y-a-Y-a-dagger} imply that $aT_t^{(s)}(Y)$ and $T_t^{(s)}(Y)a^\dagger$ are Hilbert--Schmidt. Equations~\eqref{eq-quantum-plancherel}, \eqref{eq-estimate-partial-z-chi}, and \eqref{eq-estimate-partial-overlinez-chi} yield
        \begin{align*}
            &\quad \norm{T_t^{(s)} (Y)}_{\cH}^2 \\
            &= \norm{T_t^{(s)} (Y)}_2^2 + \norm{ a \, T_t^{(s)} (Y) }_2^2 + \norm{ T_t^{(s)} (Y) \, a^\dagger }_2^2 \\
            &\le \left( \norm{T_t^{(s)}}^2 + \frac{c_t^2}{4\pi^2} \left( \norm{(z + \partial_z) G_t}^2_{L^2 (\md \mathbf{z} \md \boldsymbol{\xi})} + \norm{( \overline{z} + \partial_{\overline{z}} ) G_t}^2_{L^2 (\md \mathbf{z} \md \boldsymbol{\xi})} \right) \right) \norm{Y}_2^2.
        \end{align*}
        Hence $T_t^{(s)}:\cB_2(\mathsf h)\to\cH$ is bounded. Proposition~\ref{proposition-compactness-calligraphic-H} supplies the compact embedding $\cH\hookrightarrow\cB_2(\mathsf h)$, so $T_t^{(s)}$ is compact on $\cB_2(\mathsf h)$ for every $t>0$.
    \end{proof}

    \begin{corollary} \label{corollary-compact-resolvent-full-spectrum}
        For every $s\in[0,1]$, the induced generator $L^{(s)}$ and its adjoint $L^{(s)*}$ have compact resolvent. Consequently, if $\lambda$ and $\mu$ are the eigenvalues of the drift matrix $\mathbf Z$, then Theorem~\ref{theorem-other-eigenvalues-L} and Proposition~\ref{proposition-exhaustion-point-spectrum} yield
        \begin{equation*}
            \sigma(L^{(s)})
            =\{n\lambda+m\mu:n,m\in\N\},
            \quad
            \sigma(L^{(s)*})
            =\{n\overline\lambda+m\overline\mu:n,m\in\N\}.
        \end{equation*}
        These conclusions include the GNS and anti-GNS embeddings when their spectral gaps vanish.
    \end{corollary}
    \begin{proof}
        By \cite[Theorem II.4.29]{engel2000one}, immediate compactness of $(T_t^{(s)})_{t\geq0}$ implies that $L^{(s)}$ has compact resolvent. The same conclusion holds for $L^{(s)*}$ by taking adjoints. The displayed identities follow from the point-spectrum computation in Section~\ref{section-eigenvalues-of-L} and the standard spectral theory of operators with compact resolvent.
    \end{proof}

    \section*{Acknowledgements}

    The authors thank an anonymous referee for valuable suggestions concerning the organization and exposition of the paper.

    \begin{appendices}

    \section{Embedded Polynomials and the Domain of the Induced Generator} \label{section-proof-theorem-L-calL}

    This appendix proves the transfer identity in Theorem~\ref{theorem-L-and-calL-relation}; the proof of Theorem~\ref{theorem-base-eigenvalue-L} follows along the way. We begin with bounded operators in the domain of $\cL$, for which the identity is immediate.

    \begin{lemma} \label{lemma-L-calL-when-X-in-domain-of-calL}
        If $X\in\dom\cL$, then $\rho^{s/2}X\rho^{(1-s)/2}\in\dom L^{(s)}$ and
        \begin{equation} \label{eq-proof-section-important-expression}
            L^{(s)} (\rho^{s/2} X \rho^{(1-s)/2}) = \rho^{s/2} \cL (X) \rho^{(1-s)/2}.
        \end{equation}
    \end{lemma}
    \begin{proof}
        Let $X\in\dom\cL$ and $Y\in\cB_2(\mathsf h)$. Then
        \begin{align*}
            & \quad \lim_{t \rightarrow 0^+} \left\langle Y, \frac{T_t (\rho^{s/2} X \rho^{(1-s)/2}) - \rho^{s/2} X \rho^{(1-s)/2}}{t} \right\rangle_2 \\
            &= \lim_{t \rightarrow 0^+} \Tr \left( \rho^{(1-s)/2} Y^* \rho^{s/2} \frac{\cT_t (X) - X}{t} \right) \\
            &= \Tr \left( \rho^{(1-s)/2} Y^* \rho^{s/2} \cL (X) \right) = \left\langle Y, \rho^{s/2} \cL (X) \rho^{(1-s)/2} \right\rangle_2,
        \end{align*}
        The weak and strong generators of a strongly continuous semigroup coincide \cite[Chapter~2, Theorem~1.3]{pazy2012semigroups}, which completes the proof.
    \end{proof}

    We will extend \eqref{eq-proof-section-important-expression} to every polynomial $X$ in $a$ and $a^\dagger$, with $\cL(X)$ interpreted algebraically through \eqref{eq-gaussian-gksl-generator}.

    \subsection{Weyl Orbits and Generator-Domain Arguments}

    Sandwiching a polynomial in $a$ and $a^\dagger$ between $\rho^{s/2}$ and $\rho^{(1-s)/2}$ produces a Hilbert--Schmidt operator. We use this fact to differentiate sandwiched Weyl orbits in the Hilbert--Schmidt norm and to place their derivatives in the domain of $L^{(s)}$.

    The operator $za^\dagger-\overline z a$ generates the Weyl group $(W(rz))_{r\in\R}$, and Lemma~\ref{lemma-polynomial-rho-hilbert-schmidt} makes each of its sandwiched powers Hilbert--Schmidt. The sandwiched Weyl orbit therefore has the following Hilbert--Schmidt power-series expansion.

    \begin{lemma} \label{lemma-Wrz-rho-1/2-power-series}
        For every $r\in\R$ and $z\in\C$,
        \begin{equation} \label{eq-weyl-group-expansion-polynomial-a-a-dagger}
            \rho^{s/2} W(rz) \rho^{(1-s)/2} = \sum_{k=0}^{\infty} \frac{\rho^{s/2} r^k (z a^\dagger - \overline{z} a)^k \rho^{(1-s)/2}}{k!}
        \end{equation}
        with convergence in $\cB_2(\mathsf h)$.
    \end{lemma}

    \begin{proof}
        We first prove absolute convergence in the Hilbert--Schmidt norm.
        \begin{align*}
            &\quad \norm{\sum_{k=0}^\infty \frac{\rho^{s/2} r^k (z a^\dagger - \overline{z} a)^k \rho^{(1-s)/2}}{k!} }_2 \\
            &\le \sum_{k=0}^\infty \frac{\abs{r}^k}{k!} \left( \sum_{n=0}^\infty \norm{ \rho^{s/2} (z a^\dagger - \overline{z} a)^k \rho^{(1-s)/2} e_n }^2 \right)^{1/2} \\
            &= (1 - \me^{-\beta})^{1/2} \sum_{k=0}^\infty \frac{\abs{r}^k}{k!} \left( \sum_{n=0}^\infty \me^{-(1-s) \beta n} \norm{ \me^{-s \beta N / 2} (z a^\dagger - \overline{z} a)^k e_n }^2 \right)^{1/2} \\
            &\le (1 - \me^{-\beta})^{1/2} \sum_{k=0}^\infty \frac{\abs{r}^k}{k!} \left( \sum_{n=0}^\infty \me^{-(1-s) \beta n} \me^{-s \beta (n-k)} ( 2 \abs{z} )^{2k} \frac{(n+k)!}{n!} \right)^{1/2}  \\
            &= (1 - \me^{-\beta})^{1/2} \sum_{k=0}^\infty \frac{(2 \abs{r z})^{k}}{k!} \me^{s \beta k /2} \left( \sum_{n=0}^\infty \me^{-\beta n} \frac{(n+k)!}{n!} \right)^{1/2} \\
            &= \sum_{k=0}^\infty \frac{ \left( 2 \abs{rz} \me^{s \beta / 2} (1 - \me^{-\beta})^{-1/2} \right)^{k} }{(k!)^{1/2}} < + \infty.
        \end{align*}
        For every $n\in\N$, $e_n$ is an analytic vector for the generator of $(W(rz))_{r\in\R}$ because
        \begin{equation*}
            \sum_{k=0}^{\infty} \frac{\norm{r^k (z a^\dagger - \overline{z} a)^k e_n}}{k!} \le \sum_{k=0}^{\infty} \frac{  \abs{2 r z}^k \sqrt{\frac{(n+k)!}{n!}}}{k!} < +\infty.
        \end{equation*}
        Hence $\rho^{(1-s)/2}e_n=(1-\me^{-\beta})^{(1-s)/2}\me^{-(1-s)\beta n/2}e_n$ is also analytic, and
        \begin{equation*}
            \rho^{s/2} W(rz) \rho^{(1-s)/2} e_n = \sum_{k=0}^{\infty} \frac{ \rho^{s/2} r^k (z a^\dagger - \overline{z} a)^k \rho^{(1-s)/2}}{k!} e_n, \quad \forall n \in \N.
        \end{equation*}
        Thus the two sides of \eqref{eq-weyl-group-expansion-polynomial-a-a-dagger} agree on the dense finite-particle space. Both are bounded operators, so they agree on $\mathsf h$; the Hilbert--Schmidt convergence proved above completes the argument.
    \end{proof}

    For each $z\in\C$, the map $r\mapsto W(rz)$ is strongly continuous. After sandwiching by powers of $\rho$, it is differentiable in the Hilbert--Schmidt norm; we prove the more general form needed later.

    \begin{proposition} \label{proposition-derivative-of-sandwiched-Wrz}
        Let $ X, Y $ be polynomials in $ a $ and $ a^\dagger $. Then, the map
        \begin{equation*}
            r \mapsto \rho^{s/2} X W(rz) Y \rho^{(1-s)/2}
        \end{equation*}
        is differentiable in the Hilbert--Schmidt norm, with derivative
        \begin{equation} \label{eq-derivative-of-sandwiched-Wrz}
            \frac{\md}{\md r} \rho^{s/2} X W(rz) Y \rho^{(1-s)/2} = \rho^{s/2} X W(rz) (z a^\dagger - \overline{z} a) Y \rho^{(1-s)/2}.
        \end{equation}
    \end{proposition}
    \begin{proof}
        We first record a continuity observation. If $(K_r)_r\subset\cB_2(\mathsf h)$ is continuous in the Hilbert--Schmidt norm, then so are
        \begin{equation*}
            r \mapsto W(rz) K_r, \quad r \mapsto K_r W(rz)
        \end{equation*}
        Indeed, for any orthonormal basis $(e_n)_{n\in\N}$ of $\mathsf h$ and any sequence $r_k\to r$,
        \begin{align}
            & \quad \norm{ W(r_k z) K_{r_k} - W(r z) K_r }_2 \\
            &= \norm{ W(r_k z) K_{r_k} - W(r_k z) K_r + W(r_k z) K_r - W(r z) K_r }_2 \nonumber \\
            &\le \norm{ W(r_k z) } \norm{ K_{r_k} - K_r }_2 + \left( \sum_{n \ge 0} \norm{ ( W(r_k z) - W(r z) ) K_r e_n }^2  \right)^{1/2} \rightarrow 0,
        \end{align}
        by strong continuity of the Weyl group and dominated convergence. Continuity under right multiplication follows by taking adjoints.

        We next prove that $r\mapsto\rho^{s/2}XW(rz)Y\rho^{(1-s)/2}$ is Hilbert--Schmidt continuous. If $0\leq s<1$, then
        \begin{equation*}
            \rho^{s/2} X W(rz) Y \rho^{(1-s)/2} = \rho^{s/2} W(rz) \cdot W(-rz) X W(rz) Y \rho^{(1-s)/2}.
        \end{equation*}
        By Lemma~\ref{lemma-exchangeability-X-Wrz}, $W(-rz)XW(rz)$ is a polynomial in $a$ and $a^\dagger$ whose coefficients depend polynomially on $r$. Hence $W(-rz)XW(rz)Y\rho^{(1-s)/2}$ is Hilbert--Schmidt and varies continuously with $r$. The preceding observation about left multiplication gives the required continuity. When $s=1$, write instead
        \begin{equation*}
            \rho^{1/2} X W(rz) Y = \rho^{1/2} X W(rz) Y W(-rz) \cdot W(rz),
        \end{equation*}
        and use continuity under right multiplication.

        Replacing $Y$ by $(za^\dagger-\overline z a)Y$ in the same argument shows that
        \begin{equation} \label{eq-derivative-sandwiched-Wrz-0}
            r \mapsto \rho^{s/2} X W(r z) (z a^\dagger - \overline{z} a) Y \rho^{(1-s)/2}
        \end{equation}
        is Hilbert--Schmidt continuous.

        Let $u,v\in\cD$. Since the finite-particle space $\cD$ is invariant under $\rho$ and under polynomials in $a$ and $a^\dagger$,
        \begin{equation*}
            Y \rho^{(1-s)/2} v \in \cD \subset \dom (z a^\dagger - \overline{z} a).
        \end{equation*}
        The scalar fundamental theorem of calculus gives
        \begin{align}
            &\quad \left\langle u, \left( \rho^{s/2} X W( (r + t) z) Y \rho^{(1-s)/2} - \rho^{s/2} X W(rz) Y \rho^{(1-s)/2} \right)  v \right\rangle \nonumber \\
            &= \int_0^t \left\langle u, \rho^{s/2} X W((r + \tau) z) (z a^\dagger - \overline{z} a) Y \rho^{(1-s)/2} v \right\rangle \md \tau. \label{eq-derivative-sandwiched-Wrz-1}
        \end{align}
        The integrand in \eqref{eq-derivative-sandwiched-Wrz-1} is Hilbert--Schmidt continuous in $\tau$, so
        \begin{equation*}
            \int_0^t \rho^{s/2} X W((r + \tau) z) (z a^\dagger - \overline{z} a) Y \rho^{(1-s)/2} \md \tau
        \end{equation*}
        is a well-defined Bochner integral in $\cB_2(\mathsf h)$. Bounded linear functionals commute with Bochner integration, and \eqref{eq-derivative-sandwiched-Wrz-1} therefore gives
        \begin{align}
            &\quad \left\langle u, \left( \rho^{s/2} X W( (r + t) z) Y \rho^{(1-s)/2} - \rho^{s/2} X W(rz) Y \rho^{(1-s)/2} \right)  v \right\rangle \nonumber \\
            &= \left\langle u, \left( \int_0^t \rho^{s/2} X W((r + \tau) z) (z a^\dagger - \overline{z} a) Y \rho^{(1-s)/2} \md \tau \right) v \right\rangle \label{eq-derivative-sandwiched-Wrz-2}
        \end{align}
        Both operators inside the matrix elements in \eqref{eq-derivative-sandwiched-Wrz-2} are Hilbert--Schmidt and hence bounded on $\mathsf h$. Their difference has vanishing matrix elements on the dense set $\cD\times\cD$, so it vanishes identically. Thus, in $\cB_2(\mathsf h)$,
        \begin{align}
            &\quad \rho^{s/2} X W((r + t) z) Y \rho^{(1-s)/2} - \rho^{s/2} X W(rz) Y \rho^{(1-s)/2} \nonumber \\
            &= \int_0^t \rho^{s/2} X W( (r+\tau) z ) (z a^\dagger - \overline{z} a) Y \rho^{(1-s)/2} \, \md \tau. \label{eq-difference-by-integrating-derivative-Wz}
        \end{align}
        The integral is a Bochner integral in $\cB_2(\mathsf h)$.

        Dividing the left-hand side of \eqref{eq-difference-by-integrating-derivative-Wz} by $ t $ and subtracting the proposed derivative gives
        \begin{align*}
            &\quad \norm{ \frac{1}{t} \left( \rho^{\frac{s}{2}} X W((r + t) z) Y \rho^{\frac{1-s}{2}} - \rho^{\frac{s}{2}} X W(rz) Y \rho^{\frac{1-s}{2}} \right) - \rho^{\frac{s}{2}} X W(rz) (z a^\dagger - \overline{z} a) Y \rho^{\frac{1-s}{2}} }_2 \\
            &= \norm{ \frac{1}{t} \int_0^t \rho^{\frac{s}{2}} X W((r+\tau) z) (z a^\dagger - \overline{z} a) Y \rho^{\frac{1-s}{2}} \md \tau - \rho^{\frac{s}{2}} X W(rz) ( z a^\dagger - \overline{z} a ) Y \rho^{\frac{1-s}{2}} }_2 \\
            &= \norm{ \frac{1}{t} \int_0^t \rho^{\frac{s}{2}}  X \left( W((r + \tau) z ) - W (r z) \right) (z a^\dagger - \overline{z} a) Y \rho^{\frac{1-s}{2}} \md \tau }_2 \\
            &\le \frac{1}{\abs{t}} \int_{I_t} \norm{ \rho^{\frac{s}{2}} X ( W((r + \tau) z) - W(r z) ) (z a^\dagger - \overline{z} a) Y  \rho^{\frac{1-s}{2}} }_2 \md \tau \\
            &\le \sup_{\tau \in I_t} \norm{ \rho^{\frac{s}{2}} X ( W((r + \tau) z) - W(r z) ) (z a^\dagger - \overline{z} a) Y  \rho^{\frac{1-s}{2}} }_2 \rightarrow 0,
        \end{align*}
        as $t\to0$, by the Hilbert--Schmidt continuity in \eqref{eq-derivative-sandwiched-Wrz-0}. Here
        \begin{equation*}
            I_t := [ \min \{0, t\}, \max\{0, t\} ]
        \end{equation*}
        because $t$ may be negative. This proves \eqref{eq-derivative-of-sandwiched-Wrz}.
    \end{proof}

    We next use the Weyl action in Theorem~\ref{theo-gaussian-qms-explicit-weyl-operators} to compute $L^{(s)}$ on the embedded first-order polynomials. The following weak generator formula is adapted from the proof of \cite[Theorem~2]{agredo2021gaussian}. The unbounded operators $p$ and $q$ in \eqref{eq-Yz} explain why the Weyl operator itself need not belong to $\dom\cL$.

    \begin{proposition} \label{prop-generator-wz-very-weak-form}
        For $\xi,\eta\in\cE$,
        \begin{equation*}
            \lim_{t \rightarrow 0^+} \left\langle \xi, \frac{\cT_t (W(z)) - W(z)}{t} \eta \right\rangle = \langle \xi, W(z) Y(z) \eta \rangle,
        \end{equation*}
        where
        \begin{align}
            Y(z) &= (Z z a^\dagger - \overline{Z z} a) + \frac{1}{2} ( \overline{z} Z z - \overline{Z z} z ) - \frac{1}{2} \Re{\overline{z} C z} + \mi \Re{\overline{\zeta} z} \nonumber \\
            &= \mi \sqrt{2} (\Im{Z z} q - \Re{Z z} p) + \frac{1}{2} ( \overline{z} Z z - \overline{Z z} z ) - \frac{1}{2} \Re{\overline{z} C z} + \mi \Re{\overline{\zeta} z}. \label{eq-Yz}
        \end{align}
    \end{proposition}

    Although $W(z)$ need not lie in $\dom\cL$, Theorem~\ref{theo-Wz-rho-one-half-domain-of-L} shows that the sandwiched Weyl operator $\rho^{s/2}W(z)\rho^{(1-s)/2}$ belongs to $\dom L^{(s)}$. We first need two auxiliary lemmas.

    The first lemma is the standard density of rank-one operators generated by a dense subspace.
    \begin{lemma} \label{lemma-density-rank-one-operators}
        The space $\myspan\{\ketbra{\xi}{\eta}:\xi,\eta\in\cE\}$ is dense in $\cB_2(\mathsf h)$.
    \end{lemma}

    \begin{lemma} \label{lemma-uniform-bddness-in-t}
        The map
        \begin{equation} \label{eq-map-uniformly-bounded-in-t}
            t \mapsto t^{-1} \left(T_t^{(s)} (\rho^{s/2} W(z) \rho^{(1-s)/2}) - \rho^{s/2} W(z) \rho^{(1-s)/2} \right)
        \end{equation}
        is uniformly bounded in the Hilbert--Schmidt norm for $t>0$.
    \end{lemma}
    \begin{proof}
        Equation~\eqref{eq-introduction-gaussian-weyl-action} gives
        \begin{align}
            &\quad \norm{t^{-1} \left(T_t^{(s)} (\rho^{s/2} W(z) \rho^{(1-s)/2}) - \rho^{s/2} W(z) \rho^{(1-s)/2} \right)}_2 \nonumber \\
            &= \norm{t^{-1} \rho^{s/2} \left( \cT_t (W (z)) - W(z) \right) \rho^{(1-s)/2} }_2 \nonumber \\
            &= \norm{t^{-1} \rho^{s/2} \left( \phi_t (z) W(\me^{t Z} z) - \phi_t (z) W(z) + \phi_t (z) W(z)- W(z) \right) \rho^{(1-s)/2} }_2 \nonumber \\
            &\le \norm{t^{-1} \rho^{s/2} (W(\me^{t Z} z) - W(z)) \rho^{(1-s)/2}}_2 \nonumber \\
            &\quad + t^{-1} \abs{\phi_t (z) - 1} \norm{\rho^{s/2} W(z) \rho^{(1-s)/2}}_2. \label{ineq-Wetz-Wz}
        \end{align}
        In the last step we used $\abs{\phi_t(z)}\leq1$. The explicit formula for $\phi_t$ shows that $t^{-1}\abs{\phi_t(z)-1}$ is uniformly bounded for $t>0$. It remains to control the first term in \eqref{ineq-Wetz-Wz}.

        The proof of \cite[Proposition~20.14]{parthasarathy1992introduction} gives
        \begin{equation} \label{eq-differentiaion-exponential-vector}
            \frac{\md}{\md t} e (f + \me^{t Z} z) = ( Z \me^{t Z} z ) a^\dagger e(f + \me^{t Z} z),
        \end{equation}
        for $f\in\C$. Equation~\eqref{eq-differentiaion-exponential-vector} and Lemma~\ref{lemma-a-Wz-a-dagger-Wz-commutation-relations} then yield
        \begin{align}
            &\quad \frac{\md}{\md t} W(\me^{t Z} z) e(f) \nonumber \\
            &= \frac{\md}{\md t} \exp\!\left\{ - \frac{1}{2} \abs{\me^{t Z} z}^2 - \overline{\me^{t Z} z} f\right\} e(f + \me^{t Z} z) \nonumber \\
            &= \exp\!\left\{ - \frac{1}{2} \abs{\me^{t Z} z}^2 - \overline{\me^{t Z} z} f\right\} \left( - \Re{\overline{\me^{t Z} z} Z \me^{t Z} z} - \overline{Z \me^{t Z} z} f \right) e (f + \me^{t Z} z) \nonumber \\
            &\quad + \exp\!\left\{ - \frac{1}{2} \abs{\me^{t Z} z}^2 - \overline{\me^{t Z} z} f\right\} ( Z \me^{t Z} z ) a^\dagger e(f + \me^{t Z} z) \nonumber \\
            &= W(\me^{t Z} z) \left( - \Re{\overline{\me^{t Z} z} Z \me^{t Z} z} - \overline{Z \me^{t Z} z} a \right) e(f) +  \left( Z \me^{t Z} z \right) a^\dagger W(\me^{t Z} z) e(f) \nonumber \\
            &= W(\me^{t Z} z) \left( Z \me^{t Z} z a^\dagger - \overline{Z \me^{t Z} z} a + \mi \Im{ \overline{\me^{t Z} z} Z \me^{t Z} z} \right) e(f). \label{eq-derivative-WetZz-ef}
        \end{align}
        Integrating \eqref{eq-derivative-WetZz-ef} gives
        \begin{align}
            &\quad \rho^{s/2} ( W(\me^{t Z} z) - W(z) ) \rho^{(1-s)/2} \nonumber \\
            &= \int_0^t \rho^{s/2} W(\me^{\tau Z} z) \left( Z \me^{\tau Z} z a^\dagger - \overline{Z \me^{\tau Z} z} a + \mi \Im{ \overline{\me^{\tau Z} z} Z \me^{\tau Z} z} \right) \rho^{(1-s)/2} \md \tau. \label{eq-difference-of-WetZz-Wz-by-integral}
        \end{align}
        Set
        \begin{equation*}
            M_z := \max_{\tau \ge 0} \{ \abs{Z \me^{\tau Z} z}, \abs{ \overline{\me^{\tau Z} z } Z \me^{\tau Z} z }  \}.
        \end{equation*}
        For $0\leq s<1$, \eqref{eq-difference-of-WetZz-Wz-by-integral} implies
        \begin{align*}
            &\quad t^{-1} \norm{ \rho^{s/2} (W(\me^{t Z} z) - W(z)) \rho^{(1-s)/2} }_2 \\
            &\le M_z \norm{\rho^{s/2}} \left( \norm{a^\dagger \rho^{(1-s)/2} }_2 + \norm{ a \rho^{(1-s)/2} }_2 + \norm{ \rho^{(1-s)/2} }_2 \right)
        \end{align*}
        Remark~\ref{remark-X-rho-s-trace-class} shows that $a^\dagger\rho^{(1-s)/2}$ and $a\rho^{(1-s)/2}$ are Hilbert--Schmidt, proving the required uniform bound for $0\leq s<1$. When $s=1$, rewrite \eqref{eq-derivative-WetZz-ef} as
        \begin{equation*}
            \frac{\md}{\md t} W(\me^{t Z} z) e(f) = \left( Z \me^{t Z} z a^\dagger - \overline{Z \me^{t Z} z} a - \mi \Im{ \overline{\me^{t Z} z} Z \me^{t Z} z } \right) W(\me^{t Z} z) e(f),
        \end{equation*}
        and apply the same argument with the unbounded factors on the left.
    \end{proof}

    \begin{theorem} \label{theo-Wz-rho-one-half-domain-of-L}
        For every $z\in\C$, the operator $\rho^{s/2}W(z)\rho^{(1-s)/2}$ belongs to $\dom L^{(s)}$.
    \end{theorem}
    \begin{proof}
        We first compute the time derivative of $T_t^{(s)}(\rho^{s/2}W(z)\rho^{(1-s)/2})$ weakly. For an exponential vector $e(f)$,
        \begin{align*}
            \rho^{s/2} e(f) &= (1-\me^{-\beta})^{s/2} \sum_{\ell=0}^\infty \me^{- s \beta \ell /2} \ketbra{e_\ell}{e_\ell} \sum_{k=0}^\infty \frac{f^k}{\sqrt{k!}} e_k \\
            &= (1-\me^{-\beta})^{s/2} \sum_{k=0}^\infty \frac{f^k (\me^{-s \beta/2})^k}{\sqrt{k!}} e_k \\
            &= (1-\me^{-\beta})^{s/2} e(\me^{- s\beta/2} f).
        \end{align*}
        Hence $\rho^{s/2}\xi$ and $\rho^{(1-s)/2}\eta$ belong to $\cE$ whenever $\xi,\eta\in\cE$. Proposition~\ref{prop-generator-wz-very-weak-form} gives
        \begin{align}
            &\quad \lim_{t \rightarrow 0^+} \left\langle \xi, \frac{T_t^{(s)} ( \rho^{s/2} W(z) \rho^{(1-s)/2}) - \rho^{s/2} W(z) \rho^{(1-s)/2}}{t} \eta \right\rangle \nonumber \\
            &= \lim_{t \rightarrow 0^+} \left\langle \rho^{s/2} \xi, \frac{\cT_t(W(z))-W(z)}{t} \rho^{(1-s)/2} \eta   \right\rangle \nonumber \\
            &= \langle \xi, \rho^{s/2} W(z) Y(z) \rho^{(1-s)/2} \eta \rangle, \label{eq-existence-of-limit-in-xi-eta-form}
        \end{align}
        where $Y(z)$ is defined in \eqref{eq-Yz}. Lemma~\ref{lemma-polynomial-rho-hilbert-schmidt} also shows that $\rho^{s/2}W(z)Y(z)\rho^{(1-s)/2}$ is Hilbert--Schmidt.

        Equivalently,
        \begin{align}
            &\quad \lim_{t \rightarrow 0^+} \left\langle \ketbra{\xi}{\eta}, \frac{T_t^{(s)} ( \rho^{s/2} W(z) \rho^{(1-s)/2}) - \rho^{s/2} W(z) \rho^{(1-s)/2}}{t} \right\rangle_2  \\
            &= \lim_{t \rightarrow 0^+} \left\langle \ketbra{\xi}{\eta}, \rho^{s/2} W(z) Y(z) \rho^{(1-s)/2} \right\rangle.
        \end{align}
        By linearity, the same limit holds against every operator in $\myspan\{\ketbra{\xi}{\eta}:\xi,\eta\in\cE\}$. Lemmas~\ref{lemma-density-rank-one-operators} and~\ref{lemma-uniform-bddness-in-t} then imply
        \begin{align*}
            \frac{T_t^{(s)} ( \rho^{s/2} W(z) \rho^{(1-s)/2}) - \rho^{s/2} W(z) \rho^{(1-s)/2}}{t} \rightarrow \rho^{s/2} W(z) Y(z) \rho^{(1-s)/2}
        \end{align*}
        weakly in $\cB_2(\mathsf h)$. Since the weak and strong generators of a strongly continuous semigroup coincide, $\rho^{s/2}W(z)\rho^{(1-s)/2}\in\dom L^{(s)}$, and
        \begin{equation} \label{eq-action-L-s-on-wrapped-Wz}
            L^{(s)} \left( \rho^{s/2} W(z) \rho^{(1-s)/2} \right) = \rho^{s/2} W(z) Y(z) \rho^{(1-s)/2},
        \end{equation}
        where $Y(z)$ is given by \eqref{eq-Yz}.
    \end{proof}

    \begin{lemma} \label{lemma-differentiability-L-W-rz-rho}
        The map
        \begin{equation} \label{eq-map-r-to-L-wrapped-Wrz}
            r \mapsto L^{(s)} (\rho^{s/2} W(rz) \rho^{(1-s)/2})
        \end{equation}
        is differentiable in the Hilbert--Schmidt norm.
    \end{lemma}
    \begin{proof}
        Proposition~\ref{proposition-derivative-of-sandwiched-Wrz} gives
        \begin{equation*}
            \frac{\md}{\md r} \rho^{s/2} W(rz) (Zz a^\dagger - \overline{Zz} a) \rho^{(1-s)/2} = \rho^{s/2} W(rz) (z a^\dagger - \overline{z} a) (Zz a^\dagger - \overline{Zz} a) \rho^{(1-s)/2}
        \end{equation*}
        in the Hilbert--Schmidt norm. Theorem~\ref{theo-Wz-rho-one-half-domain-of-L} therefore yields
        \begin{align}
            &\quad \frac{\md}{\md r} L^{(s)} \left( \rho^{s/2} W(rz) \rho^{(1-s)/2} \right) = \frac{\md}{\md r} \left( \rho^{s/2} W(rz) Y(rz) \rho^{(1-s)/2} \right) \nonumber \\
            &= \rho^{s/2} W (rz) (Zz a^\dagger - \overline{Zz} a) \rho^{(1-s)/2} \nonumber \\
            &\quad + r \rho^{s/2} W(rz) (za^\dagger - \overline{z}a) (Zz a^\dagger - \overline{Zz} a)\rho^{(1-s)/2} \nonumber \\
            &\quad + r \rho^{s/2} \left( \overline{z} Z z - \overline{Zz} z - \Re{\overline{z} Cz} \right) W (rz) \rho^{(1-s)/2} \nonumber\\
            &\quad + \frac{1}{2} r^2 \rho^{s/2} \left( \overline{z} Z z - \overline{Zz} z - \Re{\overline{z} Cz} \right) W(rz) (z a^\dagger - \overline{z} a) \rho^{(1-s)/2}  \nonumber \\
            &\quad + \mi \Re{\overline{\zeta} z} \rho^{s/2}  W(rz) \rho^{(1-s)/2} \nonumber \\
            &\quad + \mi r \Re{\overline{\zeta} z} \rho^{s/2} W(rz) (z a^\dagger - \overline{z} a) \rho^{(1-s)/2}, \label{eq-explicit-derivative-L-wrapped-Wrz}
        \end{align}
        which explicitly evaluates the derivative in \eqref{eq-map-r-to-L-wrapped-Wrz}.
    \end{proof}

    Repeated application of Proposition~\ref{proposition-derivative-of-sandwiched-Wrz} gives the following corollary.

    \begin{corollary} \label{corollary-infinitely-differentiability-L-W-rz-rho}
        The map in \eqref{eq-map-r-to-L-wrapped-Wrz} is infinitely differentiable in $r$.
    \end{corollary}

    We can now interchange $L^{(s)}$ with differentiation in $r$.

    \begin{proposition} \label{theo-exchangeability-L-derivative}
        We have
        \begin{equation} \label{eq-exchange-differentiaion-generator}
            L^{(s)} \left( \frac{\md }{\md r} \left( \rho^{s/2} W (rz) \rho^{(1-s)/2} \right) \right) = \frac{\md}{\md r} L^{(s)} \left( \rho^{s/2} W(rz) \rho^{(1-s)/2} \right).
        \end{equation}
    \end{proposition}
    \begin{proof}
        Let $(s_n)\subset\R^+$ satisfy $s_n\to0$. Proposition~\ref{proposition-derivative-of-sandwiched-Wrz} gives
        \begin{align*}
            \dom L^{(s)} &\ni \frac{\rho^{s/2} W( (r+s_n) z) \rho^{(1-s)/2} - \rho^{s/2} W(rz) \rho^{(1-s)/2}}{s_n} \\
            &\rightarrow \frac{\md}{\md r} \rho^{s/2} W (rz) \rho^{(1-s)/2}
        \end{align*}
        in $\cB_2(\mathsf h)$. Lemma~\ref{lemma-differentiability-L-W-rz-rho} also gives
        \begin{align*}
            &\quad L^{(s)} \left( \frac{ \rho^{s/2} W((r+s_n)z) \rho^{(1-s)/2} - \rho^{s/2} W(rz) \rho^{(1-s)/2}}{s_n} \right) \\
            &\rightarrow \frac{\md}{\md r} L^{(s)} \left( \rho^{s/2} W(rz) \rho^{(1-s)/2} \right) \in \cB_2 (\mathsf{h}).
        \end{align*}
        Closedness of $L^{(s)}$ now implies
        \begin{equation*}
            \frac{\md}{\md r} \left(\rho^{s/2} W(rz) \rho^{(1-s)/2}\right) = \rho^{s/2} W(rz) (z a^\dagger - \overline{z} a) \rho^{(1-s)/2} \in \dom L^{(s)}
        \end{equation*}
        and \eqref{eq-exchange-differentiaion-generator} holds.
    \end{proof}

    \subsection{From Weyl Orbits to Embedded Polynomials}

    We now determine the action of $L^{(s)}$ on the embedded creation and annihilation operators.

    \begin{lemma} \label{lemma-L-a-rho-a-dagger-rho}
        We have
        \begin{align}
                L^{(s)} \left( \rho^{s/2} a^\dagger \rho^{(1-s)/2} \right) &= (- \gamma + \mi \Omega) \rho^{s/2} a^\dagger \rho^{(1-s)/2} + \mi \overline{ \kappa} \rho^{s/2} a \rho^{(1-s)/2}, \label{eq-L-a-dagger-rho}\\
                L^{(s)} \left( \rho^{s/2} a \rho^{(1-s)/2} \right) &= (- \gamma - \mi \Omega) \rho^{s/2} a \rho^{(1-s)/2} - \mi \kappa \rho^{s/2} a^\dagger \rho^{(1-s)/2}. \label{eq-L-a-rho}
        \end{align}
    \end{lemma}
    \begin{proof}
        The invariant state is diagonal, so \eqref{eq-diagonal-rho-parameters} gives $\zeta=0$. Setting $r=0$ in \eqref{eq-explicit-derivative-L-wrapped-Wrz} and using \eqref{eq-exchange-differentiaion-generator}, we obtain
        \begin{equation} \label{eq-base-eigenvalue-generator-L-pre}
            L^{(s)} \left( \rho^{s/2} (z a^\dagger - \overline{z} a) \rho^{(1-s)/2} \right) = \rho^{s/2} (Zz a^\dagger - \overline{Z z} a) \rho^{(1-s)/2}.
        \end{equation}
        Replacing $z$ by $\mi z$ in \eqref{eq-base-eigenvalue-generator-L-pre} gives
        \begin{equation} \label{eq-base-eigenvalue-generator-L-pre-2}
            L^{(s)} \left( \rho^{s/2} (z a^\dagger + \overline{z} a) \rho^{(1-s)/2} \right) = \rho^{s/2} \left( - \mi Z (\mi z) a^\dagger - \overline{\mi Z (\mi z)} a \right) \rho^{(1-s)/2}.
        \end{equation}
        Adding \eqref{eq-base-eigenvalue-generator-L-pre} and \eqref{eq-base-eigenvalue-generator-L-pre-2} yields
        \begin{align*}
            2 L^{(s)} \left( z \rho^{s/2} a^\dagger \rho^{(1-s)/2} \right) &= \left( Z z - \mi Z (\mi z) \right) \rho^{s/2} a^\dagger \rho^{(1-s)/2} \\
            &\quad - \left( \overline{Z z} + \overline{\mi Z (\mi z)} \right) \rho^{s/2} a \rho^{(1-s)/2} \\
            &= 2 (- \gamma + \mi \Omega) z \rho^{s/2} a^\dagger \rho^{(1-s)/2} + 2 \mi \overline{ \kappa} z \rho^{s/2} a \rho^{(1-s)/2},
        \end{align*}
        which proves \eqref{eq-L-a-dagger-rho}. Subtracting \eqref{eq-base-eigenvalue-generator-L-pre} from \eqref{eq-base-eigenvalue-generator-L-pre-2} gives
        \begin{align*}
            2 L^{(s)} \left( \overline{z} \rho^{s/2} a \rho^{(1-s)/2} \right) &= \left( - \mi Z(\mi z) - Z z \right) \rho^{s/2} a^\dagger \rho^{(1-s)/2} \\
            &\quad + \left(\overline{Z z} - \overline{\mi Z (\mi z)} \right) \rho^{s/2} a \rho^{(1-s)/2} \\
            &= - 2 \mi \kappa \overline{z} \rho^{s/2} a^\dagger \rho^{(1-s)/2} + 2 \left(- \gamma - \mi \Omega \right) \overline{z} \rho^{s/2} a \rho^{(1-s)/2},
        \end{align*}
        which proves \eqref{eq-L-a-rho}.
    \end{proof}

    We can now identify the fundamental eigenvalues of $L^{(s)}$.

    \begin{proof}[Proof of Theorem \ref{theorem-base-eigenvalue-L}]
        Suppose first that $\kappa=0$, so the real-linear operator $Z$ is complex linear. The drift matrix $\mathbf Z$ is diagonalizable, with eigenvalues $-\gamma+\mi\Omega$ and $-\gamma-\mi\Omega$. Assumption~\ref{hyp-stability-Z} excludes zero. Equations~\eqref{eq-L-a-dagger-rho} and \eqref{eq-L-a-rho} give
        \begin{align*}
            L^{(s)} \left( \rho^{s/2} a^\dagger \rho^{(1-s)/2} \right) &= ( - \gamma + \mi \Omega) \rho^{s/2} a^\dagger \rho^{(1-s)/2}, \\
            L^{(s)} \left( \rho^{s/2} a \rho^{(1-s)/2} \right) &= (- \gamma - \mi \Omega) \rho^{s/2} a \rho^{(1-s)/2},
        \end{align*}
        so both drift eigenvalues are eigenvalues of $L^{(s)}$.

        Now suppose that $\kappa\neq0$. For $c_1,c_2\in\C$,
        \begin{align*}
            L^{(s)} ( c_1 \rho^{s/2} a^\dagger \rho^{(1-s)/2} + c_2 \rho^{s/2} a \rho^{(1-s)/2} ) &= ( c_1 (- \gamma + \mi \Omega) - c_2 \mi \kappa ) \rho^{s/2} a^\dagger \rho^{(1-s)/2} \\
            &\quad + ( c_2 (- \gamma - \mi \Omega) + c_1 \mi \overline{ \kappa} ) \rho^{s/2} a \rho^{(1-s)/2}.
        \end{align*}
        A nonzero linear combination is an eigenvector with eigenvalue $\lambda\in\C$ precisely when
        \begin{equation} \label{eq-relation-c-c1-c2}
            c_1 (- \gamma + \mi \Omega) - c_2 \mi \kappa = \lambda c_1, \quad c_2 (- \gamma - \mi \Omega) + c_1 \mi \overline{ \kappa} = \lambda c_2.
        \end{equation}
        The stability assumption implies that $\lambda=0$ yields only $c_1=c_2=0$. For a nontrivial solution, the determinant of the system in \eqref{eq-relation-c-c1-c2} must vanish, so
        \begin{equation*}
            \lambda^2 + 2 \gamma \lambda + \left( \gamma^2 + \Omega^2 - \abs{\kappa}^2 \right) = 0,
        \end{equation*}
        whose roots are
        \begin{equation} \label{eq-roots-eigenvalues-L}
            \lambda = - \gamma \pm \sqrt{\abs{\kappa}^2 - \Omega^2}.
        \end{equation}
        These are precisely the eigenvalues of $\mathbf Z$. If the two roots are distinct, the corresponding linear combinations are eigenvectors; if they coincide, the first-order restriction has the Jordan structure stated in the theorem.
    \end{proof}

    Theorem~\ref{theorem-base-eigenvalue-L} identifies the fundamental eigenvectors among embedded first-order polynomials. We now complete the passage to arbitrary embedded polynomials by differentiating the sandwiched Weyl orbits and proving the transfer identity in Theorem~\ref{theorem-L-and-calL-relation}.

    \begin{proposition} \label{proposition-exchangeability-L-derivative-n-th-derivative}
        For every $n\geq1$,
        \begin{equation*}
            \frac{\md^n}{\md r^n} \rho^{s/2} W(rz) \rho^{(1-s)/2} \in \dom L^{(s)},
        \end{equation*}
        and
        \begin{equation} \label{eq-exchangeability-L-derivative-n-th-derivative}
            L^{(s)} \left( \frac{\md^n }{\md r^n} \rho^{s/2} W (rz) \rho^{(1-s)/2} \right) = \frac{\md^n}{\md r^n} L^{(s)} \left( \rho^{s/2} W(rz) \rho^{(1-s)/2} \right).
        \end{equation}
    \end{proposition}
    \begin{proof}
        Proposition~\ref{proposition-derivative-of-sandwiched-Wrz} gives, for every $n\geq1$,
        \begin{equation} \label{eq-weyl-group-rho-n-derivative}
            \frac{\md^n }{\md r^n} \left( \rho^{s/2} W (rz) \rho^{(1-s)/2} \right) = \rho^{s/2} W (rz) (z a^\dagger - \overline{z} a)^n \rho^{(1-s)/2}.
        \end{equation}
        Assume inductively that
        \begin{equation*}
            \frac{\md^n}{\md r^n} \rho^{s/2} W (rz) \rho^{(1-s)/2} = \rho^{s/2} W(rz) (z a^\dagger - \overline{z} a)^n \rho^{(1-s)/2} \in \dom L^{(s)}
        \end{equation*}
        and that \eqref{eq-exchangeability-L-derivative-n-th-derivative} holds at order $n$. Let $(s_k)_{k\in\N}\subset\R^+$ satisfy $s_k\to0$. Equation~\eqref{eq-weyl-group-rho-n-derivative} gives
        \begin{align*}
            &\quad  \frac{\rho^{s/2} W((r+s_k)z) (z a^\dagger - \overline{z} a)^n \rho^{(1-s)/2} - \rho^{s/2} W(rz) (z a^\dagger - \overline{z} a)^n \rho^{(1-s)/2}}{s_k} \\
            &\rightarrow \frac{\md^{n+1}}{\md r^{n+1}} \rho^{s/2} W (rz) \rho^{(1-s)/2},
        \end{align*}
        as $k\to\infty$ in $\cB_2(\mathsf h)$. The induction hypothesis yields
        \begin{align*}
            &\quad \lim_{k \rightarrow +\infty} L^{(s)} \left( \frac{ \rho^{s/2} (W((r+s_k)z) - W(rz)) (z a^\dagger - \overline{z} a)^n \rho^{(1-s)/2}}{s_k} \right) \\
            &= \lim_{k \rightarrow +\infty} \frac{1}{s_k} \left( L^{(s)} \left( \frac{\md^n}{\md r^n} \rho^{s/2} W( (r + s_k) z) \rho^{(1-s)/2} \right) \right. \\
            &\qquad \qquad \left. - L^{(s)} \left( \rho^{s/2} \frac{\md^n}{\md r^n} W(rz) \rho^{(1-s)/2} \right) \right) \\
            &= \lim_{k \rightarrow +\infty} \frac{1}{s_k} \cdot \frac{\md^n}{\md r^n} \left( L^{(s)} \left( \rho^{s/2} W( (r + s_k) z) \rho^{(1-s)/2} \right) \right. \\
            &\qquad \qquad \left. - \frac{\md^n}{\md r^n}  L^{(s)} \left( \rho^{s/2} W(rz) \rho^{(1-s)/2} \right) \right) \\
            &= \lim_{k \rightarrow +\infty} \frac{\md^n}{\md r^n} L^{(s)} \left( \frac{\rho^{s/2} W((r+s_k)z) \rho^{(1-s)/2} - \rho^{s/2} W(rz) \rho^{(1-s)/2}}{s_k} \right) \\
            &= \frac{\md}{\md r} \frac{\md^{n}}{\md r^{n}} L^{(s)} \left( \rho^{s/2} W(rz) \rho^{(1-s)/2} \right) = \frac{\md^{n+1}}{\md r^{n+1}} L^{(s)} \left( \rho^{s/2} W(rz) \rho^{(1-s)/2} \right).
        \end{align*}
        Corollary~\ref{corollary-infinitely-differentiability-L-W-rz-rho} guarantees the existence of the derivatives on the right. Closedness of $L^{(s)}$ now gives
        \begin{equation*}
            \frac{\md^{n+1}}{\md r^{n+1}} \rho^{s/2} W(rz) \rho^{(1-s)/2} \in \dom L^{(s)}
        \end{equation*}
        and \eqref{eq-exchangeability-L-derivative-n-th-derivative} at order $n+1$. The case $n=1$ is Proposition~\ref{theo-exchangeability-L-derivative}, so induction completes the proof.
    \end{proof}

    \begin{corollary} \label{corollary-X-rho-one-half-L-domain}
        If $X$ is a polynomial in $a$ and $a^\dagger$, then $\rho^{s/2}X\rho^{(1-s)/2}\in\dom L^{(s)}$.
    \end{corollary}
    \begin{proof}
        Proposition~\ref{proposition-exchangeability-L-derivative-n-th-derivative} shows that
        \begin{equation*}
            \rho^{s/2}W(rz)(za^\dagger-\overline z a)^n\rho^{(1-s)/2}
            \in\dom L^{(s)}
            \quad(r\in\R,\ z\in\C,\ n\in\N).
        \end{equation*}
        It suffices to obtain every normally ordered monomial from these operators.

        Set $r=0$ and fix $n\in\N$. Then
        \begin{equation*}
            \rho^{s/2} (z a^\dagger - \overline{z} a)^n \rho^{(1-s)/2} = \sum_{k=0}^n z^k (-\overline{z})^{n-k} \rho^{s/2} Y_{k, n-k} \rho^{(1-s)/2},
        \end{equation*}
        where $Y_{k,n-k}$ is the sum, with coefficient one, of all monomials containing $k$ copies of $a^\dagger$ and $n-k$ copies of $a$. For example,
        \begin{equation*}
            Y_{0, n} = a^n, \quad Y_{1, n-1} = a^\dagger a^{n-1} + a a^\dagger a^{n-2} + \cdots + a^{n-2} a^\dagger a + a^{n-1} a^\dagger.
        \end{equation*}
        The functions $z^k(-\overline z)^{n-k}$, $0\leq k\leq n$, are linearly independent. Evaluating them at suitable values of $z$ and taking linear combinations therefore gives
        \begin{equation*}
            \rho^{s/2} Y_{k, n-k} \rho^{(1-s)/2} \in \dom L^{(s)}, \quad \forall \, 0 \le k \le n.
        \end{equation*}
        Normal ordering, as described in Appendix~\ref{section-property-polynomials}, gives
        \begin{equation*}
            Y_{k, n-k} = \binom{n}{k} (a^\dagger)^k a^{n-k} + R_{k, n - k},
        \end{equation*}
        where $R_{k,n-k}$ has degree at most $n-2$: every use of $aa^\dagger=a^\dagger a+\1$ either preserves the leading normal-ordered monomial or lowers total degree by two. Assume inductively that $\rho^{s/2}X\rho^{(1-s)/2}\in\dom L^{(s)}$ for every polynomial $X$ of degree below $n$. Both embedded $Y_{k,n-k}$ and embedded $R_{k,n-k}$ then belong to the domain, and subtraction gives
        \begin{equation*}
            \rho^{s/2} (a^\dagger)^k a^{n-k} \rho^{(1-s)/2} \in \dom L^{(s)}.
        \end{equation*}
        Wick's theorem expresses every polynomial as a finite linear combination of normally ordered monomials, completing the induction.
    \end{proof}

    We also need the algebraic analogue of the preceding interchange argument. In the next lemma, Proposition~\ref{prop-L-and-calL-relation-Weyl-group}, and Theorem~\ref{theorem-L-and-calL-relation}, $\cL$ is interpreted algebraically through \eqref{eq-gaussian-gksl-generator}.

    \begin{lemma} \label{lemma-exchangeability-calL-derivative-r}
        Let $X$ be a polynomial in $a$ and $a^\dagger$. Then
        \begin{align}
            \frac{\md}{\md r} \rho^{s/2} \cL \left( W(rz) X \right) \rho^{(1-s)/2} &= \rho^{s/2} \cL \left( W(rz) (z a^\dagger - \overline{z} a) X \right) \rho^{(1-s)/2} \nonumber \\
            &= \rho^{s/2} \cL \left( \frac{\md}{\md r} W(rz) X \right) \rho^{(1-s)/2}. \label{eq-exchangeability-calL-derivative-r}
        \end{align}
    \end{lemma}
    \begin{proof}
        From \eqref{eq-gaussian-gksl-generator},
        \begin{align}
            &\quad \frac{\md}{\md r} \rho^{s/2} \cL (W(rz) X) \rho^{(1-s)/2} \nonumber \\
            &= \frac{\md}{\md r} ( \mi \rho^{s/2} H W(rz) X \rho^{(1-s)/2} - \mi \rho^{s/2} W(rz) X H \rho^{(1-s)/2} \nonumber \\
            &\quad -\frac{1}{2} \sum_\ell \rho^\frac{s}{2} \left( L_\ell^* L_\ell W(rz) X - 2 L_\ell^* W(rz) X L_\ell  + W(rz) X L_\ell^* L_\ell \right) \rho^\frac{1-s}{2} ) \label{eq-cL-trick-1}  \\
            &= \rho^{s/2} \cL \left( W(rz) (z a^\dagger - \overline{z} a) X \right) \rho^{(1-s)/2} \label{eq-cL-trick-2}  \\
            &= \rho^{s/2}  \cL \left( \frac{\md}{\md r} W(rz) X \right) \rho^{(1-s)/2}, \nonumber
        \end{align}
        which proves \eqref{eq-exchangeability-calL-derivative-r}. In passing from \eqref{eq-cL-trick-1} to \eqref{eq-cL-trick-2}, we used Lemma~\ref{lemma-exchangeability-X-Wrz} and Proposition~\ref{proposition-derivative-of-sandwiched-Wrz}. For example,
        \begin{equation*}
            \frac{\md}{\md r} \rho^{s/2} H W(rz) X \rho^{(1-s)/2} = \frac{\md}{\md r} \rho^{s/2} W(rz) \cdot W(rz)^{-1} H W(rz) \cdot X \rho^{(1-s)/2}
        \end{equation*}
        places the first term in the form covered by Proposition~\ref{proposition-derivative-of-sandwiched-Wrz}.
    \end{proof}

    \begin{proposition} \label{prop-L-and-calL-relation-Weyl-group}
        For every $n\in\N$,
        \begin{equation} \label{eq-L-and-calL-relation-Weyl-group}
            L^{(s)} \left( \rho^{s/2} W(rz) (z a^\dagger - \overline{z} a)^n \rho^{(1-s)/2} \right) = \rho^{s/2} \cL \left( W(rz) (z a^\dagger - \overline{z} a)^n \right) \rho^{(1-s)/2}.
        \end{equation}
    \end{proposition}
    \begin{proof}
        Equations~\eqref{eq-action-L-s-on-wrapped-Wz} and \eqref{eq-Yz}, Proposition~\ref{theo-exchangeability-L-derivative}, and Lemma~\ref{lemma-exchangeability-calL-derivative-r} give
        \begin{align}
            &\quad L^{(s)} \left( \rho^{s/2} W(rz) (z a^\dagger - \overline{z} a) \rho^{(1-s)/2} \right) \nonumber \\
            &= \frac{\md}{\md r} L^{(s)} \left( \rho^{s/2} W(rz) \rho^{(1-s)/2} \right) = \frac{\md}{\md r} \rho^{s/2} W(rz) Y(rz) \rho^{(1-s)/2} \nonumber \\
            &= \frac{\md}{\md r} \rho^{s/2} \cL (W(rz)) \rho^{(1-s)/2} = \rho^{s/2} \cL \left( W(rz) (z a^\dagger - \overline{z} a) \right) \rho^{(1-s)/2}. \label{eq-L-cL-Wrz-derivative}
        \end{align}
        The algebraic identity $\cL(W(rz))=W(rz)Y(rz)$ proves the case $n=0$, and \eqref{eq-L-cL-Wrz-derivative} proves the case $n=1$. Suppose that \eqref{eq-L-and-calL-relation-Weyl-group} holds at order $n$. Proposition~\ref{proposition-exchangeability-L-derivative-n-th-derivative} and Lemma~\ref{lemma-exchangeability-calL-derivative-r} give
        \begin{align*}
            &\quad L^{(s)} \left( \rho^{s/2} W(rz) (z a^\dagger - \overline{z} a)^{n+1} \rho^{(1-s)/2} \right) \\
            &= L^{(s)} \left( \frac{\md}{\md r} \rho^{s/2} W (rz) (z a^\dagger - \overline{z} a)^{n} \rho^{(1-s)/2} \right) \\
            &= \frac{\md}{\md r} L^{(s)} \left( \rho^{s/2} W (rz) (z a^\dagger - \overline{z} a)^{n} \rho^{(1-s)/2} \right) \\
            &= \frac{\md}{\md r} \rho^{s/2} \cL \left( W(rz) (z a^\dagger - \overline{z} a)^{n} \right) \rho^{(1-s)/2} \\
            &= \rho^{s/2} \cL \left( W(rz) (z a^\dagger - \overline{z} a)^{n+1} \right) \rho^{(1-s)/2}.
        \end{align*}
        Induction completes the proof.
    \end{proof}

    Proposition~\ref{prop-L-and-calL-relation-Weyl-group}, evaluated at $r=0$, now proves Theorem~\ref{theorem-L-and-calL-relation}: the linear-combination argument in the proof of Corollary~\ref{corollary-X-rho-one-half-L-domain} obtains every polynomial in $a$ and $a^\dagger$ from the powers $(za^\dagger-\overline z a)^n$.

        \section{Auxiliary Phase-Space and Polynomial Identities}

        \subsection{Real-Linear Phase-Space Conventions} \label{section-appendix-identification}

        Let $S:\C\to\C$ be real linear.

        \begin{lemma} \label{lemma-two-components-repr-real-linear-maps}
            There are unique $s_1,s_2\in\C$ such that
            \begin{equation} \label{eq-two-components-repr-real-linear-maps}
                S z = s_1 z + s_2 \overline{z}, \quad \forall z \in \C,
            \end{equation}
            where
            \begin{equation*}
                s_1 = \frac{S(1) - \mi S (\mi)}{2}, \quad s_2 = \frac{S(1) + \mi S(\mi)}{2}.
            \end{equation*}
        \end{lemma}

        We identify $S$ with a real matrix $\mathbf S$ and $z$ with a vector $\mathbf z\in\R^2$ by
        \begin{equation} \label{eq-identification-real-linear-operators}
            \mathbf{S} := \begin{pmatrix}
                \Re s_1 + \Re s_2 & \Im s_2 - \Im s_1 \\
                \Im s_1 + \Im s_2 & \Re s_1 - \Re s_2
            \end{pmatrix}, \quad \mathbf{z} := \begin{pmatrix}
                \Re z \\ \Im z
            \end{pmatrix}.
        \end{equation}
        Conversely, every real $2\times2$ matrix defines a real-linear operator on $\C$ through this identification.

        We will use the following identities without further comment.

        \begin{lemma}
            Let $S$ be real linear and let $y,z\in\C$. With the corresponding identifications $\mathbf S$, $\mathbf y$, and $\mathbf z$,
            \begin{equation*}
                \mathbf{S} \mathbf{z} = \begin{pmatrix}
                    \Re S z \\
                    \Im S z
                \end{pmatrix}, \quad \mathbf{y}^T \mathbf{z} = \Re \langle y, z \rangle.
            \end{equation*}
        \end{lemma}

        Moreover,
        \begin{equation*}
            \langle y, z \rangle = \mathbf{y}^T (\1 + \mi \mathbf{J}) \mathbf{z}, \quad \langle z, y \rangle = \mathbf{y}^T (\1 - \mi \mathbf{J}) \mathbf{z}.
        \end{equation*}

        Powers and exponentials are preserved by the identification:
        \begin{equation*}
            \mathbf{S}^n \mathbf{z} = \begin{pmatrix}
                \Re S^n z \\ \Im S^n z
            \end{pmatrix}, \quad \me^{\mathbf{S}} \mathbf{z} = \begin{pmatrix}
                \Re \me^S z \\
                \Im \me^S z
            \end{pmatrix}.
        \end{equation*}

        \subsection{Commutation Relations} \label{section-appendix-commutation-relations}

        We collect the commutation relations used throughout the paper. For the first-order expressions below, exchanging two factors changes the product only by a scalar multiple of the identity.

        Lemmas~\ref{lemma-commutation-relations-a-a-dagger-powers}--\ref{lemma-a-a-dagger-exponential-of-N} are standard, and we omit their proofs.

        \begin{lemma} \label{lemma-commutation-relations-a-a-dagger-powers}
            $ [(a^\dagger)^n, a] = - n (a^\dagger)^{n-1} $ and $ [a^n, a^\dagger] = n a^{n-1} $ for any $ n \ge 1 $.
        \end{lemma}

        \begin{lemma} \label{lemma-a-Wz-a-dagger-Wz-commutation-relations}
            $ [a, W(z)] = z W(z) $ and $ [a^\dagger, W(z)] = \overline{z} W(z) $.
        \end{lemma}

        \begin{lemma} \label{lemma-a-a-dagger-exponential-of-N}
            Let $ s \in \R $. Then,
            \begin{equation} \label{eq-a-a-dagger-exponential-of-N}
                \me^{-s N} a = \me^{s} a \me^{- s N}, \quad \me^{-s N} a^\dagger = \me^{-s} a^\dagger \me^{-s N}.
            \end{equation}
        \end{lemma}

        The commutator of a first-order polynomial with a power of another first-order polynomial has the following simple form.

        \begin{lemma} \label{lemma-commutation-relation-L-l-and-X-n}
            Let $ X $ and $ Y $ be first-order polynomials in $ a $ and $ a^\dagger $. Then,
            \begin{equation} \label{eq-commutation-relation-L-l-and-X-n}
                [X, Y^n] = n [X, Y] Y^{n-1}, \quad \forall n \ge 1.
            \end{equation}
            Here $[X,Y]$ is a scalar multiple of $\1$.
        \end{lemma}
        \begin{proof}
            The identity is immediate for $n=0,1$. Write $X=b_1a+b_2a^\dagger+b_3$ and $Y=c_1a+c_2a^\dagger+c_3$. Then
            \begin{equation*}
                [X, Y] = [ b_1 a + b_2 a^\dagger + b_3, c_1 a + c_2 a^\dagger + c_3 ] = (b_1 c_2 - b_2 c_1) \1.
            \end{equation*}
            If \eqref{eq-commutation-relation-L-l-and-X-n} holds for $n$, then
            \begin{equation*}
                [X, Y^{n+1}] = [X, Y^n Y] = Y^n [X, Y] + [X, Y^n] Y = (n+1) [X, Y] Y^n,
            \end{equation*}
            and induction completes the proof.
        \end{proof}

        \subsection{Polynomials in Creation and Annihilation Operators} \label{section-property-polynomials}

        We next recall Wick normal ordering, which rewrites arbitrary polynomials in a form adapted to the degree filtration.

        A product of creation and annihilation operators is \textit{Wick normal ordered} if every creation operator stands to the left of every annihilation operator, so it has the form $a^{\dagger n}a^m$ for $n,m\in\N$. We include the elementary normal-ordering argument for completeness.

        \begin{theorem}[Wick's Theorem] \label{theo-Wicks-theorem}
            Let $X$ be a monomial of degree $n$ in $a^\dagger$ and degree $m$ in $a$. Then $X$ is a finite linear combination of Wick-normal-ordered monomials of total degree at most $n+m$.
        \end{theorem}
        \begin{proof}
            The claim holds for $a$ and $a^\dagger$. Suppose inductively that a monomial $X$ can be written as
            \begin{equation*}
                X = \sum_{k} c_k a^{\dagger n_k} a^{m_k}
            \end{equation*}
            with $n_k+m_k\leq n+m$. The corresponding assertion for $Xa$ is immediate, while Lemma~\ref{lemma-commutation-relations-a-a-dagger-powers} gives
            \begin{align*}
                X a^\dagger &= \sum_k c_k a^{\dagger n_k} a^{m_k} a^\dagger = \sum_k c_k a^{\dagger n_k} (m_k a^{m_k - 1} + a^\dagger a^{m_k}) \\
                &= \sum_k m_k c_k a^{\dagger n_k} a^{m_k-1} + \sum_k c_k a^{\dagger (n_k + 1)} a^{m_k},
            \end{align*}
            which is a finite sum of normally ordered monomials of degree at most $n+m+1$. Induction on the word length completes the proof.
        \end{proof}

        We will also use the following change of first-order variables.

        \begin{corollary} \label{corollary-wick-theorem-variant}
            Let $X$ and $Y$ be linearly independent first-order polynomials without constant terms. Every monomial $W$ in $a$ and $a^\dagger$ has a finite expansion
            \begin{equation} \label{eq-finite-sum-normal-order-terms-X-Y}
                W=\sum_k c_kX^{n_k}Y^{m_k},
            \end{equation}
            where $c_k\in\C$ and $n_k,m_k\in\N$. We use the convention $X^0Y^0=\1$.
        \end{corollary}
        \begin{proof}
            Linear independence allows us to write $a=b_1X+b_2Y$ and $a^\dagger=c_1X+c_2Y$, so \eqref{eq-finite-sum-normal-order-terms-X-Y} holds for words of length one. Suppose it holds for a monomial $W$. Lemma~\ref{lemma-commutation-relation-L-l-and-X-n} gives
            \begin{align*}
                W a &= \sum_k c_k X^{n_k} Y^{m_k} a = \sum_k c_k X^{n_k} Y^{m_k} (b_1 X + b_2 Y) \\
                &= b_1 \sum_k c_k X^{n_k} Y^{m_k} X + b_2 \sum_k c_k X^{n_k} Y^{m_k + 1} \\
                &= b_1 \sum_k c_k X^{n_k} \left( [Y^{m_k}, X] + X Y^{m_k} \right) + b_2 \sum_k c_k X^{n_k} Y^{m_k + 1} \\
                &= b_1 [Y, X] \sum_k m_k c_k X^{n_k} Y^{m_k-1} \\
                &\quad + b_1 \sum_k c_k X^{n_k+1} Y^{m_k} + b_2 \sum_k c_k X^{n_k} Y^{m_k + 1}.
            \end{align*}
            Since $[Y,X]$ is a scalar multiple of $\1$, this is again an expansion of the required form. The argument for $Wa^\dagger$ is identical, and induction completes the proof.
        \end{proof}

        Both statements extend by linearity to arbitrary polynomials. We call an operator of the form $\rho^{s/2}X\rho^{(1-s)/2}$, with $X$ a polynomial in $a$ and $a^\dagger$, an \emph{embedded polynomial}.

        \begin{lemma} \label{lemma-polynomial-rho-hilbert-schmidt}
            Every embedded polynomial $\rho^{s/2}X\rho^{(1-s)/2}$ is Hilbert--Schmidt.
        \end{lemma}
        \begin{proof}
            By linearity, it suffices to take $X$ to be a monomial of degree $m$ in $a^\dagger$ and degree $n$ in $a$. The operator $\rho^{s/2}X\rho^{(1-s)/2}$ is initially defined on $\cD$. There are $\ell\leq n$ and a polynomial $P$ of degree $m+n$ such that
            \begin{equation} \label{eq-P-a-dagger-a-e}
                X e_k = 0, \quad \forall k < \ell; \quad X e_k = \sqrt{P (k)} e_{k+m-n}, \quad \forall k \ge \ell.
            \end{equation}
            Let $x=\sum_{k\geq0}x_ke_k\in\cD$ with $\norm{x}=1$. Then $(x_k)$ is finitely supported, $\sum_k\abs{x_k}^2=1$, and \eqref{eq-P-a-dagger-a-e} gives
            \begin{align*}
                &\quad \norm{\rho^{s/2} X \rho^{(1-s)/2} x} \\
                &= \sqrt{1 - \me^{-\beta}} \norm{\me^{-s \beta N / 2} X \me^{-(1-s)\beta N /2} \left( \sum_{k = 0}^\infty x_k e_k \right) } \\
                &= \sqrt{1 - \me^{-\beta}} \norm{ \sum_{k=\ell}^\infty x_k \me^{-s \beta ( k + m - n) / 2} \sqrt{P(k)} \me^{-(1-s) \beta k /2} e_{k+m-n}} \\
                &\le \sqrt{1 - \me^{-\beta}} \, \me^{-s \beta (m-n)/2} \sum_{k = \ell}^\infty \abs{\me^{- \beta k / 2} \sqrt{P(k)} \cdot x_k} \\
                &\le \sqrt{1 - \me^{-\beta}} \, \me^{-s \beta (m-n)/2} \left( \sum_{k=0}^\infty \me^{-\beta k} P (k) \right)^{1/2} < +\infty,
            \end{align*}
            so the operator extends boundedly to $\mathsf h$. Moreover,
            \begin{equation*}
                \sum_{k = 0}^\infty \norm{\rho^{s/2} X \rho^{(1-s)/2} e_k }^2 = (1 - \me^{-\beta}) \me^{-s \beta (m - n)} \sum_{k=\ell}^\infty \abs{\me^{-\beta k} P (k)} < +\infty,
            \end{equation*}
            which proves that the extension is Hilbert--Schmidt.
        \end{proof}

        \begin{remark} \label{remark-X-rho-s-trace-class}
            Let $X$ be a polynomial in $a$ and $a^\dagger$ and let $0<s\leq1$. Since $\rho=(1-\me^{-\beta})\me^{-\beta N}$, every positive power of $\rho$ is trace class. The proof of Lemma~\ref{lemma-polynomial-rho-hilbert-schmidt} also shows that $X\rho^{s/2}$ and $\rho^{s/2}X$ are Hilbert--Schmidt. Hence $X\rho^s=(X\rho^{s/2})\rho^{s/2}$ and $\rho^sX=\rho^{s/2}(\rho^{s/2}X)$ are trace class.
        \end{remark}

        \begin{lemma} \label{lemma-X-rho-rho-hat-X}
            Let $X$ be a polynomial in $a$ and $a^\dagger$, and let $s\in\R$. On the common invariant core $\cD$, the operator $\rho^sX\rho^{-s}$ agrees with a polynomial in $a$ and $a^\dagger$ of the same degree.
        \end{lemma}
        \begin{proof}
            It suffices to consider a monomial $X$. The case $s=0$ is immediate. For $s>0$, we construct a polynomial $\widehat X$ of the same degree such that, on $\cD$,
            \begin{equation} \label{eq-X-rho-rho-hat-X}
                \rho^{s} X = \widehat{X} \rho^{s}.
            \end{equation}
            Assume that \eqref{eq-X-rho-rho-hat-X} holds for a certain $ X $. Then, by Lemma \ref{lemma-a-a-dagger-exponential-of-N}, we have
            \begin{equation*}
                \rho^{s} X a = \widehat{X} \rho^{s} a = \me^{s \beta} \widehat{X} a \rho^{s}, \quad \rho^{s} X a^\dagger =  \widehat{X} \rho^{s} a^\dagger = \me^{-s \beta}   \widehat{X} a^\dagger \rho^{s}
            \end{equation*}
            on $\cD$. Induction on the word length proves \eqref{eq-X-rho-rho-hat-X}. Composing on the right with $\rho^{-s}$ on $\cD$ gives
            \begin{equation*}
                \widehat{X} = \rho^s X \rho^{-s}.
            \end{equation*}
            The case $s<0$ is identical.
        \end{proof}

        The Weyl commutation relations also give the following conjugation formula.

        \begin{lemma} \label{lemma-exchangeability-X-Wrz}
            Let $X$ be a polynomial of degree $d$ in $a$ and $a^\dagger$, and let $r\in\R$ and $z\in\C$. Then $W(-rz)XW(rz)$ is a polynomial of degree $d$ whose coefficients depend polynomially on $r$.
        \end{lemma}
        \begin{proof}
            Write $X$ in Wick-normal-ordered form:
            \begin{equation*}
                X = \sum_{m+n \le d} c_{m, n} (a^\dagger)^m a^n.
            \end{equation*}
            Lemma~\ref{lemma-a-Wz-a-dagger-Wz-commutation-relations} gives
            \begin{equation*}
                W(-rz) \, a \, W(rz) = a + rz \1, \quad W(-rz) \, a^\dagger \, W(rz) = a^\dagger + r \overline{z} \1.
            \end{equation*}
            Since conjugation is multiplicative,
            \begin{align*}
                W(-rz) X W(rz) &= \sum_{m+n \le d} c_{m,n} \left( W(-rz) \, a^\dagger \, W(rz) \right)^m \left( W(-rz) \, a \, W(rz) \right)^n \\
                &= \sum_{m+n \le d} c_{m,n} (a^\dagger + r \overline{z} \1 )^m ( a + r z \1 )^n \\
                &= \sum_{m+n \le d} c_{m, n} \sum_{j=0}^m \sum_{k=0}^n \binom{m}{j} \binom{n}{k} (r \overline{z})^{m - j} (r z)^{n - k} (a^\dagger)^j a^k.
            \end{align*}
            For fixed $j$ and $k$, the coefficient of $(a^\dagger)^ja^k$ is
            \begin{equation*}
                \sum_{ \substack{ m \ge j, n \ge k \\ m+n \le d}} c_{m, n} \binom{m}{j} \binom{n}{k} \overline{z}^{m - j} z^{n - k} r^{m+n-j-k},
            \end{equation*}
            a polynomial in $r$ of degree at most $d-j-k$. If $j+k=d$, the only possible indices are $m=j$ and $n=k$, so the homogeneous component of degree $d$ is unchanged. Hence $W(-rz)XW(rz)$ has the same degree as $X$.
        \end{proof}

        Recall from \eqref{eq-introduction-polynomial-spaces} and the definition following it that, for $r\in\N$, $\cX_r$ is the space of $s$-embedded polynomials of degree at most $r$, and $\cX=\bigcup_{r\in\N}\cX_r$. Each $\cX_r$ is finite-dimensional. As a subspace of $\cB_2(\mathsf h)$, it is independent of $s$: Lemma~\ref{lemma-X-rho-rho-hat-X} shows that changing the embedding parameter amounts to modular conjugation, which scales $a$ and $a^\dagger$ and preserves degree.

        \begin{lemma} \label{lemma-density-of-domain-L-star-plus-L}
            The space $\cX$ is dense in $\cB_2(\mathsf h)$.
        \end{lemma}
        \begin{proof}
            Let $v\in\cB_2(\mathsf h)$ satisfy $v\perp\cX$. Lemma~\ref{lemma-Wrz-rho-1/2-power-series} implies $v\perp\rho^{s/2}W(z)\rho^{(1-s)/2}$ for every $z\in\C$. Hence the normal functional
            \begin{equation*}
                W(z) \mapsto \Tr( \rho^{(1-s)/2} v^* \rho^{s/2} W(z) )
            \end{equation*}
            vanishes on all Weyl operators; the sandwiched operator $\rho^{(1-s)/2}v^*\rho^{s/2}$ is trace class by H\"older's inequality. The linear span of the Weyl operators is $\sigma$-weakly dense in $\cB(\mathsf h)$ \cite[Corollary~2.4.15]{bratteli2012operator}. Therefore $\rho^{(1-s)/2}v^*\rho^{s/2}=0$, and faithfulness of $\rho$ gives $v=0$. Thus $\cX$ is dense.
        \end{proof}

        \section{Diagonal Standardization of the Invariant State} \label{section-non-diagonal-gaussian-states}

        Let $\rho$ be a Gaussian state with mean $\omega\in\C$ and covariance operator $S$, and let $\mathbf S$ be the real matrix representing $S$. We call $\mathbf S$ the \emph{covariance matrix} of $\rho$. Williamson's theorem \cite{williamson1936algebraic} gives a number $\nu\geq1$ and a symplectic matrix $\mathbf G$ such that
        \begin{equation} \label{eq-symplectic-diagonalization-covariance-matrix}
            \mathbf{S} = \mathbf{G}^T \begin{pmatrix} \nu & 0 \\ 0 & \nu \end{pmatrix} \mathbf{G}.
        \end{equation}

        Equivalently, every faithful one-mode Gaussian state is unitarily equivalent to a thermal state of inverse temperature $\beta>0$:
        \begin{equation} \label{eq-gaussian-state-diagonalization}
            \rho = W(\omega) B(\mathbf{M}) (1 - \me^{-\beta}) \me^{-\beta N} B(\mathbf{M})^{-1} W(\omega)^{-1}.
        \end{equation}
        Here $\mathbf M:=\mathbf G^{-1}$, $B(\mathbf M)$ is the unitary Bogoliubov transformation associated with $\mathbf M$ \cite{bruneau2007bogoliubov,parthasarathy2021pedagogical}, and $\coth(\beta/2)=\nu$. Moreover, $B(\mathbf M)^{-1}=B(\mathbf M^{-1})$, and the Bogoliubov transformation preserves the canonical commutation relation:
        \begin{equation*}
            [ B (\mathbf{M}) a B (\mathbf{M})^*, B (\mathbf{M}) a^\dagger B (\mathbf{M})^* ] = [a, a^\dagger].
        \end{equation*}

        Let $M$ be the real-linear operator represented by $\mathbf M$, write $Mz=m_1z+m_2\overline z$, and let $\boldsymbol\omega$ be the real vector representing $\omega$.

        For a Gaussian QMS with invariant state $\rho$, the invariance equations characterize the mean and covariance in terms of $\Omega$, $\kappa$, $\zeta$, $u_\ell$, and $v_\ell$ \cite[Theorem~8]{agredo2021gaussian}. The next lemma specializes those equations to a state diagonal in the number basis.

        \begin{lemma} \label{lemma-diagonal-rho-parameters}
            The normal invariant Gaussian state $\rho$ is diagonal in the number basis if and only if \eqref{eq-diagonal-rho-parameters} and \eqref{eq-diagonal-rho-parameters-beta} hold.
        \end{lemma}
        \begin{proof}
            A Gaussian state is diagonal in the number basis precisely when it is invariant under every phase rotation $\me^{\mi\theta N}$. Rotation invariance of its quantum characteristic function is equivalent to
            \begin{equation} \label{eq-equivalent-condition-rho-diagonal}
                \mathbf{\omega} = 0, \quad \mathbf{S} = \nu \1 .
            \end{equation}
            Faithfulness gives $\nu>1$, so there is a unique $\beta>0$ with $\nu=\coth(\beta/2)$, and the state is $(1-\me^{-\beta})\me^{-\beta N}$. The invariance equations from \cite[proof of Theorem~8]{agredo2021gaussian} are
            \begin{equation} \label{eq-lyapunov-equations}
                \boldsymbol{\zeta} = \mathbf{Z}^T \boldsymbol{\omega}, \quad \mathbf{Z}^T \mathbf{S} + \mathbf{S} \mathbf{Z} + \mathbf{C} = 0.
            \end{equation}
            Substituting \eqref{eq-equivalent-condition-rho-diagonal} and the definitions of $Z$ and $C$ from \eqref{eq-introduction-drift-operator} and \eqref{eq-introduction-diffusion-operator} gives
            \begin{equation*}
                \zeta = 0, \quad \kappa = \mi \tanh \left( \frac{\beta}{2} \right) \left( \sum_\ell u_\ell v_\ell \right) , \quad \me^{-\beta} \sum_\ell \abs{v_\ell}^2 = \sum_\ell \abs{u_\ell}^2,
            \end{equation*}
            which is equivalent to \eqref{eq-diagonal-rho-parameters}--\eqref{eq-diagonal-rho-parameters-beta}.
        \end{proof}

        We use the standardization from \cite[Section~2]{girotti2024gaussian}. Define
        \begin{equation} \label{definition-standard-QMS}
            \cS_t (X) := U^* \cT_t (U X U^*) U, \quad \forall X \in \cB(\mathsf{h}),
        \end{equation}
        where
        \begin{equation*}
            U := W(\omega) B(\mathbf{M}).
        \end{equation*}
        The action of $(\cS_t)_{t\geq0}$ on Weyl operators gives the following standardization result \cite{girotti2024gaussian}.

        \begin{lemma} \label{lemma-drift-matrix-of-standardized-QMS}
            The semigroup $(\cS_t)_{t\geq0}$ is a Gaussian QMS with diagonal invariant state $\rho_0:=(1-\me^{-\beta})\me^{-\beta N}$ and drift matrix $\mathbf M^{-1}\mathbf Z\mathbf M$.
        \end{lemma}

        We call $(\cS_t)_{t\geq0}$ the \emph{standardized QMS} associated with $(\cT_t)_{t\geq0}$; its invariant state has zero mean and scalar covariance.

        Since $\rho=U\rho_0U^*$, the map
        \begin{equation*}
            \cU : \cB_2 (\mathsf{h}) \rightarrow \cB_2 (\mathsf{h}), \quad \cU (Y) = U^* Y U, \quad Y \in \cB_2 (\mathsf{h})
        \end{equation*}
        is unitary.

        \begin{proposition} \label{proposition-induced-semigroup-of-standardized-qms}
            We have $\cU T_t^{(s)}\cU^*=S_t^{(s)}$, where $(S_t^{(s)})_{t\geq0}$ is the semigroup induced by $(\cS_t)_{t\geq0}$ with respect to $\rho_0$.
        \end{proposition}
        \begin{proof}
            For $X\in\cB(\mathsf h)$, direct computation gives
            \begin{align*}
                \cU T_t^{(s)} \cU^* ( \rho_0^{s/2} X \rho^{(1-s)/2}_0 ) &= \cU T_t^{(s)} (  U U^* \rho^{s/2} U X U^* \rho^{(1-s)/2} U U^*  ) \\
                &= U^* \rho^{s/2} U U^* \cT_t ( U X U^*) U U^* \rho^{(1-s)/2} U \\
                &= \rho_0^{s/2} \cS_t (X) \rho^{(1-s)/2}_0 = S_t^{(s)} (\rho_0^{s/2} X \rho_0^{(1-s)/2}).
            \end{align*}
            Density of $\rho_0^{s/2}\cB(\mathsf h)\rho_0^{(1-s)/2}$ in $\cB_2(\mathsf h)$ proves the identity on the entire Hilbert--Schmidt space.
        \end{proof}

        Consequently, the induced generators are unitarily equivalent, as are their adjoints and symmetric parts. This equivalence preserves spectra, point spectra, multiplicities, compactness, compact resolvent, and the numerical spectral gap. The standardized and original drift matrices are similar and therefore have the same eigenvalues.

        \begin{lemma} \label{lemma-elementary-inequalities-spectral-gap}
            The invariant-state parameters satisfy
            \begin{equation} \label{eq-kms-spectral-gap-strict-positive}
                \gamma>\abs{\kappa},
            \end{equation}
            and
            \begin{equation} \label{eq-gns-spectral-gap-positive}
                \gamma\geq\cosh(\beta/2)\abs{\kappa}.
            \end{equation}
            Equivalently,
            \begin{equation} \label{eq-gns-spectral-gap-positive-equivalent-form}
                \frac{\gamma+\abs{\kappa}}{\gamma-\abs{\kappa}}
                \leq\frac{\cosh(\beta/2)+1}{\cosh(\beta/2)-1}
                =\coth^2(\beta/4).
            \end{equation}
        \end{lemma}
        \begin{proof}
            The diagonal invariant-state relations give
            \[
                \sum_\ell\abs{u_\ell}^2=\me^{-\beta}\sum_\ell\abs{v_\ell}^2,
                \qquad
                \kappa=\mi\tanh(\beta/2)\sum_\ell u_\ell v_\ell.
            \]
            The Cauchy--Schwarz inequality therefore yields
            \begin{align*}
                \abs{\kappa}\cosh(\beta/2)
                &=\sinh(\beta/2)\abs{\sum_\ell u_\ell v_\ell}\\
                &\leq\sinh(\beta/2)
                \left(\sum_\ell\abs{u_\ell}^2\right)^{1/2}
                \left(\sum_\ell\abs{v_\ell}^2\right)^{1/2}\\
                &=\frac{1-\me^{-\beta}}{2}\sum_\ell\abs{v_\ell}^2=\gamma.
            \end{align*}
            This proves \eqref{eq-gns-spectral-gap-positive}. Since $\beta>0$, it implies \eqref{eq-kms-spectral-gap-strict-positive} when $\kappa\neq0$; for $\kappa=0$, the latter follows from $\gamma>0$. Rearranging \eqref{eq-gns-spectral-gap-positive} gives \eqref{eq-gns-spectral-gap-positive-equivalent-form}.
        \end{proof}

        \begin{remark} \label{remark-complete-positivity-Lyapunov}
        Let $\mathbf C$ be the matrix representing the diffusion operator $C$ and set
        \begin{equation} \label{eq-definition-C-Z}
            \mathbf C_Z:=\mathbf C-\mi(\mathbf Z^T\mathbf J+\mathbf J\mathbf Z).
        \end{equation}
        \cite[Proposition~6]{fagnola2025spectral} shows that $\mathbf C_Z\geq0$, and complex conjugation therefore gives $\overline{\mathbf C_Z}\geq0$. Combining the invariant-state Lyapunov equation with \eqref{eq-definition-C-Z} yields
        \begin{align}
            \mathbf Z^T(\mathbf S+\mi\mathbf J)+(\mathbf S+\mi\mathbf J)\mathbf Z
            &=-\mathbf C_Z, \label{eq-Lyapunov-plus-iJ}\\
            \mathbf Z^T(\mathbf S-\mi\mathbf J)+(\mathbf S-\mi\mathbf J)\mathbf Z
            &=-\overline{\mathbf C_Z}. \label{eq-Lyapunov-minus-iJ}
        \end{align}
        These identities are the matrix form of the one-mode estimate $\gamma\geq\cosh(\beta/2)\abs{\kappa}$ and supply the positivity arguments below.
        \end{remark}

    \section{Characteristic-Function Calculations}
    \label{section-characteristic-function-calculations}

    This appendix contains the Gaussian calculations used in Theorem~\ref{theorem-characteristic-function}. For a complex matrix $X$, $\norm{X}$ denotes the operator norm induced by the Euclidean norm, and $\Re X:=(X+X^*)/2$ denotes its Hermitian part. For the complex symmetric matrices $\mathbf R_t$ and $\mathbf G_t$, the Hermitian part agrees with the entrywise real part. When $\Re\mathbf B_t$ appears as an off-diagonal block of $\Re\mathbf G_t$, it denotes the entrywise real part of $\mathbf B_t$.

    The two longer calculations from Section~\ref{section-immediate-compactness} are recorded here so that the main line of the characteristic-function argument remains visible.

    \proofweightedweylpairs

    \proofcharacteristickernel

    We begin by simplifying the matrices introduced in the proof of Theorem~\ref{theorem-characteristic-function}. The matrix $\mathbf C_t-\mathbf B_t\mathbf A_t^{-1}\mathbf B_t^T$ is the Schur complement of the block $\mathbf A_t$ in $\mathbf G_t$.

    \begin{lemma} \label{lemma-simplification-A-t-B-t-C-t}
        The quantities $\mathbf A_t$, $\mathbf B_t$, $\mathbf C_t$, and $c_t$ defined in \eqref{eq-original-A-t}--\eqref{eq-original-c-t} have the forms \eqref{eq-A-t-by-M-t}--\eqref{eq-final-c-t}, respectively. Moreover,
        \begin{align}
            \mathbf{A}_t^{-1} \mathbf{B}_t^T &= - 2\mathbf{\Theta}_{1-s} (\me^{t \mathbf{Z}} + \me^{-t \mathbf{Z}^T} )^{-1} \mathbf{\Theta}_{1-s}^T, \label{eq-simplified-A-t-inv-B-t-T} \\
            \mathbf{C}_t - \mathbf{B}_t \mathbf{A}_t^{-1} \mathbf{B}_t^T &= \coth(\beta/4)\mathbf{\Theta}_{1-s} ( \1 -\me^{t \mathbf{Z}} \me^{t \mathbf{Z}^T} ) ( \1 + \me^{t \mathbf{Z}} \me^{t \mathbf{Z}^T} )^{-1} \mathbf{\Theta}_{1-s}^T. \label{eq-simplified-schur-complement}
         \end{align}
    \end{lemma}
    \begin{proof}
        Recall from \eqref{definition-Theta-s}, \eqref{eq-definition-J-s} and \eqref{eq-definition-R-t} that
        \begin{align*}
            \mathbf{\Theta}_s &= \cosh\!\left(\frac{1-2s}{4}\beta\right)\1-\sinh\!\left(\frac{1-2s}{4}\beta\right)\mi\mathbf J, \\
            \mathbf{J}_s &= \me^{\frac{s}{2} \beta} (\1 + \mi \mathbf{J}) - \me^{\frac{1-s}{2} \beta} (\1 - \mi \mathbf{J}), \\
            \mathbf{R}_t &= \mathbf{J}_s^T ( \frac{1 - \me^{-\beta/2}}{2} \1 - (\mathbf{Q}_t + \1)^{-1} ) \mathbf{J}_s.
        \end{align*}
        We use the following two identities repeatedly:
        \begin{align}
            (\1 + \mi \mathbf{J})\mathbf{\Theta}_{\frac{1}{2} - s}^T &= \me^{\frac{s}{2} \beta} (\1 + \mi \mathbf{J}), \label{eq-Theta-identity-1} \\
            \sinh(\beta/4) ( \mathbf{\Sigma} - \mi \mathbf{J} )\mathbf{\Theta}_{\frac{1}{2} - s}^T &= \mathbf{\Theta}_s. \label{eq-Theta-identity-2}
        \end{align}
        The first is the transpose of \eqref{eq-identity-exp-s-beta-one-half-Theta}.

        The recurring matrix products can be expressed in terms of the symmetric matrices $\me^{t\mathbf Z^T}\me^{t\mathbf Z}$ and $\me^{t\mathbf Z}\me^{t\mathbf Z^T}$:
        \begin{align*}
            \left( \mathbf{Q}_t + \1 \right)^{-1} &= \left( (\coth(\beta/2) + 1) \1 + \csch(\beta/2) \me^{t \mathbf{Z}^T} \me^{t \mathbf{Z}} \right)^{-1} \\
            &= \sinh(\beta/2) (\me^{\beta/2} + \me^{t \mathbf{Z}^T} \me^{t \mathbf{Z}} )^{-1}, \\
            \mathbf{J}_s \mathbf{R}_t^{-1} \mathbf{J}_s^T &= \mathbf{J}_s \mathbf{J}_s^{-1} \left( \frac{1- \me^{-\beta/2}}{2} \1 - (\mathbf{Q}_t + \1)^{-1} \right)^{-1} (\mathbf{J}_s^T)^{-1} \mathbf{J}_s^T \\
            &= \frac{2}{1 - \me^{-\beta/2}} \left( \1 - (\me^{\beta/2} + 1) ( \me^{t \mathbf{Z}^T} \me^{t \mathbf{Z}} + \me^{\beta/2} )^{-1} \right)^{-1} \\
            &= \frac{2}{1 - \me^{-\beta/2}} ( \me^{t \mathbf{Z}^T} \me^{t \mathbf{Z}} + \me^{\beta/2} ) ( \me^{t \mathbf{Z}^T} \me^{t \mathbf{Z}} - \1 )^{-1},
        \end{align*}
        Consequently,
        \begin{equation} \label{eq-Q-t-1-J-R-J}
            (\mathbf{Q}_t + \1)^{-1} \mathbf{J}_s \mathbf{R}_t^{-1} \mathbf{J}_s^T = ( 1 + \me^{\beta/2} ) ( \me^{t \mathbf{Z}^T} \me^{t \mathbf{Z}} - \1 )^{-1},
        \end{equation}
        and
        \begin{equation} \label{eq-Q-t-1-plus-Q-t-1-J-R-J-Q-t-1}
            (\mathbf{Q}_t + \1)^{-1} + (\mathbf{Q}_t + \1)^{-1} \mathbf{J}_s \mathbf{R}_t^{-1} \mathbf{J}_s^T (\mathbf{Q}_t + \1)^{-1} = \sinh(\beta/2) ( \me^{t \mathbf{Z}^T} \me^{t \mathbf{Z}} - \1 )^{-1}.
        \end{equation}
        Using \eqref{eq-Theta-identity-1} and \eqref{eq-Theta-identity-2} in \eqref{eq-original-B-t}, we obtain

        \begin{align*}
            \mathbf{B}_t &= \csch(\beta/4) \mathbf{\Theta}_{1 - s} \me^{t \mathbf{Z}} ( \mathbf{Q}_t + \1 )^{-1} \left( - \me^{\frac{s}{2} \beta} (\1 + \mi \mathbf{J}) \right. \nonumber \\
            &\quad \left. + \mathbf{J}_s \mathbf{R}_t^{-1} \mathbf{J}_s^T (\mathbf{Q}_t + \1)^{-1} (\mathbf{Q}_t - \mi \mathbf{J}) \mathbf{\Theta}_{\frac{1}{2}-s}^T \right), \\
            &= - \csch(\beta/4) \mathbf{\Theta}_{1-s} \me^{t \mathbf{Z}} \left( ( (\mathbf{Q}_t + \1)^{-1} + (\mathbf{Q}_t + \1)^{-1} \mathbf{J}_s \mathbf{R}_t^{-1} \mathbf{J}_s^T (\mathbf{Q}_t + \1)^{-1} ) \right. \\
            &\qquad \left. \cdot (\1 + \mi \mathbf{J}) - (\mathbf{Q}_t + \1)^{-1} \mathbf{J}_s \mathbf{R}_t^{-1} \mathbf{J}_s^T \right) \mathbf{\Theta}_{\frac{1}{2}-s}^T \\
            &= \sinh(\beta/2) \csch(\beta/4) \mathbf{\Theta}_{1-s} ( \me^{t \mathbf{Z}^T} - \me^{- t \mathbf{Z}} )^{-1} (\mathbf{\Sigma} - \mi \mathbf{J}) \mathbf{\Theta}_{\frac{1}{2} - s}^T \\
            &=2\coth(\beta/4)\mathbf{\Theta}_{1-s} ( \me^{t \mathbf{Z}^T } - \me^{- t \mathbf{Z}} )^{-1}\mathbf{\Theta}_{1-s}^T.
        \end{align*}

        Similarly, \eqref{eq-original-C-t} gives
        \begin{align*}
            \mathbf{C}_t &= \mathbf{\Sigma} - \csch(\beta/4) \mathbf{\Theta}_{1-s} \me^{t \mathbf{Z}} ( \mathbf{Q}_t + \1 )^{-1} \me^{t \mathbf{Z}^T} \csch(\beta/4) \mathbf{\Theta}_{1-s}^T \nonumber \\
            &\quad - \csch(\beta/4) \mathbf{\Theta}_{1-s} \me^{t \mathbf{Z}} (\mathbf{Q}_t + \1)^{-1} \mathbf{J}_s \mathbf{R}_t^{-1} \mathbf{J}_s^T ( \mathbf{Q}_t + \1 )^{-1} \me^{t \mathbf{Z}^T} \csch(\beta/4) \mathbf{\Theta}_{1-s}^T \\
            &= \mathbf{\Sigma} - \sinh(\beta/2) \csch(\beta/4) \mathbf{\Theta}_{1-s} \me^{t \mathbf{Z}} ( \me^{t \mathbf{Z}^T} \me^{t \mathbf{Z}} - \1 )^{-1} \me^{t \mathbf{Z}^T} \csch(\beta/4) \mathbf{\Theta}_{1-s}^T \\
            &= \sinh(\beta/2) \csch(\beta/4) \mathbf{\Theta}_{1-s} \left( \frac{1}{2} \1 + \me^{t \mathbf{Z}} \me^{t \mathbf{Z}^T} (\1 - \me^{t \mathbf{Z}} \me^{t \mathbf{Z}^T} )^{-1} \right) \csch(\beta/4) \mathbf{\Theta}_{1-s}^T \\
            &=\coth(\beta/4)\mathbf{\Theta}_{1-s} (\1 + \me^{t \mathbf{Z}} \me^{t \mathbf{Z}^T} ) (\1 - \me^{t \mathbf{Z}} \me^{t \mathbf{Z}^T} )^{-1}\mathbf{\Theta}_{1-s}^T.
        \end{align*}
        From \eqref{eq-original-c-t},
        \begin{align*}
            c_t &= \frac{4 \sqrt{\coth(\beta/4)}}{(1-\me^{-\beta})^{\frac{1}{2}} \sqrt{\det(\mathbf{Q}_t + \1) \det \mathbf{R}_t}} \\
            &= \frac{4}{(1 - \me^{-\beta/2})  \abs{ \det \mathbf{J}_s}  \sqrt{ \det(\mathbf{Q}_t + \1) \det ( \frac{1 - \me^{-\beta/2}}{2}\1 - (\mathbf{Q}_t + \1)^{-1} ) } } \nonumber \\
            &= \left( (\me^{\beta/2} - 1) \sqrt{\det \left( (1 - \me^{-\beta/2})(\mathbf{Q}_t + \1) / 2 - \1 \right)} \right)^{-1} \nonumber \\
            &= \coth(\beta/4) \left( \det (\1 - \me^{t \mathbf{Z}^T} \me^{t \mathbf{Z}} ) \right)^{-1/2}.
        \end{align*}
        It remains to simplify $\mathbf A_t$. From \eqref{eq-original-A-t},

        \begin{align*}
            \mathbf{A}_t &= \1 - \me^{s \beta} (\1 - \mi \mathbf{J}) ( \mathbf{Q}_t + \1 )^{-1} (\1 + \mi \mathbf{J}) - \mathbf{\Theta}_{\frac{1}{2} - s} (\mathbf{Q}_t + \mi \mathbf{J}) \\
            &\quad \cdot (\mathbf{Q}_t + \1)^{-1} \mathbf{J}_s \mathbf{R}_t^{-1} \mathbf{J}_s^T ( \mathbf{Q}_t + \1 )^{-1} (\mathbf{Q}_t - \mi \mathbf{J}) \mathbf{\Theta}_{\frac{1}{2} - s}^T \\
            &= \mathbf{\Theta}_{\frac{1}{2} - s} \left( \1 - (\1 - \mi \mathbf{J}) (\mathbf{Q}_t + \1)^{-1} (\1 + \mi \mathbf{J}) \right. \\
            &\qquad \left. - ( \mathbf{Q}_t + \mi \mathbf{J} ) ( \mathbf{Q}_t + \1 )^{-1} \mathbf{J}_s \mathbf{R}_t^{-1} \mathbf{J}_s^T (\mathbf{Q}_t + \1)^{-1} (\mathbf{Q}_t - \mi \mathbf{J}) \right) \mathbf{\Theta}_{\frac{1}{2} - s}^T.
        \end{align*}

        The middle factor simplifies as follows:
        \begin{align*}
            &\quad \1 - (\1 - \mi \mathbf{J}) (\mathbf{Q}_t + \1)^{-1} (\1 + \mi \mathbf{J}) \\
            &\qquad - ( \mathbf{Q}_t + \mi \mathbf{J} ) ( \mathbf{Q}_t + \1 )^{-1} \mathbf{J}_s \mathbf{R}_t^{-1} \mathbf{J}_s^T (\mathbf{Q}_t + \1)^{-1} (\mathbf{Q}_t - \mi \mathbf{J}) \\
            &= \1 - \mathbf{J}_s \mathbf{R}_t^{-1} \mathbf{J}_s^T - (\1 - \mi \mathbf{J}) (\mathbf{Q}_t + \1)^{-1} (\1 + \mi \mathbf{J}) \\
            &\quad - (\1 - \mi \mathbf{J}) (\mathbf{Q}_t + \1)^{-1} \mathbf{J}_s \mathbf{R}_t^{-1} \mathbf{J}_s^T ( \mathbf{Q}_t + \1 )^{-1} (\1 + \mi \mathbf{J}) \\
            &\quad + (\1 - \mi \mathbf{J}) (\mathbf{Q}_t + \1)^{-1} \mathbf{J}_s \mathbf{R}_t^{-1} \mathbf{J}_s^T + \mathbf{J}_s \mathbf{R}_t^{-1} \mathbf{J}_s^{T} (\mathbf{Q}_t + \1)^{-1} (\1 + \mi \mathbf{J}) \\
            &= \1 + (2 (\mathbf{Q}_t + \1)^{-1} - \1 ) \mathbf{J}_s \mathbf{R}_t^{-1} \mathbf{J}_s^T + \mi [ \mathbf{J}_s \mathbf{R}_t^{-1} \mathbf{J}_s^T (\mathbf{Q}_t + \1)^{-1}, \mathbf{J} ] \\
            &\quad - (\1 - \mi \mathbf{J}) ( (\mathbf{Q}_t + \1)^{-1} + (\mathbf{Q}_t + \1)^{-1} \mathbf{J}_s \mathbf{R}_t^{-1} \mathbf{J}_s^T (\mathbf{Q}_t + \1)^{-1} ) (\1 + \mi \mathbf{J}) \\
            &= \coth(\beta/4) (\1 + \me^{t \mathbf{Z}^T}\me^{t \mathbf{Z}} ) ( \1 - \me^{t \mathbf{Z}^T}\me^{t \mathbf{Z}} )^{-1} + \mi (1 + \me^{\beta/2}) [\mathbf{J}, (\1 - \me^{t \mathbf{Z}^T}\me^{t \mathbf{Z}})^{-1}] \\
            &\quad + \sinh(\beta/2) (\1 - \mi \mathbf{J}) (\1 - \me^{t \mathbf{Z}^T}\me^{t \mathbf{Z}})^{-1} (\1 + \mi \mathbf{J}) \\
            &= \coth(\beta/4) (\1 + \me^{t \mathbf{Z}^T}\me^{t \mathbf{Z}} ) ( \1 - \me^{t \mathbf{Z}^T}\me^{t \mathbf{Z}} )^{-1} + \mi (1 + \me^{\beta/2}) [\mathbf{J}, (\1 - \me^{t \mathbf{Z}^T}\me^{t \mathbf{Z}})^{-1}] \\
            &\quad + \frac{1}{2} \sinh(\beta/2) (\1 - \mi \mathbf{J}) (\1 + \me^{t \mathbf{Z}^T}\me^{t \mathbf{Z}}) (\1 - \me^{t \mathbf{Z}^T}\me^{t \mathbf{Z}})^{-1} (\1 + \mi \mathbf{J}) \\
            &= \cosh^2(\beta/4) \coth(\beta/4) (\1 + \me^{t \mathbf{Z}^T}\me^{t \mathbf{Z}} ) ( \1 - \me^{t \mathbf{Z}^T}\me^{t \mathbf{Z}} )^{-1} \\
            &\quad + \frac{1}{2} \sinh(\beta/2) \mathbf{J} (\1 + \me^{t \mathbf{Z}^T}\me^{t \mathbf{Z}} ) ( \1 - \me^{t \mathbf{Z}^T}\me^{t \mathbf{Z}} )^{-1} \mathbf{J} \\
            &\quad + \mi (1 + \cosh(\beta/2)) [\mathbf{J}, (1 - \me^{t \mathbf{Z}^T} \me^{t \mathbf{Z}})^{-1}] \\
            &= \frac{1}{2} \sinh(\beta/2) (\mathbf{\Sigma} + \mi \mathbf{J}) (\1 + \me^{t \mathbf{Z}^T}\me^{t \mathbf{Z}} ) ( \1 - \me^{t \mathbf{Z}^T}\me^{t \mathbf{Z}} )^{-1} (\mathbf{\Sigma} - \mi \mathbf{J}).
        \end{align*}
        Applying \eqref{eq-Theta-identity-2}, we obtain

        \begin{align*}
            \mathbf{A}_t &= \frac{1}{2} \sinh(\beta/2) \mathbf{\Theta}_{\frac{1}{2} - s} (\mathbf{\Sigma} + \mi \mathbf{J}) \\
            &\quad \cdot (\1 + \me^{t \mathbf{Z}^T} \me^{t \mathbf{Z}}) ( \1 - \me^{t \mathbf{Z}^T} \me^{t \mathbf{Z}} )^{-1} (\mathbf{\Sigma} - \mi \mathbf{J}) \mathbf{\Theta}_{\frac{1}{2} - s}^T \\
            &=\coth(\beta/4)\mathbf{\Theta}_{1-s}  (\1 + \me^{t \mathbf{Z}^T} \me^{t \mathbf{Z}}) ( \1 - \me^{t \mathbf{Z}^T} \me^{t \mathbf{Z}} )^{-1}\mathbf{\Theta}_{1-s}^T.
        \end{align*}

        Consequently,
        \begin{equation*}
            \mathbf{A}_t^{-1} =\tanh(\beta/4)\mathbf{\Theta}_{1-s} ( \1 - \me^{t \mathbf{Z}^T} \me^{t \mathbf{Z}} ) (\1 + \me^{t \mathbf{Z}^T} \me^{t \mathbf{Z}} )^{-1}\mathbf{\Theta}_{1-s}^T,
        \end{equation*}
        from which \eqref{eq-simplified-A-t-inv-B-t-T} and \eqref{eq-simplified-schur-complement} follow.
    \end{proof}

    \begin{lemma} \label{lemma-simplification-characteristic-function-sandwiched-K}
        We have
        \begin{align}
            \mathbf{C}_t - \mathbf{B}_t \mathbf{D}_t^{-1} \mathbf{B}_t^T &= \mathbf{\Sigma}, \label{eq-charc-sandwiched-K-identity-1} \\
            \mathbf{\Sigma} - \csch^2(\beta/4)\mathbf{\Theta}_s \mathbf{D}_t^{-1}\mathbf{\Theta}_s^T &= \mathbf{Q}_t, \label{eq-charc-sandwiched-K-identity-2} \\
            \mathbf{B}_t \mathbf{D}_t^{-1}\mathbf{\Theta}_s^T &= -\mathbf{\Theta}_{1-s} \me^{t \mathbf{Z}}. \label{eq-charc-sandwiched-K-identity-3}
        \end{align}
    \end{lemma}
    \begin{proof}
        Recall from \eqref{eq-definition-D-t} that $ \mathbf{D}_t = \mathbf{A}_t + \mathbf{\Sigma} $. Hence
        \begin{align*}
            \mathbf{D}_t &=\coth(\beta/4)\mathbf{\Theta}_{1-s} ( \1 + \me^{t \mathbf{Z}^T}  \me^{t \mathbf{Z}} ) (\1 - \me^{t \mathbf{Z}^T}  \me^{t \mathbf{Z}})^{-1}\mathbf{\Theta}_{1-s}^T \\
            &\quad+\coth(\beta/4)\mathbf\Theta_{1-s}\mathbf\Theta_{1-s}^T \\
            &=2\coth(\beta/4)\mathbf{\Theta}_{1-s} (\1 - \me^{t \mathbf{Z}^T} \me^{t \mathbf{Z}})^{-1}\mathbf{\Theta}_{1-s}^T.
        \end{align*}
        Consequently,
        \begin{equation*}
            \sqrt{\det \mathbf{D}_t} = 2 \coth(\beta/4) \left( \det (\1 - \me^{t \mathbf{Z}^T} \me^{t \mathbf{Z}} ) \right)^{-1/2}
        \end{equation*}
        and
        \begin{equation} \label{eq-inverse-D-t}
            \mathbf{D}_t^{-1} = \frac12\tanh(\beta/4)\mathbf{\Theta}_{1-s} (\1 - \me^{t \mathbf{Z}^T} \me^{t \mathbf{Z}} )\mathbf{\Theta}_{1-s}^T.
        \end{equation}
        Using \eqref{eq-inverse-D-t}, we compute
        \begin{align*}
            &\quad \mathbf{C}_t -\mathbf{B}_t \mathbf{D}_t^{-1} \mathbf{B}_t^T \\
            &=\coth(\beta/4)\mathbf{\Theta}_{1-s} (\1 + \me^{t \mathbf{Z}} \me^{t \mathbf{Z}^T}) ( \1 - \me^{t \mathbf{Z}} \me^{t \mathbf{Z}^T} )^{-1}\mathbf{\Theta}_{1-s}^T \\
            &\quad+2\coth(\beta/4)\mathbf{\Theta}_{1-s} ( \1 - \me^{-t \mathbf{Z}^T} \me^{- t \mathbf{Z}} )^{-1}\mathbf{\Theta}_{1-s}^T \\
            &=\coth(\beta/4)\mathbf\Theta_{1-s}\mathbf\Theta_{1-s}^T
            =\coth(\beta/4)\1=\mathbf\Sigma,
        \end{align*}
        which gives \eqref{eq-charc-sandwiched-K-identity-1}.

        The remaining identities follow directly from \eqref{eq-inverse-D-t}:
        \begin{equation*}
            \mathbf{\Sigma} -\csch^2(\beta/4)\mathbf{\Theta}_s \mathbf{D}_t^{-1}\mathbf{\Theta}_s^T = \coth(\beta/4) \1 - \csch(\beta/2) (\1 - \me^{t \mathbf{Z}^T} \me^{t \mathbf{Z}}) = \mathbf{Q}_t,
        \end{equation*}
        and
        \begin{equation*}
            \mathbf{B}_t \mathbf{D}_t^{-1}\mathbf{\Theta}_{s}^T
            =2\coth(\beta/4)\cdot\frac12\tanh(\beta/4)
            \mathbf{\Theta}_{1-s}( \me^{t \mathbf{Z}^T } - \me^{- t \mathbf{Z}} )^{-1}
            ( \1 - \me^{t \mathbf{Z}^T} \me^{t \mathbf{Z}} )
            =-\mathbf{\Theta}_{1-s}\me^{t\mathbf Z}.
        \end{equation*}
    \end{proof}

    \begin{lemma}
    \label{lemma-estimate-bounds-of-M-t}
        For every $s\in[0,1]$ and $t>0$,
        \begin{equation} \label{eq-estimate-bounds-of-M-t}
            \norm{\me^{t\mathbf\Theta_{1-s}\mathbf Z\mathbf\Theta_{1-s}^T}}<1.
        \end{equation}
    \end{lemma}
    \begin{proof}
        Definition~\eqref{definition-Theta-s} and the identity $\mathbf S=\coth(\beta/2)\1$ give
        \begin{equation} \label{eq-Theta-square-Lyapunov-combination}
        \mathbf\Theta_r^2=\frac12\sech(\beta/2)
        \Big(\sinh((1-r)\beta)(\mathbf S-\mi\mathbf J)
        +\sinh(r\beta)(\mathbf S+\mi\mathbf J)\Big).
        \end{equation}
        Taking $r=1-s$ and using \eqref{eq-Theta-basic-identities}, we obtain
        \begin{equation} \label{eq-dissipativity-transformed-Z}
        \begin{aligned}
        &\qquad \left(\mathbf\Theta_{1-s}\mathbf Z\mathbf\Theta_{1-s}^T\right)^*
          +\mathbf\Theta_{1-s}\mathbf Z\mathbf\Theta_{1-s}^T = \mathbf\Theta_{1-s}^T
        \left(\mathbf Z^T\mathbf\Theta_{1-s}^2
        +\mathbf\Theta_{1-s}^2\mathbf Z\right)
        \mathbf\Theta_{1-s}^T\\
        &\quad=-\frac12\sech(\beta/2)\mathbf\Theta_{1-s}^T
        \Big(\sinh(s\beta)\overline{\mathbf C_Z}
        +\sinh((1-s)\beta)\mathbf C_Z\Big)
        \mathbf\Theta_{1-s}^T\leq0.
        \end{aligned}
        \end{equation}
        The last equality follows from \eqref{eq-Lyapunov-plus-iJ}--\eqref{eq-Lyapunov-minus-iJ}. Hence $\mathbf\Theta_{1-s}\mathbf Z\mathbf\Theta_{1-s}^T$ is dissipative. It is also stable because $\mathbf\Theta_{1-s}^T=\mathbf\Theta_{1-s}^{-1}$, so it is similar to the stable matrix $\mathbf Z$.

        Equivalently, $-\mathbf\Theta_{1-s}\mathbf Z\mathbf\Theta_{1-s}^T$ is an accretive hypocoercive matrix in the terminology of \cite[Definition~5]{achleitner2023hypocoercivity}. \cite[Proposition~12(a)]{achleitner2023hypocoercivity} gives constants $c,\varepsilon>0$ and an odd positive integer $m$ such that
        \[
            \norm{\me^{t\mathbf\Theta_{1-s}\mathbf Z\mathbf\Theta_{1-s}^T}}
            =1-ct^m+O(t^{m+1}),
            \qquad 0\leq t<\varepsilon.
        \]

        Thus the norm is strictly smaller than one for every sufficiently small $t>0$. Dissipativity also gives a contraction for every $t\geq0$. For an arbitrary $t>0$, choose $n\in\N$ so large that $t/n$ lies in the interval of strict contraction. The semigroup property then yields
        \[
            \norm{\me^{t\mathbf\Theta_{1-s}\mathbf Z\mathbf\Theta_{1-s}^T}}
            \leq
            \norm{\me^{(t/n)\mathbf\Theta_{1-s}\mathbf Z\mathbf\Theta_{1-s}^T}}^n<1,
        \]
        which proves \eqref{eq-estimate-bounds-of-M-t}.
    \end{proof}

    \begin{lemma} \label{lemma-inverse-real-part-bound}
        If $X$ is a complex matrix with $\norm{X}<1$, then
        \begin{equation} \label{eq-real-part-inverse-bound}
            \Re(\1+X)^{-1}\geq\frac{1}{1+\norm{X}}\,\1.
        \end{equation}
    \end{lemma}
    \begin{proof}
        The identities
        \[
            2(\1+X)^{-1}=\1+(\1-X)(\1+X)^{-1}
        \]
        and
        \begin{equation} \label{eq-real-part-Cayley-identity}
            \Re\!\left((\1-X)(\1+X)^{-1}\right)
            =(\1+X^*)^{-1}(\1-X^*X)(\1+X)^{-1}
        \end{equation}
        imply
        \[
            \Re\!\left((\1-X)(\1+X)^{-1}\right)
            \geq\frac{1-\norm{X}}{1+\norm{X}}\,\1.
        \]
        Substitution into the first identity gives \eqref{eq-real-part-inverse-bound}.
    \end{proof}

    \begin{lemma} \label{lemma-strict-positivity-real-parts-R-t}
        For every $t>0$, $\Re\mathbf R_t>0$.
    \end{lemma}
    \begin{proof}
        From \eqref{eq-definition-R-t}, \eqref{eq-definition-J-s}, and \eqref{eq-Theta-basic-identities}, direct simplification gives
        \begin{align}
        \mathbf R_t={}&(\me^{\beta/2}-1)^2\csch^2(\beta/4)\sinh(\beta/2) \nonumber \\ 
        &\quad \cdot \mathbf{J}^T
        \left(
        \left(\me^{\beta/2}\1
        +\me^{t\mathbf\Theta_{1-s}\mathbf Z^T\mathbf\Theta_{1-s}^T}
         \me^{t\mathbf\Theta_{1-s}\mathbf Z\mathbf\Theta_{1-s}^T}\right)^{-1}
         - \frac{1}{1+\me^{\beta/2}}\1
        \right)\mathbf J. \label{eq-R-t-contraction-form}
        \end{align}
        For completeness, this follows by writing
        \[
            \mathbf Q_t+\1=\csch(\beta/2)
            \left(\me^{\beta/2}\1+\me^{t\mathbf Z^T}\me^{t\mathbf Z}\right)
        \]
        and conjugating the matrix in parentheses with $\mathbf\Theta_{1-s}$; the factors $\mathbf\Theta_{1-s}$ cancel against those in
        $\mathbf J_s=(\me^{\beta/2}-1)\csch(\beta/4)\mi\mathbf J\mathbf\Theta_s$.

        By \eqref{eq-estimate-bounds-of-M-t}, the norm of the product of exponentials in \eqref{eq-R-t-contraction-form} is strictly smaller than one. Applying \eqref{eq-real-part-inverse-bound} with
        \[
        X=\me^{-\beta/2}
        \me^{t\mathbf\Theta_{1-s}\mathbf Z^T\mathbf\Theta_{1-s}^T}
        \me^{t\mathbf\Theta_{1-s}\mathbf Z\mathbf\Theta_{1-s}^T}
        \]
        gives
        \begin{align*}
        &\Re\Big(\me^{\beta/2}\1
        +\me^{t\mathbf\Theta_{1-s}\mathbf Z^T\mathbf\Theta_{1-s}^T}
         \me^{t\mathbf\Theta_{1-s}\mathbf Z\mathbf\Theta_{1-s}^T}\Big)^{-1}\\
        &\qquad\geq
        \frac{1}{\me^{\beta/2}+
        \norm{\me^{t\mathbf\Theta_{1-s}\mathbf Z^T\mathbf\Theta_{1-s}^T}
        \me^{t\mathbf\Theta_{1-s}\mathbf Z\mathbf\Theta_{1-s}^T}}}\,\1
        >\frac{1}{1+\me^{\beta/2}}\1.
        \end{align*}
        The scalar prefactor in \eqref{eq-R-t-contraction-form} is positive and $\mathbf J$ is invertible. Hence $\Re\mathbf R_t>0$.
    \end{proof}

    We now establish the strict positivity of the real part of the block matrix $ \mathbf{G}_t $ governing the Gaussian kernel.

    \begin{lemma} \label{lemma-strict-positivity-real-parts-G-t}
        For every $t>0$, $\Re\mathbf G_t>0$. Consequently,
        \[
            \Re\mathbf A_t>0,\qquad \Re\mathbf C_t>0,
            \qquad
            \Re\mathbf C_t-(\Re\mathbf B_t)(\Re\mathbf A_t)^{-1}(\Re\mathbf B_t)^T>0.
        \]
    \end{lemma}
    \begin{proof}

        Recall from \eqref{eq-definition-G-t} that
        \[
            \mathbf G_t=
            \begin{pmatrix}
                \mathbf C_t&-\mathbf B_t\\
                -\mathbf B_t^T&\mathbf A_t
            \end{pmatrix},
        \]
        and that \eqref{eq-G-t-Cayley-transform} expresses this matrix as a positive scalar multiple of a Cayley transform. We first verify that the block matrix in that formula is a strict contraction, as required in Lemma~\ref{lemma-inverse-real-part-bound}.

        For every complex matrix $X$,
        \begin{equation} \label{eq-block-matrix-norm}
            \left\|\begin{pmatrix}0&X\\X^T&0\end{pmatrix}\right\|=\norm{X}.
        \end{equation}
        Indeed, after taking the adjoint and multiplying, the square of the norm on the left is the norm of
        \[
            \operatorname{diag}\bigl(\overline X X^T,X^*X\bigr).
        \]
        Both diagonal blocks have norm $\norm{X}^2$, which proves \eqref{eq-block-matrix-norm}.

        Since $\mathbf\Theta_{1-s}^T=\mathbf\Theta_{1-s}^{-1}$, similarity invariance of the exponential gives
        \[
            \mathbf\Theta_{1-s}\me^{t\mathbf Z}\mathbf\Theta_{1-s}^T
            =\me^{t\mathbf\Theta_{1-s}\mathbf Z\mathbf\Theta_{1-s}^T}.
        \]
        Taking $X=\mathbf\Theta_{1-s}\me^{t\mathbf Z}\mathbf\Theta_{1-s}^T$ in \eqref{eq-block-matrix-norm}, and then using Lemma~\ref{lemma-estimate-bounds-of-M-t}, gives
        \begin{equation*}
        \left\|\begin{pmatrix}
                0&\mathbf\Theta_{1-s}\me^{t\mathbf Z} \mathbf\Theta_{1-s}^T \\ 
                \mathbf\Theta_{1-s}\me^{t\mathbf Z^T}\mathbf\Theta_{1-s}^T&0
            \end{pmatrix}\right\| =\left\|\me^{t\mathbf\Theta_{1-s}\mathbf Z\mathbf\Theta_{1-s}^T}\right\|<1.
        \end{equation*}

        Now let $F$ be any complex matrix with $\norm{F}<1$. Applying the identity \eqref{eq-real-part-Cayley-identity} from Lemma~\ref{lemma-inverse-real-part-bound} with $X=-F$ gives
        \begin{equation} \label{eq-real-part-G-t-Cayley}
            \Re\!\left((\1+F)(\1-F)^{-1}\right) =(\1-F^*)^{-1}(\1-F^*F)(\1-F)^{-1} >0.
        \end{equation}
        Indeed, $\1-F$ is invertible and
        \[
            \1-F^*F\geq(1-\norm{F}^2)\1>0,
        \]
        so the right-hand side of \eqref{eq-real-part-G-t-Cayley} is an invertible congruence of a strictly positive matrix. Applying \eqref{eq-real-part-G-t-Cayley} to the block matrix in \eqref{eq-G-t-Cayley-transform}, and using $\coth(\beta/4)>0$, proves $\Re\mathbf G_t>0$.

        Finally, the block form \eqref{eq-definition-G-t} shows that $\Re\mathbf A_t$ and $\Re\mathbf C_t$ are principal submatrices of $\Re\mathbf G_t$; hence both are strictly positive. Applying the Schur-complement criterion to the lower-right block $\Re\mathbf A_t$ gives
        \[
            \Re\mathbf C_t-(\Re\mathbf B_t)(\Re\mathbf A_t)^{-1}(\Re\mathbf B_t)^T>0,
        \]
        as claimed.

    \end{proof}

    \end{appendices}

    \printbibliography

    % \textbf{Data Availability}: This manuscript does not have associated data.

    % \bigbreak

    % \textbf{Conflict of interest}: The authors declare that there is no conflict of interest.

\end{document}